\documentclass [12 pt, twoside]{article}
\usepackage{amsmath,amsthm,amssymb,amsfonts,amscd,makeidx,amsopn, setspace}
\usepackage[all]{xy}
\makeindex

\theoremstyle{plain}

\newtheorem{thm}{Theorem}[section]

\newtheorem{pro}[thm]{Proposition}

\newtheorem{lem}[thm]{Lemma}

\newtheorem{cla}[thm]{Claim}
\newtheorem{cor}[thm]{Corollary}
\newtheorem{con}[thm]{Conjecture}

\newtheorem{step}{Step}
\newtheorem{case}{Case}

\theoremstyle{definition}

\newtheorem{dfn}[thm]{Definition}

\newtheorem{nt}[thm]{Notation}

\newtheorem{dfnlm}[thm]{Definition-Lemma}
\newtheorem{rem}[thm]{Remark}

\newtheorem{exa}[thm]{Example}
\newtheorem{setup}[thm]{Set up}

\theoremstyle{remark}

\newcommand{\Z}{\mathbb{Z}}
\newcommand{\N}{\mathbb{N}}
\newcommand{\C}{\mathbb{C}}

\newcommand{\Q}{\mathbb{Q}}
\newcommand{\PS}{\mathbb{P}}
\newcommand{\F}{\mathbb{F}}

\DeclareMathOperator{\Pic}{Pic} \DeclareMathOperator{\Ker}{Ker}
\DeclareMathOperator{\rk}{rk} \DeclareMathOperator{\W}{Weil}
\DeclareMathOperator{\Coker}{Coker}

\DeclareMathOperator{\codim}{Codim}
\DeclareMathOperator{\Exc}{Exc} \DeclareMathOperator{\Spec}{Spec}
 
\DeclareMathOperator{\proj}{Proj} \DeclareMathOperator{\Bl}{Bl}
\DeclareMathOperator{\Bs}{Bs}
\DeclareMathOperator{\Def}{Def}
\DeclareMathOperator{\Sing}{Sing}
\DeclareMathOperator{\Hom}{Hom}
\DeclareMathOperator{\Ext}{Ext}
\DeclareMathOperator{\HHom}{\mathcal{H}om}
\DeclareMathOperator{\EExt}{\mathcal{E}xt}

\DeclareMathOperator{\PProj}{\underline{Proj}}
\DeclareMathOperator{\im}{Im}
\DeclareMathOperator{\projL}{proj Lim}
\DeclareMathOperator{\Hilb}{\underline{Hilb}}
\DeclareMathOperator{\nklt}{nklt}
\DeclareMathOperator{\Ann}{Ann}
\DeclareMathOperator{\Tor}{Tor}
\DeclareMathOperator{\Red}{Red}

\begin{document}
\bibliographystyle{alpha}

\title{The topology of terminal quartic $3$-folds}
\author{Anne-Sophie KALOGHIROS}
\date{}

\maketitle
\tableofcontents

\section{Introduction}
\label{sec:introduction}
\subsection{Main results}
\label{sec:main-results}

Consider a possibly singular complex projective variety $Y$ of dimension
$3$. The $3$-fold $Y$ is said to be $\Q$-factorial if an integral
multiple of every Weil divisor on $Y$ is a Cartier
divisor. $\Q$-factoriality is a subtle property, which is not local in
the analytic topology. For instance, an ordinary
double point is not locally analytically $\Q$-factorial, yet a nodal quartic
$3$-fold $Y_4^3 \subset \PS^4$ is known to be $\Q$-factorial if it has
less than $8$ ordinary double points. The topology of a quartic
$3$-fold $Y$ with isolated singularities is well understood when $Y$ is
$\Q$-factorial. In this case, the Grothendieck-Lefschetz theorem 
states that every Weil (or Cartier)
divisor on $Y$ is the restriction of a Weil (or Cartier) divisor defined on
$\PS^4$. However, if the $3$-fold $Y$ is not $\Q$-factorial, very
little is known about its topology. In this thesis, I study the topology of
some mildly singular quartic $3$-folds in $\PS^4$.

Let $Y=Y_4^3 \subset \PS^4$ be a quartic $3$-fold. I assume that $Y$ has
\emph{terminal} singularities. Notice in addition that $Y$ is a
locally complete intersection: it is therefore \emph{Gorenstein}. These
assumptions ensure, on the one hand, that the canonical class of $Y$
is well defined as a Cartier divisor and, on the other, that pull-backs of the
canonical class of $Y$ to any resolution are well behaved. More
precisely, let $\widetilde{Y} \to Y$ be a resolution of
singularities. If $Y$ is Gorenstein, the canonical sheaf $\omega_Y$ is
locally free. Since $k(Y)=k(\widetilde{Y})$, a local generating
section $s$ of $\omega_Y$ near a
singular point can be regarded as a rational differential on
$\widetilde{Y}$. If $Y$ has terminal
singularities, the section $s$ is regular as a rational
differential on $\widetilde{Y}$ and vanishes along the exceptional
locus of $\widetilde{Y} \to Y$. If $Y$ has terminal Gorenstein
singularities, $Y$ is $\Q$-factorial if and only if it is factorial;
that is, if its local rings are unique factorisation domains. In the
case of quartic $3$-folds, this means that $Y$ is $\Q$-factorial if
and only if $\dim H^2(Y, \Z)=\dim
H_4(Y, \Z)$. The Grothendieck-Lefschetz theorem on Picard groups states
that $\Pic Y \simeq \Pic \PS^4 \simeq \Z [\mathcal{O}_Y(1)]$. Yet, no
such result holds for the
group of Weil divisors $H_4(Y, \Z)$. Let $\sigma(Y)= b_4(Y)- b^2(Y)=
b_4(Y)-1$ be the
\emph{defect} of $Y$. The defect of $Y$ measures how far $Y$ is from being
$\Q$-factorial or, in other words, to what extent Poincar\'e duality 
fails on $Y$.

Such quartic $3$-folds are a special case of Fano $3$-folds with
terminal Gorenstein singularities and Picard rank $1$. A variety $X$
is Fano
 if its anticanonical sheaf $\omega_X^{-1}$ is ample. The defect of a
 Fano $3$-fold $Y$ with terminal Gorenstein
singularities can be defined as above. $\Q$-factorial Fano $3$-folds
with terminal singularities and
Picard rank $1$ play a crucial role in Mori
theory: they arise as one of the possible end products of the Minimal
Model Program for non-singular varieties. Moreover, it is known that terminal
Gorenstein $\Q$-factorial Fano $3$-folds are
deformations of non-singular ones \cite{MR1944132}. Non-singular Fano
$3$-folds of Picard rank $1$
have been classified \cite{MR463151,MR503430}. Yet, very
little is understood about the topology of non $\Q$-factorial Fano $3$-folds.

The defect of some very simple quartic $3$-folds with isolated
singularities is already non-zero. For instance, if $X$ is a
sufficiently general quartic $3$-fold containing a plane, then it has
$9$ ordinary double points and is not $\Q$-factorial. Similarly, a
general determinantal quartic $3$-fold has $20$ ordinary double points
and is not $\Q$-factorial. Finally, consider the linear system $\Sigma$ of
quartics spanned by the monomials $\{x_0^4, x_1^4, (x_4^2x_3+x_2^3)x_0,
x_3^3x_1, x_4^2x_1^2\}$ on $\PS^4$. A general quartic $X\in \Sigma$ is
not $\Q$-factorial and yet it has a unique singular
point $P=(0:0:0:0:1)$, which is a c$A_1$ point \cite{MR2091681}.

The study of quartics with ordinary double points suggests that any
 bound on the defect of terminal quartic $3$-folds should be at least $15$. Indeed, a quartic
 $3$-fold $Y_4^3 \subset \PS^4$ with no worse than ordinary double
 points is known to have at most $45$ nodes
 \cite{MR848512,MR712934}. Up to projective equivalence, there is a
 unique quartic with $45$ nodes: the Burkhardt quartic, given by the equation
\[
 \{x_0^4 -x_0 (x_1^3+x_2^3+x_3^3+x_4^3)+ 3x_1x_2x_3x_4 = 0 \}. 
\]
It is easy to show that the defect of the Burkhardt quartic is at least $15$
(Section~\ref{sec:mixed-hodge-theory}). 
It is a natural
question to ask how many topological types of quartic $3$-folds with
terminal Gorenstein singularities there are. I provide a partial
answer to this question by bounding and studying the defect of such
quartic $3$-folds.   

The first main result of this work is a bound on the defect of
terminal quartic $3$-folds.
\begin{thm}[Main Theorem $1$]
\label{thm:18}
Let $Y_4^3 \subset \PS^4$ be a quartic $3$-fold with terminal singularities. The defect of
$Y$ is at most:
\begin{enumerate}
\item $8$ if $Y$ does not contain a plane or a quadric,
\item $11$ if $Y$ contains a quadric but no plane,
\item $15$ if $Y$ contains a plane.
\end{enumerate}
\end{thm}
As Example~\ref{exa:1} shows, this bound is attained by the Burkhardt
quartic, a quartic with no worse than ordinary double
points.

The second main result is a geometric ``motivation'' of the global
topological property of
$\Q$-factoriality. If $Y$ is
not $\Q$-factorial, by definition, $Y$ contains a special surface. More
precisely, $Y$ contains a Weil non $\Q$-Cartier divisor. I show
that this special surface belongs to a finite list. In particular,
its degree is bounded. 
\begin{thm}[Main Theorem $2$]
 \label{thm:19}
Let $Y_4^3 \subset \PS^4$ be a terminal Gorenstein quartic $3$-fold. Then one of the following holds:
\begin{enumerate} 
\item $Y$ is $\Q$-factorial.
\item $Y$ contains a plane $\PS^2$.
\item $Y$ contains an irreducible reduced quadric $Q$. 
\item  $Y$ contains an anticanonically embedded del Pezzo surface of
  degree $4$. 
\item $Y$ has a structure of Conic Bundle over $\PS^2$, $\F_0$
  or $\F_2$.
\item $Y$ contains a rational scroll $E \to C$  over a curve $C$ whose
  genus and degree appear in the table on page \pageref{table}.
\end{enumerate}  
\end{thm}
 
Typical examples of non $\Q$-factorial varieties are $3$-folds
that contain planes or quadrics, which are not Cartier. I show that, in
accordance with geometric intuition, a
quartic $3$-fold has to contain a surface of low degree if it is not
$\Q$-factorial.

\subsection{Outline of the thesis}
\label{sec:outline-thesis}

I have divided this thesis in six Sections.

Section~\ref{sec:mixed-hodge-theory} reviews relevant material and
results from Mixed Hodge Theory. The
notion of defect of hypersurfaces was first introduced by Clemens
\cite{MR690465} in an attempt, based on Deligne's Mixed Hodge Theory
\cite{MR0498552}, to extend Griffiths' results on
the Hodge theory of hypersurfaces \cite{MR0260733} to mildly singular
varieties. The
traditional approach to the determination of the defect of
hypersurfaces, or of Fano $3$-folds, has
focused on understanding their mixed Hodge structures, as in
\cite{MR1047124,MR1358982}. Such an approach heavily
relies on explicit computations of cohomology groups of specific
varieties: it is impractical to determine a sharp bound on the
defect of quartic $3$-folds with terminal singularities.

In Section~\ref{sec:categ-weak-fano}, I define the category of
weak$\ast$ Fano $3$-folds and
I show that the Minimal Model Program can be run in this category. For any
Fano $3$-fold $Y$ with terminal Gorenstein singularities, there is
a small $\Q$-factorialisation $X \to
Y$ \cite{MR924674}. If $Y$ does not contain a plane, $X$ belongs to
the category of weak$\ast$ Fano $3$-folds. Bounding the defect of $Y$
is equivalent to determining the maximal number of divisorial
contractions involved in a Minimal Model Program on $X$. 

In Section~\ref{sec:bound-defect-some-1}, I bound the defect of any
terminal Gorenstein Fano $3$-fold that does not contain a
plane. I show
that if the anticanonical model $Y$ of a weak$\ast$ Fano $3$-fold $X$ does not
contain a quadric, this property is preserved when running a Minimal
Model Program on $X$. This enables me to improve the bound on the
defect when the $3$-fold $Y$ does not contain a quadric.

In Section~\ref{sec:fano-3-folds}, I study quartic
$3$-folds containing a plane. The $3$-fold obtained by
blowing up $Y$ along this plane has a natural structure of del Pezzo
fibration of degree $3$. Using Corti's results \cite{MR1426888} to
relate the defect of $Y$ to the number of reducible fibres of this
cubic fibration, I bound the number of reducible fibres: this completes the
proof of Theorem~\ref{thm:18}.

Section~\ref{sec:deformation-theory} recalls several results on the
deformation theory of so-called generalised Fano $3$-folds. These
results were mainly obtained by
Namikawa, Koll\'ar and Mori \cite{MR1489117,MR1149195}. Namikawa
defines \emph{generalised Fano $3$-folds} and shows that any generalised
Fano $3$-fold is a one parameter flat deformation of a non-singular
one. Any small (partial)
$\Q$-factorialisation of a terminal Gorenstein Fano $3$-fold and
terminal Gorenstein Fano $3$-folds themselves are generalised Fano
$3$-folds. Moreover, the degree and the Picard rank are constant in
this deformation. The key observation is that each step of the Minimal Model Program
on weak$\ast$ Fano
$3$-folds may be deformed, in a suitable sense, to an extremal contraction of a generalised
Fano $3$-fold of Picard rank $2$. 

In Section~\ref{sec:takeuchi-game}, I extend several ideas expressed by
Takeuchi in
\cite{MR1022045}. Any divisorial
contraction of the Minimal Model Program on a small
$\Q$-factorialisation $X$ of $Y$  induces an extremal divisorial
contraction on $Z$,
a Picard rank $2$ partial $\Q$-factorialisation of $Y$. I study
extremal divisorial contractions with Cartier exceptional divisor on
partial $\Q$-factorialisations of terminal Gorenstein Fano $3$-folds.
Any such divisorial contraction can then be deformed
to a divisorial contraction on a
non-singular generalised Fano $3$-fold of Picard rank $2$. Following
Takeuchi's approach, I then write systems of Diophantine
equations associated to
each extremal contraction of the Minimal Model Program on $X$. These
equations translate numerically the properties of the extremal contraction. 
Such systems have very few
solutions. The explicit study of the system associated to the first
divisorial contraction on $X$ yields a finite number of
solutions. Each solution exhibits a possible surface that has to
be contained in the quartic $3$-fold $Y$. This establishes
Theorem~\ref{thm:19} and gives a geometric ``motivation'' of
$\Q$-factoriality in the case of quartic $3$-folds.  

Similar methods may be used to bound the defect of
any terminal Gorenstein Fano $3$-fold with Picard rank $1$ and degree
greater than $4$. I have
bounded the defect of the Fanos with Picard rank $1$ that
contain no plane. An explicit study of terminal Gorenstein Fano
$3$-folds that contain a plane is also possible.

\begin{con}
  \begin{enumerate}
  \item
 The defect of a non $\Q$-factorial terminal Gorenstein
 $Y_{2,3}\subset \PS^5$ with Picard rank $1$ is at most $8$,  
  \item
 The defect of a non $\Q$-factorial terminal Gorenstein
 $Y_{2,2,2}\subset \PS^6$  with Picard rank $1$ is at most $8$,
\item
 The defect of a Picard rank $1$ non $\Q$-factorial terminal Gorenstein
 Fano $3$-fold of genus $g \geq 6$ is at most $\big[\frac{12-g}{2}\big]+5$.
\end{enumerate}
\end{con}

However, bounding the defect of the double cover of a sextic
divisor in $\PS^3$ is likely to be
more complicated. The study of such double sextics, whether they contain
planes or surfaces whose normalisations are planes,
would be significantly more difficult than the quartic case. I
conjecture that the defect of a terminal Gorenstein double sextic with
Picard rank $1$ is at most $18$.

The methods I develop in this thesis can also be applied to terminal
Gorenstein Fano $3$-folds of larger Picard rank.

It would be natural to consider wider classes of singularities such as
cDV singularities, canonical  Gorenstein or even non-canonical
singularities. It is likely that the analysis required will be much
finer. Throughout this work, I made extensive use of the fact that
crepant
contractions are small and assumed that the anticanonical
class was Cartier. 
I do not believe that the methods I have developped may be extended
directly to these more general settings.

It is known that non-singular quartics are non-rational, while the Burkhardt
quartic is rational. The defect
could be used, in some cases, to determine whether a quartic $3$-fold
$Y_4^3\subset \PS^4$ is rational or not. Non
$\Q$-factoriality does not necessarily imply rationality. Indeed, a
general quartic containing a plane is non-rational.  Yet, if $Y$ does not contain a
plane, the defect of $Y$ can only be high when the end
product of the Minimal Model Program on its
$\Q$-factorialisation is rational. This could provide a rationality criterion for some
values of the defect. Examples of rational non
$\Q$-factorial quartics can be obtained by running
a Minimal Model Program in reverse, using Takeuchi's numerical
constraints.

I should like to
thank the range of mathematicians with whom I have had many useful
conversations during my PhD. I should like to thank Miles Reid, Burt
Totaro and Ivan Smith for their comments on preliminary
versions of this thesis. I should especially like to
thank my supervisor Alessio Corti, for his guidance, his encouragement
and for teaching me some very beautiful Mathematics.

\section{Mixed Hodge Theory}
\label{sec:mixed-hodge-theory}

Clemens introduced the notion of defect of hypersurfaces
\cite{MR690465} in an attempt to generalise Griffiths' results
\cite{MR0260733} on the Hodge theory of non-singular hypersurfaces to
hypersurfaces with ordinary double points. 

I recall in this section his definition of
the defect of hypersurfaces and the subsequent generalisations of this
notion to the framework of terminal Gorenstein Fano $3$-folds. The
defect is expressed in terms of local cohomology groups. 

I show how the defect of a terminal Gorenstein Fano $3$-fold $Y$ is related to the mixed Hodge structure on the
cohomology $H^{\bullet}(Y)$ and, in particular, to the weight
filtration. The defect depends not only on the Hodge filtration, but
also on the weight filtration: it is a
global topological invariant. It cannot be determined by local
analytic methods. 

\subsection[Defect of a Fano $3$-fold]{The notion of defect of a Fano $3$-fold}
\label{sec:notion-defect-fano-1}

Let $Y$ be a non $\Q$-factorial Fano $3$-fold with terminal Gorenstein
singularities.

Kawamata shows \cite{MR924674} that a
$3$-fold $Y$ with terminal
singularities has a (not necessarily unique) small projective
$\Q$-factorialisation. He proves:
 
\begin{pro}[$\Q$-factorialisation]
\label{pro:4}
Let $Y$ be an algebraic threefold with only terminal
singularities. Then, there is a birational morphism $f \colon X \to Y $
such that $X$ is terminal and $\Q$-factorial, $f$ is an
isomorphism in codimension $1$ and $f$ is projective.
\end{pro}

Let $X$ be a small crepant projective $\Q$-factorialisation of $Y$. 

\begin{dfn}
\label{dfn:8}
The \emph{defect} of $Y$ is
\begin{eqnarray*} 
\sigma(Y)=\rk(\W(Y)/ \Pic(Y))= \dim H_4(Y)-\dim H^2(Y)
\\
=\dim \Pic(X) - \dim \Pic(Y). 
\end{eqnarray*}
\end{dfn}
\begin{rem}
  \label{rem:21}
In the case of hypersurfaces in projective space with isolated
singularities, such as a quartic $3$-fold $Y_4^3 =\{f=0\} \subset
\PS^4$ with isolated singularities, it is known that the Euler characteristic depends only on the
degree and on local invariants related to the singularities. More
precisely, the Euler characteristic may in theory be computed once one
understands the semi-simple part of the
monodromy operator, i.e. the spectrum of $f$ in the neighbourhood of
the singularities. The defect, or rather $b_4(Y)$, is however a more subtle invariant:
 it also depends on the position of the singularities.  
\end{rem}
\begin{lem} \cite{MR1358982}
  \label{lem:10}
Let $Y$ be a terminal Gorenstein Fano $3$-fold. Denote by $\Sigma$ the
singular locus $\Sigma=\Sing(Y)$ and by $U$ its complement $U =
Y\smallsetminus \Sigma$. The defect of $Y$ is
\[
\sigma(Y)= \dim_{\C}\Coker[H^3(U, \C) \to H^4_{\Sigma}(Y, \C)].
\]  
\end{lem}
\begin{proof}
Denote by $\Sigma= \{P_1, \cdots , P_n\}$ the singular set of $Y$.
Consider the local cohomology exact sequence of mixed Hodge structures:
\begin{equation}
  \label{eq:8}
 H^3(U) \to H^{4}_{\Sigma}(Y) \to H^4(Y) \to H^4(U)\to  H^5_{\Sigma}(Y)
\end{equation}

Let $\{U_i\}_{1\leq i\leq n}$ be mutually disjoint open neighbourhoods
of the singular points $P_i$. By excision, $H^5_{\Sigma}(Y)= \oplus
H^5_{\Sigma}(U_i)=\oplus  H^5_{\{P_i\}}(Y)$. In 
Section~\ref{sec:mixed-hodge-struct-3}, I show that $H^5_{\{p_i\}}(Y)=0$.  
Hence, $\sigma (Y)= b_4(Y)-b_4(U)$. 

By definition
$U$ is non-singular, hence $H^4(U, \C)$ is dual to the cohomology group
with compact support $H^2_c(U)$. The group $H^2_c(U)$ is isomorphic to
the relative cohomology $H^2(Y, \Sigma) \simeq H^2(Y)$ because the
singularities of $Y$ are isolated.  
\end{proof}

The value of $b_n(Y)$ for a hypersurface $Y$ of $\PS^n$ is
known to be related to the dimension of some linear systems of
homogeneous polynomials that vanish at the singular points of $Y$
\cite{MR690465}. Dimca makes these relations explicit \cite{MR1047124} and
shows that the defect depends on the mixed Hodge
structure of the local cohomology groups at the singular points.

I recall his results in the case of $Y \subset \PS^4$ a quartic
$3$-fold with isolated singularities.

Let $\Sigma = \{ P_1, \ldots , P_n\}$ be the singular set of $Y$ and
$U=Y \smallsetminus \Sigma$. Let $i$ (resp. $j$) be the inclusion $Y
\to \PS^4$ (resp. $U \to \PS^4$). The \emph{primitive} part of
the cohomology of $Y$ (resp. $U$) is $H^{\bullet}_{\text{prim}}(Y)= \Coker
[i^{\ast}\colon H^{\bullet}(\PS^4)\to  H^{\bullet}(Y)]$ (resp. $H^{\bullet}_{\text{prim}}(U)=
\Coker(j^{\ast})$). The natural inclusion $k \colon U \to Y$ induces a
morphism $k_{\text{prim}}\colon H^{\bullet}_{\text{prim}}(Y)\to
H^{\bullet}_{\text{prim}}(U)$, and carries isomorphically the
non-primitive part of $H^{\bullet}(Y)$ into the non-primitive part of
$H^{\bullet}(U)$ except in the top dimension.

The Poincar\'e residue map \[R \colon
H^k((\PS^4 \smallsetminus \Sigma) \smallsetminus (Y \smallsetminus \Sigma)) \to
H^{k-1} (Y\smallsetminus \Sigma)\] defines a $(-1, -1)$ isomorphism of Hodge structures
from $ H^k(\PS^4 \smallsetminus Y)=H^k((\PS^4 \smallsetminus \Sigma)
\smallsetminus (Y \smallsetminus \Sigma))$ to the primitive part of
the cohomology of $U= Y\smallsetminus \Sigma$, $H^{k-1}_{\text{prim}}(U)$. From the exact
sequence 
\begin{eqnarray*}
H^{3}(U) \stackrel{\theta}\to H^{4}_{\Sigma}(Y) \to H^{4}(Y)  
\end{eqnarray*}
and the Poincar\'e isomorphism $R\colon H^{4}(\PS^4 \smallsetminus \Sigma)
\simeq H^{3}_{\text{prim}}(U)$, we deduce the exact sequence:
\begin{eqnarray}
\label{eq:7}
H^4(\PS^4 \smallsetminus Y) \stackrel{\delta} \to H^4_{\Sigma}(Y)\to H^4_{\text{prim}}(Y)\to 0,  
\end{eqnarray}
where $\delta= \theta\circ R$. 
The local cohomology $H^{\bullet}_{\Sigma}(Y)$ has a natural mixed
Hodge structure inherited from that of the cohomology of the pair $(Y,
Y-\Sigma)$ (Section~\ref{sec:append-mixed-hodge}). Let
$\{U_i\}_{1\leq i\leq n}$ be mutually disjoint open neighbourhoods of the
singular points $P_i$; by excision,
$H^k_{\Sigma}(Y)= \oplus H^k_{\{P_i\}}(U_i)$ and $H^4_{\Sigma}(V)$ is
computable.

Denote by $F$ the Hodge filtration. Recall
that the Hodge filtration in the case of hypersurfaces with terminal
singularities coincides with  the filtration by the order of the pole
on $\PS^4-Y$ \cite{MR0260733,MR1072821}.  
The polar filtration on the De Rham complex $A^{\bullet}=H^0(\PS^4-Y,
\Omega^{\bullet}_{\PS^4-Y})$ is given by:
\[
F^t A^j=\{ \omega \in A^j \mid \omega \mbox{ has a pole along $Y$ of
  order at most $j-t$}\},
\] 
for $j-t\geq 0$ and $F^t A^j=0$ for $j-t<0$.
Recall Griffiths' explicit description of $A^4$. Consider the
differential $4$-form $\Omega \in \Omega^4_{\PS^4}$
\[
\Omega= \sum_{0\leq i \leq 4}(-1)^i x_i dx_0 \wedge\cdots \wedge
\widehat{dx_i} \wedge \cdots \wedge dx_4. 
\]
If $Y$ is the hypersurface $\{f=0\} \subset \PS^4$, any element $\omega \in A^4$ may be written
\[
\omega= \frac{P \Omega}{f^t}
\]
where $P$ is a homogeneous polynomial of degree $4t-5$. If $f$ does
not divide $P$, the order of the pole of $\omega $ along $Y$ is
precisely $t$. 

Let $t$
be the highest natural number such that $F^t H^4_\Sigma(Y)=
H^4_\Sigma(Y)$ and call $\delta^t$ the linear map:
\[
\delta^t \colon S_{4(4-t)-5} \to H^4_{\Sigma}(Y).
\]
This map is, strictly speaking, a composition of the above $\delta$ with
the natural map $S_{4(4-t)-5} \simeq F^t A^4 \to F^t H^4(Y \smallsetminus \Sigma)$.

\begin{dfnlm}\cite{MR1047124}
\label{dfnlm:2}
 Let $Y \subset \PS^4 $ be a quartic $3$-fold. The defect of $Y$ is:
\[
\sigma(Y)= \dim H^4_{\Sigma}(Y) - \codim \Ker \delta^t= \dim H^4_{\text{prim}}(Y).   
 \]
\end{dfnlm}
 \begin{rem}
\label{rem:46}
Note that if the quartic hypersurface $Y\subset \PS^4$ has no worse
than ordinary double points, the computations in
Section~\ref{sec:mixed-hodge-struct-3} show that $t=2$ and \eqref{eq:7} gives a lower
bound for the defect of $Y$ \cite{MR1194180}: 
\[
\sigma(Y)\geq N -\frac{b_3(V)}{2}= N-30,
\]   
where $N$ is the number of double points and $V$ is a non-singular
quartic hypersurface in $\PS^4$. 
  \end{rem}
\begin{proof}
This expression agrees in the case of quartics with the expression of
the defect given in Lemma~\ref{lem:10}. Indeed, the map
$\delta^t$ factors through $H^3(U, \C)$. More precisely, Dimca shows
\cite{MR1047124} that the map $\delta$ factors as:
\[
\delta \colon H^4(\PS^4 \smallsetminus Y)\stackrel{R}\to H^3(U) \to H^4_{\Sigma}(Y), 
\]
where the second map is the connecting morphism in the sequence of
local cohomology.
\end{proof}
\begin{rem}
  \label{rem:19}
The defect of a quartic hypersurface with terminal
singularities is expressed in
terms of residuation from $A^{\bullet}= H^0(\PS^4\smallsetminus Y,
\Omega^{\bullet})$. This extends Griffiths' results on the Hodge theory
of non-singular hypersurfaces.
\end{rem}
The equivalent definitions of the defect of a quartic hypersurface
with terminal singularities show that the only singularities affecting
the defect of a quartic $3$-fold are \emph{essential
  singularities}. 
\begin{dfn}
  \label{dfn:14}
A singularity $P_i\in Y$ is \emph{essential} if the local cohomology group
$H^{4}_{\{P_i\}}(Y)$ is not trivial.
\end{dfn} 
The local cohomology groups $H^{k}_{\{P_i\}}(Y)$ associated to
isolated singularities are computed in Section~\ref{sec:mixed-hodge-struct-3}. 
\begin{rem}\cite{MR1047124}
  \label{rem:40}
Given a homogeneous polynomial $h \in S_{4(4-t)-5}$, $\delta^t(h)=0$ means
that $h$ satisfies certain linear conditions $\mathcal{C}$ at the essential singular points.
Denote by $D$ the linear system in $S_{4(4-t)-5}$ defined by these
conditions. The defect of $Y$ is the difference between the number of
linear conditions $\mathcal{C}$ and the codimension of $D$ in
$S_{4(4-t)-5}$. The defect does not depend only on the linear system $D$,
but also on the number of conditions $\mathcal{C}$ used to define it. The conditions
$\mathcal{C}$ are independent if and only if the defect is zero.
\end{rem}

The expression of the defect given in \ref{dfnlm:2} allows direct computations
of the defect in some special cases, such as nodal hypersurfaces. 

\begin{lem}\cite{MR1047124}
 \label{lem:16}
Let $Y \subset \PS^4$ be a quartic $3$-fold whose only essential
singularities are nodes. Then the defect of $Y$ is given by:
\begin{gather*}
\sigma(Y)=h^{2,2}(H^4_{\text{prim}}(Y))\\
= \dim H^3_{\Sigma}(Y)- \codim \{h \in
S_3=S_{4(4-2)-5},\text{$h(P)=0$ for any node of $Y$} \}\\
=\sharp\text{\{nodes\}}- \sharp \text{\{conditions imposed on
  cubics by vanishing at the nodes of $Y$\}}.
\end{gather*} 
\end{lem}
\subsection[Mixed Hodge structure and defect]{Mixed Hodge structure of $H^{\bullet}(Y)$ and defect}
\label{sec:mixed-hodge-struct}

I have explained in Section~\ref{sec:notion-defect-fano-1} that the
defect is related to the cohomology groups $H^3(Y)$. I adapt ideas
formulated by
Namikawa and Steenbrink for Calabi-Yau $3$-folds \cite{MR1358982} to give a geometric
description of the mixed Hodge structure on $H^3(Y)$. 

The cohomology group $H^3(Y)$ is naturally endowed with a mixed Hodge structure
supported in weights $2$ and $3$. There are at least two interpretations of this
mixed Hodge structure. 

First, the weight $3$ part of
$H^3(Y)$ is isomorphic to the pure Hodge structure
$H^3(\widetilde{Y})$ on a good resolution $\widetilde{Y}$ of $Y$. 
 
The weight $2$
part of $H^3(Y)$ is non-zero if there exists a non trivial relation
among the $H^2(E_i)$ in $H^2(\widetilde{Y})$ (where $E_i$
denotes the exceptional divisor above the singular point $P_i$). This
corresponds to the geometric intuition that the defect of a $3$-fold
is not only related to the analytic local type of the singularities,
but also to their relative position. 

Second, if there is a smoothing $\mathcal{Y}\to \Delta$ of $Y$,
examining the relationship between the mixed Hodge structure of $Y$ and
the Hodge structure of a fibre $\mathcal{Y}_t$ for $t\neq 0$ expresses the defect in terms of the
vanishing cohomology. The weight $3$ part of $H^3(Y)$ is
isomorphic to $H^3(\mathcal{Y}_t)$, while its weight $2$ part depends
on the cohomology of vanishing cycles.     

\subsection{Defect and cohomology of a good resolution of $Y$}
\label{sec:defect-cohom-good}
Let $Y$ be a terminal Gorenstein Fano $3$-fold. Denote by
$\Sigma=\{P_1, \cdots, P_n\}$ its singular set. Let
$f\colon\widetilde{Y}\to Y$ be a good resolution of $Y$ and $E= \sum
E_i$ be the simple normal crossings exceptional divisor of $f$. For each $i$,
$\{E_i^j\}_j$ are the irreducible components of $E_i$.

\begin{dfnlm}\cite{MR1358982}
Let $Y_i$ be a contractible Stein open neighbourhood of $P_i$ and
$\sigma^{\text{an}}(P_i)= \W(Y_i)/\Pic(Y_i)$. The \emph{analytic local defect} at
$P_i$ is:
\[
\sigma^{\text{an}}(P_i)= \dim( H^1(Y_i, \mathcal{O}_{Y_i}^{\ast})/\bigoplus
\Z[E_{i}^j])=\dim H^3_{\{P_i\}}(Y)= \dim H^4_{\{P_i\}}(Y).
\]  
\end{dfnlm}
\begin{proof}
The expressions for $H^3_{\{P_i\}}(Y)$ and $ H^4_{\{P_i\}}(Y)$ are
determined in Section~\ref{sec:mixed-hodge-struct-3}.
\end{proof}
\begin{rem}
  \label{rem:26}
In particular, if $P_i$ is an ordinary double point, $\sigma^{\text{an}}(P_i)=1$.
\end{rem}
\begin{pro}\cite{MR1358982}
\label{pro:3}
Let $Y$ be a terminal Gorenstein Fano $3$-fold and $\Sigma=\{P_1,
\cdots, P_n\}$ be the singular set of $Y$. The weight filtration
on $ H^3(Y)$ has the following description:
\begin{eqnarray*}
  \label{eq:1}
Gr_k^W H^3(Y)= (0) \quad \mbox{for}\quad k \neq 2,3 \\
\dim W_2 H^3(Y) = \sum \sigma^{\text{an}}(P_i)-\sigma(Y)  
\end{eqnarray*}
\end{pro}
\begin{proof}
As is explained in Section~\ref{sec:mixed-hodge-struct-1}, the long
exact sequence \eqref{eq:2} 
\[
 \cdots H^l(Y, \Q) \to H^l(\widetilde{Y}, \Q)\oplus H^l(\Sigma, \Q) \to
H^l(E, \Q) \to  H^{l+1}(Y, \Q) \to \cdots
\]
 is compatible with the weight
filtration. 
In particular, the following sequences are exact:
\begin{eqnarray*}
0 \to Gr^W_3 H^3(Y) \to H^3(\tilde{Y}) \to 0 \\
0 \to Gr^W_2 H^2(Y) \to H^2(\tilde{Y}) \to \oplus Gr^W_2 H^2(E_i)\to
Gr^W_2 H^3(Y) \to 0 \\
0 \to \oplus Gr^W_1 H^2(E_i)\to Gr^W_1 H^3(Y) \to 0 \\
0 \to \oplus Gr^W_0 H^2(E_i)\to Gr^W_0 H^3(Y) \to 0 
\end{eqnarray*}

Section~\ref{sec:mixed-hodge-struct-2} shows that
$H^2(E)=Gr^W_2 H^2(E)$ is of pure weight $2$. 
The exact sequence of local cohomology groups associated to $\Sigma$ is:
\[
 \ldots \to H^2(Y)\to H^2(U) \stackrel{\alpha}\to  H^3_{\Sigma}(Y) \to
 H^3(Y) \to H^3(U)
 \to \ldots.
\]
The sequence
\[
0 \to H^2_{E}(\widetilde{Y}) \to H^2(E)\to H^3_{\Sigma}(Y)\to 0
\]
is exact and compatible with the weight filtration (Section~\ref{sec:mixed-hodge-struct-3}). 

The local
cohomology group $H^3_{\Sigma}(Y)$ is purely of weight $2$, while
$H^3(U)$ is purely of weight $3$. Moreover,
\[
\im(\alpha)= \im [H^2(\tilde{Y}) \to H^2(U) \to  H^3_{\Sigma}(Y)],
\]
and  $W_2 H^3(Y)= H^3_{\Sigma}(Y)/ \im(\alpha)= \Coker
[H^2(\tilde{Y})\to H^3_{\Sigma}(Y)]$.
The $3$-fold $Y$ is a Fano $3$-fold and has rational singularities. By
Kawamata-Viehweg vanishing (Theorem~\ref{thm:11}) and the Leray
spectral sequence,
$H^2(\tilde{Y}, \mathcal{O}_{\tilde{Y}})=(0)$ and $H^2(\tilde{Y})=
H^1(\tilde{Y}, \mathcal{O}_{\tilde{Y}}^{\ast}) \otimes \C$.
The result then follows from the isomorphism:
\[
\W(Y)/\Pic(Y)\simeq \im \big[H^1(\tilde{Y},
\mathcal{O}_{\widetilde{Y}}^{\ast}) \to \bigoplus H^1(Y_i,
\mathcal{O}_{Y_i}^{\ast})/\oplus_j \Z[E_i^j]  \big].
\]
This isomorphism is a direct application of the Kawamata-Viehweg
vanishing theorem and
of the Leray spectral sequence on $Y_i$ \cite{MR924674, MR1149195}.
\end{proof} 

\begin{rem}
\label{rem:22}
Notice that this gives the following formula for the defect of a
terminal Gorenstein Fano $3$-fold:
\[
\sigma(Y)= b_3(\widetilde{Y})-b_3(Y)+ \sum \sigma^{\text{an}}(P_i).
\] 
\end{rem}
\begin{exa}[Case of a nodal $3$-fold]
Assume that the singularities of $Y$ are ordinary double points
(nodes). Denote by $\mu \colon \widetilde{Y}\to Y$ the  blow up of
$Y$ at the nodes $P_i$ and by $\nu \colon\widehat{Y}\to Y$ a small
resolution of the nodes of $Y$. 
  \begin{enumerate}
  \item The $\mu$-exceptional locus is a disjoint union of
non-singular quadrics $Q_i=\PS^1 \times \PS^1$. In that case, for each $i$,
$\sigma^{\text{an}}(P_i)=1$. Hence, if the defect is strictly less than the number
of nodes, $W_2 H^3(Y)\neq (0)$. 
The sequence associated to the
resolution $\widetilde{Y}\to Y$ (Section~\ref{sec:mixed-hodge-struct-1},\eqref{eq:2}) reads
\[ 
\cdots \to H^2(\widetilde{Y})\to H^2(E)=\oplus H^2(Q_i) \to W_2 H^3(Y) \to 0 
.\] 
Hence, if $W_2 H^3(Y) \neq (0)$, there is a non trivial relation between the $Q_i$-s in $H^2(\tilde{Y},
\C)$. In this case, the nodes fail to impose independent linear
relations on the vanishing of cubics (Lemma~\ref{lem:16}).  
\item Note for future reference that the second Betti numbers of $Y,
  \widehat{Y}$ and $\widetilde{Y}$ satisfy the following relations:
  \begin{equation}
    \label{eq:47}
b_2(\widehat{Y})= b_2(Y)+ \sigma(Y) \quad \mbox{and}\quad b_2(\widetilde{Y})=
b_2(\widehat{Y})+ N,     
  \end{equation}
where $N$ is the number of nodes.
\end{enumerate}
\end{exa}
\subsection{Defect and cohomology of a smoothing of $Y$}
\label{sec:defect-cohom-smooth}

Clemens compared the cohomology of a nodal $3$-fold $Y$ to the cohomology
of a smoothing of $Y$ \cite{MR690465}. This lead him to introduce the
notion of defect of $Y$: whereas the presence of nodes does not affect
the second integral cohomology, the fourth homology is altered and
Poincar\'e duality fails on $Y$. I recall how the defect relates to
the cohomology of a smoothing of $Y$.
 
\begin{dfn}
\label{dfn:16}  
Let $Y$ be a terminal Gorenstein $3$-fold.
A \emph{smoothing} of $Y$ is a proper flat map $f \colon
\mathcal{Y}\to \Delta$ from an analytic space $\mathcal{Y}$ to a
$1$-dimensional complex disc $\Delta$ such that $f^{-1}(0)=Y$ and
$f^{-1}(t)= \mathcal{Y}_t$ is a non-singular $3$-fold for all $t\in
\Delta\smallsetminus \{0\}$.
\end{dfn}

In Section~\ref{sec:deformation-theory}, I recall several results on the
deformation theory of Fano $3$-folds. In particular, Namikawa shows
\cite{MR1489117} that if $Y$ is a nodal Fano $3$-fold
with Picard rank $1$,
there exists a smoothing $f\colon \mathcal{Y}\to \Delta$ of $Y$.
The fibre $\mathcal{Y}_t$ is a non-singular Fano
$3$-fold for $t \neq 0$. If $Y$ is a terminal Gorenstein Fano $3$-fold
with Picard rank $1$, there exists a $1$-parameter proper flat
deformation $f\colon \mathcal{Y}\to \Delta$ such that
$\mathcal{Y}_{t_0}$ is non-singular for some $t_0\in \Delta$ and
$\mathcal{Y}_t$ is a terminal Gorenstein Fano $3$-fold for all $t\in
\Delta$. In such a deformation, the plurigenera are
constant. The relation between the cohomology groups of
$Y$ and those of a fibre $\mathcal{Y}_t$ provides further information
on the defect of $Y$.

\begin{lem}\cite{MR0239612,MR0485870}
  \label{lem:27}
Let $Y$ be a normal projective $3$-fold with terminal Gorenstein
singularities. Suppose that $Y$ has a smoothing $f \colon
\mathcal{Y}\to \Delta$. Then, $\Pic Y \simeq \Pic \mathcal{Y}_t$.
\end{lem}
\begin{proof}[Sketch Proof]  

Let $\{P_1, \cdots , P_n\}$ be the singular set of $Y$ and denote by
$(U_i,P_i)$ a small Stein open neighbourhood of the singular point
$P_i$. Terminal Gorenstein $3$-fold singularities are isolated
hypersurface singularities. Let $f_i\colon \mathcal{U}_i\to \Delta$ be
a $1$-parameter flat deformation of $(U_i,p_i)$. Let $B_i$ be a small
ball of $\mathcal{U}_i$ with centre $P_i$ and radius $\epsilon
>0$. For $\mathcal{\eta}$ sufficiently small, all fibres
$(\mathcal{U}_i)_t$ for $|t|<\eta$ intersect $B_i$ transversally and
$(\mathcal{U}_i)_t\cap B_i= B_{i,t}$ is diffeomorphic to the Milnor
fibre of $P_i$ for $t\neq 0$. If the singularity at $P_i$ is given in local
analytic coordinates by $\{f(x,y,z,t)=0\}$, for $f$ a polynomial with an
isolated critical point at the origin, the Milnor fibre is $\{f(x,
y,z,t)=1\}$. The Milnor fibre is homotopic to a bouquet of $3$-spheres
\cite{MR0239612}; its cohomology is supported in degrees $0$ and $3$. There is
a homeomorphorphism between $U_i\smallsetminus \{P_i\}$ and
$(\mathcal{U}_i)_t\smallsetminus B_{i,t}$, hence the
exact sequence of relative cohomology shows that for $|t|<\eta$:
\begin{eqnarray*}
H^i(U_i,\Z) \simeq H^i((\mathcal{U}_i)_t, \Z) \quad \mbox{ for } i \neq
3,4, \\
0 \to H^3(U_i)\to H^3((\mathcal{U}_i)_t)\to H^3(B_{i,t})\to
H^4(U_i)\to H^4((\mathcal{U}_i)_t)\to 0.
\end{eqnarray*}
There is a homeomorphism between $Y\smallsetminus \{P_1,\cdots,
P_n\}$ and $\mathcal{Y}_t \smallsetminus \bigoplus_{1\leq i\leq n}
B_{i_t}$, so that for $|t|$ sufficiently small:
\begin{eqnarray*}
H^i(Y,\Z) \simeq H^i(\mathcal{Y}_t, \Z) \quad \mbox{ for } i \neq
3,4, \\
0 \to H^3(Y)\to H^3(\mathcal{Y}_t)\to \bigoplus_{1\leq i\leq n} H^3(B_{i,t})\to
H^4(Y)\to H^4(\mathcal{Y}_t)\to 0.
\end{eqnarray*}
\end{proof}
\begin{lem}
\label{lem:7}
If $\mathcal{Y}\to \Delta$ is a smoothing of $Y$, the defect of $Y$ satisfies:
\[
\sigma(Y)= b_3(Y)-b_3(\mathcal{Y}_t)+ \sum h^3(B_{i,t}).
\]  
\end{lem} 
\begin{proof}
The sequence of mixed Hodge structures
\[ 
0 \to H^3(Y) \to H^3(\mathcal{Y}_t) \to \bigoplus H^3(B_i)\to
H^4(Y)\to H^4(\mathcal{Y}_t) \to 0
\]
is exact. It follows that 
\[
 b_3(Y)+\sum h^3(B_{i,t})-b_3(\mathcal{Y}_t)=b_4(Y)-b_4(\mathcal{Y}_t).
\]
By Poincar\'e duality on the non-singular $3$-fold $\mathcal{Y}_t$,
$b_4(\mathcal{Y}_t)=b_2(\mathcal{Y}_t)$. 
As the Picard rank is constant in the smoothing, $b_2(\mathcal{Y}_t)=
b_2(Y)$ and the result follows.
\end{proof}
\begin{rem}
\label{rem:15} Notice that $h^3(B_{i,t})=\dim H^3_{\{P_i\}}(Y)$. 
\end{rem}
\begin{rem}
\label{rem:10}
For an ordinary double point, $B_{i,t}$ is, by construction, a
$3$-sphere $\{x^2+y^2+z^2+t^2=|t|\}$, so that $h^3(B_{i,t})=1$.
In particular, if $Y$ is a quartic with no worse than ordinary double
points and if there exists a smoothing of $Y$, the following holds:
\[
\sigma(Y)= b_3(Y)+ N - 60, 
\] 
where $N$ is the number of nodes. Indeed, the degree of
a Fano $3$-fold is constant in a flat family, and the third Betti
number of a non-singular quartic $3$-fold is $60$.
\end{rem}

\begin{rem}
\label{rem:24}
As the third Betti numbers of non-singular Fano $3$-folds with Picard
rank $1$ are known \cite{MR1668579}, Lemma~\ref{lem:7} provides
a bound on the number of essential singular points, once a bound on
the defect of terminal Gorenstein Fano $3$-folds is known.
\end{rem}  

\begin{rem}
\label{rem:27}
This analysis could be carried out in the more general
case of terminal or canonical singularities. However, in such cases, one
would have to study the complex of vanishing cohomology (or
vanishing cycles). This complex is not, in general, supported in degrees $0$
and $3$.  
\end{rem}

I have presented several definitions and formulae related to the
defect of Gorenstein Fano $3$-folds, based on their Mixed Hodge
Theory. Lengthy computations are necessary in order to explicitly determine the defect
of a particular terminal Gorenstein quartic or Fano $3$-fold. These
expressions rely on the analytic local type of the
singularities and are too unwieldy to yield a bound on the
defect of terminal quartic $3$-folds.

\subsection{Mixed Hodge structures on
  $H^{\bullet}(Y)$ and $H^{\bullet}_{\Sigma}(Y)$}  
\label{sec:append-mixed-hodge}

In the previous sections, I used several results on the Mixed Hodge
Theory of various cohomology groups associated to a Gorenstein
terminal Fano $3$-fold $Y$. In this Section, I indicate how these results
are obtained. I do not provide complete proofs, but give an overview
of the relevant aspects of the theory.  

Let $X$ be a complex algebraic variety. Deligne endows the cohomology
groups $H^{\bullet}(X, \Z)$ with a functorial mixed Hodge structure
\cite{MR0498552}. His formalism extends the Hodge theory of non-singular
projective complex varieties.

Let $X$ be a complex algebraic variety. Hironaka's
resolution theorem states that there
exists a resolution of singularities $f \colon \widetilde{X} \to
X$. The map $f$ is proper and birational and $\widetilde{X}$ is
non-singular. The exceptional divisor of $f$ has simple normal
crossings.
Deligne introduces the method of simplicial cohomological descent,
which uses resolution of singularities to
define a mixed Hodge structure on the cohomology of $X$.
Whereas earlier applications of resolution theorems had implications
on the hypercohomology of complexes of sheaves, cohomological descent
yields results at the level of these complexes.

I use the techniques of cubic hyper-resolutions,
developed in \cite{MR972983} to determine the mixed Hodge structure
on a terminal Gorenstein Fano $3$-fold $Y$.

\subsubsection{Mixed Hodge structure on $H^{\bullet}(Y)$}
\label{sec:mixed-hodge-struct-1}

Let $Y$ be a terminal Gorenstein Fano $3$-fold and let
$\Sigma=\{P_1, \ldots, P_n \} $ be the singular set of $Y$.

There is a good resolution $f \colon \widetilde{Y}
\to Y$ of $Y$. The morphism $f$ is proper and birational, and its
exceptional set is a divisor $E= \sum E_i$ with simple
normal crossings. As above, $E_i$ denotes $f^{-1}(P_i)$ and $\{E_i^j\}_j$
are the irreducible components of $E_i$.

The diagram 
\[
\xymatrix{
\sum E_i \ar[r]^{\tilde{i}} \ar[d] & \tilde{Y} \ar[d]^{f}\\
\Sigma \ar[r]^{i}& Y
}
\]
induces the long exact sequence in cohomology:
\begin{equation}
\label{eq:2}
 \cdots\to H^l(Y, \Q) \to H^l(\widetilde{Y}, \Q)\oplus H^l(\Sigma, \Q) \to
H^l(E, \Q) \to  H^{l+1}(Y, \Q) \to \cdots
\end{equation}
This long exact sequence is compatible with the Hodge and weight
filtrations.

The above diagram is said to be of cohomological descent. It induces a resolution of the complex $\Q_{Y}$
\[ 
0 \to
\Q_Y 
\to \mathbb{R}f_{\ast} \Q_{\widetilde{Y}} \oplus
\mathbb{R}i_{\ast} \Q_{\Sigma} \to
\mathbb{R}(f\tilde{i})_{\ast} \Q_{\sum E_i} \to 0.
\] 
In the present case, this resolution reads :
\[ 0 \to
\Q_Y 
\to \mathbb{R}f_{\ast} \Q_{\widetilde{Y}} \bigoplus \oplus \Q_{\{P_i\}}
 \to
\mathbb{R}(f\tilde{i})_{\ast} \Q_{\sum E_i}\to 0.
\]  
The cohomology groups appearing in \eqref{eq:2} satisfy the following
set of properties. These properties follow either from definitions or from standard results in Mixed Hodge
Theory \cite{MR0498552,MR972983} because $\widetilde{Y}$ is
non-singular and
projective and $Y$ is projective with isolated singularities.    
\begin{enumerate}
\item $H^l(\Sigma, \Q)=0$ for $l \neq 0$ and $H^0(\Sigma, \Q)= \oplus \Q_{\{P_i\}}$.
\item For all $l$, $H^l(\widetilde{Y}, \Q)$ is a pure Hodge structure of
  weight $l$.
\item For all $j>l$, $Gr_j^W H^l(Y, \Q)=0 $. 
\item For all $l$, $H^l(E, \Q)= \oplus H^l(E_i, \Q)$.
\end{enumerate}

Proposition~\ref{pro:3} relates the defect of $Y$ to the weight
filtration $H^3(Y)$. The cohomology of $\widetilde{Y}$ has a pure
Hodge structure, hence \eqref{eq:2} shows that understanding 
the weight filtration on the cohomology of $E$ suffices to determine
the defect of $Y$.

In the next subsection, I recall how the weight filtration on
$H^{\bullet}(E)$ is obtained.

\subsubsection{Mixed Hodge structure on $H^{\bullet}(E, \C)$}
\label{sec:mixed-hodge-struct-2}

By definition, the
irreducible components $E_i^j$ of $E_i$ are non-singular and intersect
tranversally. 

The mixed Hodge structure of $H^{\bullet}(E_i)$
is determined by a Mayer-Vietoris cubic hyper-resolution. From now
onwards, I drop the index $i$ and write $E$ for $E_i$.

Let $E^j$, $1\leq j \leq r$ be the irreducible components of
$E$. For any $p \leq r$, let $E_{p-1}$ denote the disjoint union of all
$p$-fold intersections of components of $E$.
\[
E_{p-1}= \amalg_{1 \leq j_0 <\cdots <j_p \leq r} (E^{j_0} \cap
\cdots \cap E^{j_r}).
\]
The projective varieties $E_{p-1}$ and the maps 
\[
d_k^p\colon E_p \to E_{p-1}
\]
 for $k=1,\cdots, p$ induced by the inclusions
\[
E^{j_1}\cap \ldots \cap E^{j_{p+1}} \to E^{j_1}\cap\ldots \cap
E_{j_{k-1}}\cap E_{j_{k+1}}\cap \ldots \cap E^{j_{p+1}},  
\]
and the natural augmentation maps $a_p\colon E_p \to E$ define a
strict simplicial variety $E_{\bullet}$, that is a simplicial resolution of $E$ \cite{MR972983}. 

The Mayer-Vietoris hyper-resolution of $E$ is fairly simple. Very few
terms are non-zero as $E$ is a simple normal crossings divisor in a $3$-fold. In particular, the
$p$-fold intersections  of components of $E$ are empty for $p>3$, so
that $E_p$ is non-empty only for $0 \leq p \leq 2 $.   

Each $E_{p}$ is non-singular and projective, hence there exists a
cohomological Hodge complex $((\Q_{E_{p}}, W),(\Omega_{E_{p}},
F,W))$ on $E_p$. These Hodge complexes yield a cohomological mixed Hodge complex $K= ((K_{\Q},
W), (K_{\C}, W, F))$ on $E$ and, in particular, the following resolution of
$\Q_{E}$:
\[
0\to \Q_{E} \to a_{0 \ast}\Q_{E_{0}}\to a_{1 \ast}\Q_{E_{1}}\to  a_{2
  \ast}\Q_{E_{2}} \to 0.
\]
I omit results related to the Hodge filtration. The weight
spectral sequence reads:
\[
{}_wE_1^{p,q}=H^q(E_{p}, \C).
\]   
The spectral sequence $({}_wE_r, d_r)$ abuts to $H^{p+q}(E)$. It
degenerates at ${}_wE_2$ and $d_1^{p,q}$ is induced
by $\sum_{k=0}^{p-1} (-1)^{p+k}(d_{k}^{p-1})^{\ast}$ \cite{MR972983}. 
In particular, 
\[
Gr_p^W H^{p+q}(Y)={}_wE_2^{p,q}= H({}_wE_1^{p-1, q} \to {}_wE_1^{p,
  q}\to {}_wE_1^{p+1, q}).
\]
In the present case, $E_{p}=\emptyset$ for $p>2$, $E_{0}$ has dimension $2$, $E_{1}$
has dimension $1$ and  $E_{2}$ is a set of points.
The only non-zero terms in the spectral sequence associated to the
weight filtration are ${}_wE_1^{0,q}$ for $0 \leq q \leq 4$,
${}_wE_1^{1,r}$ for $0 \leq r \leq 2$ and ${}_wE_1^{2,0}$.

The definition of the map $d_1^{\ast}: H^{\bullet}(E_0)\to
H^{\bullet}(E_1)$ implies that $H^2(E)$ is of pure weight $2$. 

I give the result of the computation of $H^4(E)$, which is used in the next
subsection:
\[
H^4(E, \C)=Gr_4^W H^4(E)\simeq \oplus \C[E^j].
\]  

\subsubsection{Mixed Hodge structure on the local cohomology $H^{\bullet}_{\Sigma}(X, \C)$}
\label{sec:mixed-hodge-struct-3}

I use the same notation as above. The cohomology of the
pair $H^{\bullet}(Y, Y \smallsetminus \Sigma)$, and therefore the local cohomology
$H^{\bullet}_{\Sigma}(Y)$, carries a mixed Hodge structure. These
mixed Hodge structures are
independent of the choice of resolution \cite{MR713277}. 

I recall some results on local cohomology groups associated to
isolated singularities.

\begin{thm}[Goresky-MacPherson vanishing]\cite{MR713277}
\label{thm:20}
Let $Y_i$ be a contractible Stein neighbourhood of the isolated singular
point $P_i$ and $\widetilde{Y_i}$ a good resolution with exceptional
divisor $E_i$. The restriction map 
\[
H^{k}(\widetilde{Y_i}, \Q) \to  H^k(\widetilde{Y_i} \smallsetminus E_i,
\Q)\simeq H^k(Y_i \smallsetminus \{P_i\}, \Q)
\]  
is surjective for $k \leq 2$ and the zero map for $ k \geq 3$.
\end{thm}
\begin{pro}[Fujiki duality]\cite{MR602463}
\label{pro:9}
 There is a duality isomorphism:  
\begin{equation}
\label{eq:9}
H^k_E(\widetilde{X}) \simeq \Hom(H^{6-k}(E), \Q(-3)),  
\end{equation}
where $\Q(-3)$ is the mixed Hodge structure on $\C$ with rational
lattice $\frac{1}{(2\pi i)^3}\Q$ purely of weight $(3,3)$.
\end{pro}

\begin{pro}\cite{MR713277}
  \label{pro:10}
The local cohomology groups $H^k_{\{P_i\}}(Y_i)$ fit in the following exact
sequences:
\begin{eqnarray}
\label{eq:34}
0 \to H^k_{E_i}(\widetilde{Y_i}) \to H^k(E_i) \to
H^{k+1}_{\{P_i\}}(Y_i)\to 0 \quad \text{for} \quad k<3
\\
0 \to H^3_{E_i}(\widetilde{Y_i}) \to H^3(E_i) \to 0\\
0 \to H^{k}_{\{P_i\}}(Y_i) \to H^k_{E_i}(\widetilde{Y_i}) \to
H^k(E_i)\to 0 \quad \text{for}\quad k>3.
\end{eqnarray} 
\end{pro}
\begin{rem}
\label{rem:42}
The exact sequences in Proposition~\ref{pro:10} are obtained
by applying Goresky-MacPherson vanishing to the local cohomology exact
sequence associated to $H^{\bullet}_{E_i}(\widetilde{Y_i})$, noting
that as $Y_i$ is contractible, $H^{k}(\widetilde{Y_i})= H^{k}(E_i)$.  
\end{rem}

\begin{rem}
  \label{rem:41}
Fujiki shows \cite{MR602463} that the maps $\alpha_{k} \colon
H^{k}_{E_i}(\widetilde{Y_i}) \to H^{k}(E_i)$ satisfy $\alpha_{k}=
^{t}\alpha_{6-k}$ under the duality \eqref{eq:9}. In particular, this
shows that for $k<3$, the exact sequences:
\begin{eqnarray*}
  \label{eq:32}
  0 \to H^k_{E_i}(\widetilde{Y_i}) \to H^k(E_i) \to
H^{k+1}_{\{P_i\}}(Y_i)\to 0
\end{eqnarray*}
and
\begin{eqnarray*}
  0 \to H^{6-k}_{\{P_i\}}(Y_i) \to H^{6-k}_{E_i}(\widetilde{Y_i}) \to
H^{6-k}(E_i)\to 0
\end{eqnarray*}
are dual. The local cohomology groups $H^{k+1}_{\{P_i\}}(Y_i)$ and
$H^{6-k}_{\{P_i\}}(Y_i)$ are dual to each other. 
\end{rem}
\begin{lem}
\label{lem:28}
 The local cohomology of $Y$ at $P_i$ is:
 \begin{gather*}
H^{k}_{\{P_i\}}(Y)=(0), \quad \mbox{for}\quad k\neq 1,3,4,6\\
H^1_{\{P_i\}}(Y)\simeq H^6_{\{P_i\}}(Y) \simeq \C [E_i]\\
H^3_{\{P_i\}}(Y)\simeq H^4_{\{P_i\}}(Y) \simeq H^2(E_i)/ \sum
\C[E_i^j]
\end{gather*}
\end{lem}
\begin{proof}
First, notice that, by excision, $H^{k}_{\{P_i\}}(Y_i)= H^k_{\{P_i\}}(Y)$.
The local cohomology exact sequence shows that
$H^{0}_{\{P_i\}}(Y_i)=(0)$. As Remark~\ref{rem:41} shows, the result
need only be checked for $k=1,2$. 

The result is clear for $k=1$ by the
exact sequence \eqref{eq:34}, because $H^0_{E_i}(\widetilde{Y_i})$ is
dual to $H^6(E_i)=(0)$. 

For $k=2$, as the singularities are
assumed to be rational $H^1(\widetilde{Y_i},
\mathcal{O}_{\widetilde{Y_i}})=(0)$ and $H^{1}(E_i)\simeq
H^{1}(\widetilde{Y_i})=(0)$, the result follows from the exact
sequence \eqref{eq:34}.

The exact sequence \eqref{eq:34} gives:
\[
0 \to H^2_{E_i}(\widetilde{Y_i}) \to H^2(E_i)\to H^3_{\{P_i\}}(Y)\to 0.
\]
The $H^k(E_i)$ were determined in Section~\ref{sec:mixed-hodge-struct-2} and, by
Fujiki duality, $H^2_{E_i}(Y_i)$ is dual to $H^4(E_i)$, so that:
\[
 H^3_{\{P_i\}}(Y)= H^2(E_i)/ \bigoplus \C[E_i^j]. 
\]
The singularity $P_i \in Y$ is rational and $Y_i$ is a contractible neighbourhood:
\[
H^1(Y_i, \mathcal{O}_{Y_i}^{\ast})= H^2(Y_i,\Z)= H^2(\widetilde{Y_i},\Z)= H^2(E_i,\Z).
\]
The local cohomology of $Y$ at $P_i$ is of the form:
\[
 H^3_{\{P_i\}}(Y)= H^1(Y_i, \mathcal{O}_{Y_i}^{\ast})/ \bigoplus \C[E_i^j]. 
\]
\end{proof}

\section{The category of weak$\ast$ Fano 3-folds}
\label{sec:categ-weak-fano}

Let $Y$ be a terminal Gorentein non $\Q$-factorial Fano $3$-fold and
denote by $X$ a small $\Q$-factorialisation of $Y$. The Picard rank of
$X$ is equal to the rank of the Weil group of $Y$. If $Y$ has Picard
rank $1$ and defect
$1$, the Picard rank of $X$ is $2$. One of three things occurs: either $X$ has a
structure of del Pezzo fibration over $\PS^1$, or of conic bundle over
$\PS^2$, or there exists an extremal contraction $\phi \colon X \to
X'$, where $X'$ is a $\Q$-factorial terminal Fano $3$-fold with Picard
rank $1$. As $X$ has isolated hypersurface singularities, a direct geometric
analysis of the contraction $\phi$ is possible
(Theorem~\ref{thm:8}) and this describes the Weil non-$\Q$-Cartier divisors that can
lie on $Y$. 

One could hope to bound the defect and to obtain some information on
the Weil group of $Y$ by
running a Minimal Model Program (MMP) on $X$. The category of terminal
$\Q$-factorial $3$-folds is stable under the operations of the
MMP. In general, however, this approach is too naive: if $\phi \colon X
\to X'$ is an extremal contraction, $X'$ does not necessarily have
hypersurface singularities and nef anticanonical divisor. I show that
if $Y$ does not contain a plane, the MMP can be run on
$X$, and that the $3$-folds encountered when doing so are terminal
Gorenstein and have nef and big anticanonical divisor.

In this section, I define the category of \emph{weak$\ast$ Fano
  $3$-folds}. If $X$ is a weak$\ast$ Fano $3$-fold, its
  anticanonical model $Y$ is a terminal Gorenstein Fano $3$-fold that
  does not contain a plane. The $3$-fold $Y$ is in general not
  $\Q$-factorial. Conversely, any small $\Q$-factorialisation $X$ of a
  terminal Gorenstein Fano $3$ fold $Y$ whose anticanonical ring is
  generated in degree $1$ is a weak$\ast$ Fano unless
  $Y$ contains a plane. 

This section shows that the category of weak$\ast$ Fano $3$-folds is
stable under the operations of the MMP. If $X$ is a weak$\ast$ Fano,
the end product of the MMP on $X$ is described.

\subsection{Definitions and basic results}
\label{sec:defin-basic-results}

\begin{dfn}
  \label{dfn:1}
A $3$-fold $X$ is \emph{weak Fano} if $X$ is Gorenstein, terminal and if its anticanonical divisor $-K_X$ is nef and
big.

The \emph{anticanonical ring} of $X$ is $ R(X,-K_X)= \bigoplus_{n \in
  \N}H^0(X, -nK_X)$.

The \emph{anticanonical model} of $X $ is $Y=\proj R(X,-K_X)$.  
\end{dfn}
\begin{dfn}
  \label{dfn:2} 
$X$ is \emph{weak$\ast$} if moreover:
\begin{enumerate}
\item $X$ is $\Q$-factorial,
\item the morphism $X \to Y$ is small, \item $Y$ contains no plane $\PS^2$ with
$-K_Y|_{\PS^2} = \mathcal{O}_{\PS^2}(1)$, \item the anticanonical ring
$R=R(X,-K_X)$ is
generated by $R_1$.
\end{enumerate}
\end{dfn}
\begin{rem}
  \label{rem:2}
The category of weak$\ast$ Fano $3$-folds is stable under flops.  
\end{rem}
\begin{dfnlm}
  \label{dfnlm:1}
If $X$ is a weak Fano $3$-fold, then:
\[
  h^0(X, -mK_X)= 2m+1 +\frac{1}{12}
m(m+1)(2m+1)(-K_X)^3.
\]
Denote by $g=\frac{1}{2}(-K_X)^3 + 1$ the
\emph{genus} of $X$; in particular: 
\[
h^0(X,-K_X)=g+2.
\]
\end{dfnlm}
I recall here the statement of the
Kawamata-Viehweg vanishing theorem.
\begin{thm}[Kawamata-Viehweg vanishing,\cite{MR1658959}]
\label{thm:11}  
Let $(X, \Delta)$ be a proper Kawamata log terminal (klt) pair. Let $N$ be a $\Q$-Cartier Weil
divisor on $X$ such that $N \equiv M + \Delta$, where $M$ is a nef and
big $\Q$-Cartier $\Q$-divisor.
Then $H^i(X, \mathcal{O}_X(-N))=(0)$ for $i< \dim X$.  
\end{thm}
\begin{proof}[Proof of \ref{dfnlm:1}]
The Kawamata-Viehweg vanishing theorem applied to the pair $(X,0)$ and to the divisor
$-(m+1)K_X$ shows that: 
\[
H^i\bigl(X, \mathcal{O}_X((m+1)K_X)\bigl)=(0)
\]
for $i< 3$ and $m \geq 0$. By Serre duality, this implies \[H^j(X,
\mathcal{O}_X(-m K_X))=(0) \] for $j>0$ and $m \geq 0$ and  \[\chi
(\mathcal{O}_X(-mK_X))=h^0(X, -mK_X).\] The Riemann-Roch theorem for Gorenstein $3$-folds
\cite{MR927963} reads:
\[
   \chi(\mathcal{O}_X(-mK_X))=(1+2m)\chi (\mathcal{O}_X)+ \frac{1}{12}
m(m+1)(2m+1)(-K_X)^3.
\]
\end{proof}
\begin{lem}[Cone theorem for weak Fano $3$-folds]
  \label{lem:1}
If $X$ is a weak Fano $3$-fold then $NE(X)$ is a finite rational
polyhedron (in particular $NE(X)= \overline{NE}(X)$). If $R \subset NE(X)$ is an extremal ray,
then either $K_X \cdot R <0$, or $K_X \cdot R =0$ and there exists an
effective divisor $D$ such that $D \cdot R <0$. There is a contraction
morphism $\phi_R$ associated to each extremal ray $R$. If $\phi_R$ is
small, it is a flopping contraction.  
\end{lem}
\begin{proof}
  The $3$-fold $X$ has terminal singularities.
By the standard cone theorem \cite[Theorem 4-2-1]{MR946243}:
\[
\overline{NE}(X)= \overline{NE}(X)_{K_X \geq 0}+ \sum {C_j}
\]
where the extremal rays $C_j$ are discrete in the half space $\{
K_X <0 \}$ and can be contracted. Since $-K_X$ is nef, for any $z
\in N^1(X)$,  $K_X \cdot z \geq 0$ if and only if $K_X \cdot z
=0$. The anticanonical divisor $-K_X$ is big: for some
integer $m>0$, $m(-K_X) \sim A+D$, where $A$ is an ample divisor and
$D$ is effective. By the Nakai-Moishezon criterion for
ampleness, if $K_X \cdot z=0$ then $D \cdot z <0$. In particular, for
$0<\epsilon \ll 1 $,
\[
\overline{NE}(X) \subset {K_X}_{<0}   \cup  (K_X + \epsilon D)_{<0}.
\]
By the usual compactness argument, NE(X) is a finitely generated rational
polyhedron. Extremal rays can be contracted by the contraction theorem
\cite[Theorem 3-2-1]{MR946243}.
Finally, if $\phi_R$ is small, $R$ flips or flops. A flipping curve $\gamma$
on a terminal $3$-fold $X$ satisfies $-K_X \cdot \gamma <1$
\cite{MR796252}. The anticanonical divisor is Cartier and nef: $\phi_R$ is a flopping contraction.
\end{proof}
\begin{thm}
  \label{thm:1}
If $X$ is a weak Fano 3-fold then one of the following holds:
\begin{enumerate}
\item The linear system $|-K_X|$  has a non-singular
section. 
\item  The anticanonical model of $X$, $Y=\proj R(X, -K_X)$, is
  birational to a special complete intersection $X_{2,6} \subset
  \PS(1^4, 2,3)$ with a node. More precisely, $Y$ is defined by
equations of the form:
\[
\left\{ \begin{array}{c}
a_2=0\\
w^2+v^3+va_4 +a_6=0\\
\end{array}\right.
\]
where the coordinates $x_i$ have degree $1$, $v$ and $w$ have degree
$2$ and $3$ respectively, and where each $a_j$ is a homogeneous form of
degree $j$ in the variables $x_i$.
\end{enumerate}
\end{thm}
\begin{rem}
\label{rem:9}
If $X$ is a weak$\ast$ Fano $3$-fold, the linear system $\vert -K_X \vert$ is basepoint free. 
The anticanonical divisor of $X$ is nef and big, hence by the
basepoint free
theorem \cite[Theorem 3.3]{MR1658959}, the linear system $\vert -nK_X \vert
$ has no basepoint for $n$ sufficiently large.  
If the anticanonical ring of $X$ is
generated in degree $1$, the map
\[
H^0 (X, -K_X)^{\otimes n} \to H^0(X,-nK_X)
\]
is surjective for any $n \in \mathbb{N}$. The linear system $\vert
-K_X \vert $ is itself basepoint free. 
\end{rem}
I first prove the following lemma:
\begin{lem}
  \label{lem:2}
Let $X$ be a weak Fano $3$-fold. The general member $S$ of $\vert -K_X \vert$
is a $K3$ surface with no worse than Du Val singularities.  
\end{lem}  
\begin{dfn}
  \label{dfn:3}
  Let $X$ be a weak Fano $3$-fold. The \emph{Fano index} of $X$, $i(X)$, is
  the maximal integer such that $-K_X=i(X) H$ with
  $H$ a  nef and big Cartier divisor.
\end{dfn} 
\begin{rem}
  \label{rem:8}
The Fano index of a weak Fano $3$-fold $X$ with small anticanonical
map is the same as that of its anticanonical
model $Y$.
Indeed, the anticanonical map $f \colon X \to Y$ is
birational and small. As the divisor
$-K_X$ is relatively trivial, by the basepoint free/contraction
theorem \cite[Corollary 1.5]{MR924674}, $-K_X=
f^{\ast}D$, with $D$ ample. The map $f$ is crepant, so that $D=-K_Y$ and
 the Fano indices of $Y$ and $X$ are equal. 
\end{rem}
For a weak Fano $3$-fold $X$, by Kawamata-Viehweg vanishing and the
Riemann-Roch theorem:
\[
h^0(X, H)= 1+ \frac{2}{i(X)} + \frac{1}{12}H^3(1+i(X))(2+i(X)) \geq 2.
\]
\begin{lem}
 \label{lem:4}
$\Pic(X)$ has no torsion.
\end{lem}
\begin{proof} 
Let $D$ be a torsion divisor on $X$. By definition, there exists a smallest
integer $m>1$ such that $mD \sim 0$. The divisor $-K_X+D$ is nef and
big and therefore, by Kawamata-Viehweg vanishing,
\[
H^i(X, D)=H^i(X, K_X+(-K_X+D))=(0) \quad \text{for} \quad i>0.
\]
Moreover, $H^0(X, D)=(0)$, as otherwise $D$ itself would be linearly
equivalent to $0$. The
divisor $D$ has Euler characteristic
$\chi(\mathcal{O}_X(D))=(0)$. 

As $D$ is numerically trivial, by the Riemann-Roch theorem, 
\[
\chi(\mathcal{O}_X(D))=\chi(\mathcal{O}_X).  
\]
This contradicts $\chi(\mathcal{O}_X)=h^0(X, \mathcal{O}_X)=1$ on a
weak Fano $3$-fold. 
\end{proof}
\begin{rem}
  \label{rem:6}
  The divisor class $H$ of Definition~\ref{dfn:3} is uniquely determined.  
\end{rem}
\begin{dfn}
\label{dfn:17}
Let $X$ be a weak$\ast$ Fano $3$-fold of Fano index $i(X)$. The \emph{degree}
of $X$ is $H^3$. The \emph{anticanonical
  degree} of $X$ is $-K_X^3$.   
\end{dfn}

I use Kawamata's basepoint free technique to show that $\vert -K_X
\vert$ has a section with canonical singularities. I recall the
notion and properties of non-Kawamata log-terminal (non-klt) centres
introduced in \cite{MR1457742}. 
\begin{dfn}
\label{dfn:13}
Let $X$ be a normal variety and let $D= \sum d_i D_i$ be an effective
$\Q$-divisor such that $K_X+D$ is $\Q$-Cartier.
\begin{enumerate}
\item The \emph{non-klt locus} of $(X,D)$ is the set of points where $(X,D)$
is not klt, that is:
\[
\nklt(X,D)= \{ x \in  X: (X, D) \mbox{ is not klt at }  x\}.
\]
\item If the pair $(X,D)$ is log canonical (lc), a subvariety $W \subset X$
  is a \emph{log canonical centre} or \emph{lc centre}  for the pair
$(X,D)$ if there is a log resolution $\mu \colon \widetilde{X}\to
X$ and a prime divisor $E$ on $\widetilde{X}$ with discrepancy
coefficient $e = -1$ such that $\mu(E)=W$.
\end{enumerate}
\end{dfn}
\begin{rem}
  \label{rem:38}
Let $(X, D)$ be a pair and  $\mu \colon \widetilde{X}\to
X$ a log resolution. Denote by $E_1, \cdots , E_m $ the prime
divisors with discrepancy less than or equal to $-1$. The non-klt locus of
$(X,D)$ is:
\[ \nklt(X, D)= \mu(E_1+ \cdots E_m). \]
\end{rem}
\begin{rem}
\label{rem:39}
Let $(X,D)$ be a lc pair. A divisor $W$ is a codimension
$1$ lc centre if $D=W+D'$, where the support of $D'$ does not contain $W$.  
\end{rem}

\begin{pro}\cite{MR1457742}
\label{pro:7}
Let $X$ be a klt variety  and $(X, D)$ a lc pair. 
\begin{enumerate}
\item
Let $W_1$ and $W_2$ be lc centres of $(X,D)$. If $W$ is an irreducible
component of $W_1\cap W_2$, then $W$ is a lc centre as well. In
particular, if $x \in \nklt(X,D)$, there is a well defined minimal lc
centre containing $x$.
\item 
If $W$ is a minimal lc centre, then it is normal. 
\end{enumerate}   
\end{pro}
\begin{pro}[Subadjunction,\cite{MR1462926, Kol07}]
\label{pro:8}
Let $X$ be a normal variety with klt singularities and let $D$ be an
effective $\Q$-Cartier divisor such that the pair $(X,D)$ is log
canonical. Let $W \subset \nklt(X,D)$ be an isolated minimal lc centre. There is
an effective $\Q$-divisor $D_W$ on $W$ such that: 
\begin{enumerate}
\item The pair $(W,D_W)$ is klt,
\item 
$(K_{X}+D)_{\vert W} \sim_{\Q} K_{W}+D_W.$
\end{enumerate}
The choice of the divisor $D_{W}$ is not canonical: $D_{W}$
is the sum of a boundary divisor $B_{W}$ and of a general divisor
$M_W$ in a nef $\Q$-linear equivalence class $J(W,D)$. Write
$D=D'+D''$, where $D'$ (resp. $D''$) is the sum of components of $D$
that contain (resp. do not contain) $W$. The boundary part $B_W$ is uniquely
determined and is supported on $D''_{\vert W}$. The moduli part $J(W,D)$
is determined only by the pair $(X,D')$; the choice of $M_W \in
J(W,D)$ is not canonical.     
\end{pro}
\begin{proof}[Proof of \ref{lem:2}]
I follow ideas exposed in \cite{re,MR1181201,MR1675146}.
As is recalled in Definition-Lemma~\ref{dfnlm:1}, $h^0(X,-K_X)=g+2$,
where $g$ denotes the genus of $X$. Let $S$ be a general member of the
linear system $\vert -K_X \vert$. 
\setcounter{step}{0}
\begin{step}
The surface $S$ is connected.  
\end{step}
The sequence
\[
0 \to H^0(X, \mathcal{O}_X(-S)) \to H^0(X, \mathcal{O}_X) \to H^0(S,
\mathcal{O}_S) \to  H^1(X, \mathcal{O}_X(-S)) \to \cdots
\]
is exact. The surface $S$ is linearly equivalent to
$-K_X$. Kawamata-Viehweg vanishing shows that:
\[
H^0(X, \mathcal{O}_X) \simeq  H^0(S, \mathcal{O}_S)
\] and $S$ is connected.
\begin{setup}
The anticanonical divisor $-K_X$ is Cartier and nef and big and $X$ has terminal
singularities, hence by \cite[Proposition 2.61]{MR1658959}, there is a resolution $\mu \colon \widetilde{X}
\to X$ and a divisor with normal crossings $\sum E_i$ such that:
\begin{enumerate}
\item $ K_{\widetilde{X}}= \mu^{\ast}(K_X)+ \sum a_i E_i $, with $a_i
  \in \N$ and $a_i
  >0$ if and only if $E_i$ is exceptional;
\item $\mu^{\ast}(-K_X)-\sum p_i E_i$ is an ample $\Q$-Cartier
  divisor on $\widetilde{X}$, for suitable $p_i \in \Q, 0\leq p_i<\!\!<1$.
\end{enumerate}
A small perturbation of the rational numbers $p_i$ does not affect the
ampleness of $\mu^{\ast}(-K_X)-\sum p_i E_i $. 
We can blow up further to obtain a
new resolution $\mu \colon \widetilde{X} \to X$, that still satisfies
the conditions above, and that is a log resolution of the
linear system $\vert -K_X \vert$. Then:
\begin{eqnarray}
  \label{eq:6}
\mu^{\ast}(\vert -K_X \vert)= \vert L \vert + \sum r_i E_i, 
\end{eqnarray}
where $\vert L \vert$ is a free linear system, $r_i \in \N$ and
$r_i>0$ if $\mu(E_i)$ is in the base locus of $\vert -K_X \vert$. 
The divisor $\sum E_i$ has simple normal crossings. 

Note that $\mu$ determines a log resolution of the general member $S$
of $\vert -K_X \vert$: the divisor $\mu^{\ast}(S)+ \Exc(\mu)$ has
simple normal crossing support. Indeed, if the divisor $\sum E_i$ has
simple normal crossings on the non-singular variety $\widetilde{X}$
and if $L$ is a general section of the free linear system $\vert L \vert $,
then $L+ \sum E_i$ has simple normal crossings as well. Then:
\begin{eqnarray}
  \label{eq:26}
  \mu^{\ast}(S)=L+ \sum r_i E_i,
\end{eqnarray}
where $L$ is a member of the free linear system $\vert L \vert$.  
Reid proves that the free linear system $\vert L \vert$ is not
composed with a pencil \cite{re}.
\end{setup}
\begin{step}
The surface $S$ is irreducible, reduced and has no worse than du Val
singularities if no component of the non-klt locus of $(X,S)$ is
contained in the base locus of $\vert -K_X \vert$.    
\end{step}
If the free linear system $\vert L \vert$ is not composed with a
pencil, by Bertini's theorem, its general member is reduced and irreducible; the surface $S$ is
irreducible and reduced if the linear system $\vert -K_{X} \vert$ has no base component. 
Recall that a base
component of $\vert -K_X \vert$ is the image by $\mu$ of a divisor
$E_i$ with $a_i=0$ and $r_i \geq 1$. 

The surface $S$ has canonical
singularities if the discrepancy of every exceptional divisor
appearing in a resolution of $S$ is non-negative. The restriction of
$\mu$ to $\mu_{\ast}^{-1}(S)=L$ is a resolution of $S$. By Bertini, $S$ is
non-singular away from the base locus of $\vert -K_X \vert$, so that
$S$ has canonical singularities if $a_i-r_i \geq 0$ (that is: $a_i-r_i
>-1$, as $a_i$ and $r_i$ are natural numbers) for every exceptional
divisor with centre contained in the base locus of $\vert -K_X \vert$.

In the formalism introduced above, the surface $S$ is reduced,
irreducible and has no worse than canonical singularities, if no
component of the non-klt locus of $(X,S)$ is contained in the
base locus of $\vert -K_X \vert$. 

\begin{rem}
The rational numbers $a_i$ are non-negative because $X$ is
terminal. If $W_i= \mu(E_i)$ is a lc centre, the coefficient $r_i$ is
positive and $W_i$ is contained in the base locus of $\vert -K_X \vert$.  
\end{rem}

\begin{step}
Assume that the pair $(X,S)$ is not purely log terminal. For some $0<b
\leq 1$, the pair $(X,bS)$ is log canonical but not klt.
\end{step}

The pair $(X,S)$ is not log terminal, therefore there exists an index
$i_0=0$ such that $a_0 -r_0 \leq -1$, that is $a_0+1 \leq r_0$. Define: 
\begin{eqnarray}
  \label{eq:28}
 b = \min \{\frac{a_i+1}{r_i}\} \leq 1, 
\end{eqnarray}
where the minimum is taken over the indices $i$ such that $r_i$ is
non-zero.
The rational number $b$ is positive because $a_i \geq 0$ for all $i$.
The pair $(X,bS)$ is strictly log canonical, as the prime divisors
$E_i$ for the indices where the minimum is attained have discrepancy $-1$.

\begin{step}
There is a divisor $D \in \vert -K_X\vert$ and a rational number
$k \in \Q$, such that the non-klt locus of $(X, b(1-\epsilon)S+
\epsilon kD)$ is irreducible for all $0<\epsilon<\!\!<1$. The lc centre
$W= \nklt(X, b(1-\epsilon)S+\epsilon kD)$ is a minimal lc centre for $(X,D)$.
\end{step}

Let $E_i$ for $0 \leq i \leq d$ be the prime divisors such that: 
\[
\nklt(X, bS)= \mu(E_0)\cup \cdots \cup \mu(E_d),
\]
so that:
\begin{equation}
  \label{eq:27}
K_{\widetilde{X}}= \mu^{\ast}(K_X+bS)-bL+\sum_{i\geq
  d+1}(a_i-br_i)E_i-\sum_{0 \leq i \leq d}E_i  
\end{equation}
with $a_i-br_i>-1$ for all $i\geq d+1$.

As a slight increase of the $p_i$ does not affect the ampleness of the
$\Q$-Cartier divisor
$\mu^{\ast}(-K_X)-\sum p_i E_i$, we may assume that the minimum
\[
k= \min_{0\leq i\leq d}\{ \frac{br_i}{p_i} \}
\]
is attained only for one index $i_0=0$.

Let $D$ be an element of the linear system $\vert -K_X \vert$ such
that $\mu^{\ast}(D)-\sum p_i E_i$ is an ample divisor; from
\eqref{eq:27}, we get, for $0 <\epsilon <\!\!<1$:
\begin{multline}
  \label{eq:44}
K_{\widetilde{X}}= \mu^{\ast}(K_X+b(1-\epsilon)S+\epsilon
  kD)-b(1-\epsilon)L-\epsilon k(\mu^{\ast}(D)-\sum p_i E_i)\\
+\sum_{i\geq
  d+1}(a_i-br_i+ \epsilon(br_i-kp_i))E_i
+ \sum_{1 \leq i \leq d}(-1+
  \epsilon (br_i-kp_i))E_i -E_0.   
\end{multline}
By construction, for all $1\leq i \leq d$ and all $\epsilon >0$, \[-1+
  \epsilon (br_i-kp_i)>-1.\] 
For $\epsilon$ sufficiently small, $a_i-br_i+
  \epsilon(br_i-kp_i)>-1$ for $i \geq d+1$ and $\epsilon k(
  \mu{\ast}(-K_X)-\sum p_i E_i)$ is a boundary divisor. 

From now on, I include $L$ and a Cartier multiple of the ample divisor
$\mu{\ast}(-K_X)-\sum p_i E_i$ among
the divisors $E_i$. Denote by $b_{i, \epsilon}$ the discrepancy
coefficient of $E_i$ and set:
\[
E=E_0, \quad  F=\sum_{i \neq 0, b_{i, \epsilon}<0} -b_{i, \epsilon}
  E_i, \quad \mbox{and
  } A=\sum_{b_{i,\epsilon}\geq 0}b_{i,\epsilon}E_i.
\]
The divisor $A$ is effective and exceptional, since
$b_{i,\epsilon}>0$ for $0<\epsilon <\!\!<1$ is only possible for
$a_i>0$, i.e. for $E_i$ exceptional. The divisor $F$ is a boundary,
and \eqref{eq:44} reads:
\begin{equation}
  \label{eq:29}
  K_{\widetilde{X}}= \mu^{\ast}(K_X +b(1-\epsilon)S+\epsilon
  kD)+A-F-E.
\end{equation}
The non-klt locus of $(X, b(1-\epsilon)S+\epsilon kD)$ is the
irreducible component $\mu(E)$.

\begin{step}
If $W$ is the unique lc centre for $(X,b(1-\epsilon)S+ \epsilon kD)$, then $W$ is not contained
in $\Bs \vert -K_X \vert$.  
\end{step}
Consider the divisor:
\begin{equation}
  \label{eq:46}
N(t)= \mu^{\ast}(-tK_X)-(K_{\widetilde{X}}+F)+A-E.  
\end{equation}

From equation \eqref{eq:29} we see that: 
\[
N(t)= \mu^{\ast}(-tK_X) -\mu^{\ast}(K_X +b(1-\epsilon)S+\epsilon kD)
\]
and $N(t)$ is numerically equivalent to
$\mu^{\ast}\big(-(t+1-b+\epsilon (b-k))K_X\big)$; $N(t)$ is nef and
big for $t+1-b+\epsilon (b-k)>0$.

By the Kawamata-Viehweg vanishing theorem,
\[
H^{1}(\widetilde{X}, N(t)+K_{\widetilde{X}}+F)=(0)
\] 
for $t+1-b+\epsilon (b-k)>0$. The constant $b$ is less than or equal
to $1$; for $t=1$ and $\epsilon$ sufficiently small, 
\[
H^{1}(\widetilde{X}, \mu^{\ast}(-K_X)+A-E)=0,
\]
and the map
\[
H^0(\widetilde{X}, \mu^{\ast}(-K_X)+A) \to H^0 (E,
(\mu^{\ast}(-K_X)+A)_{\vert E})
\]
is surjective. 

By the Leray spectral sequence, $H^0(\widetilde{X},
\mu^{\ast}(-K_X)+A) \simeq H^0(X, -K_X)$, and 
any section of
$H^0(\widetilde{X}, \mu^{\ast}(-K_X)+A)$ that does not vanish on $E$
pushes forward to a section of 
$H^0(X ,-K_X)$ not vanishing identically on $W$; if
$H^0 (E,(\mu^{\ast}(-K_X)+A)_{\vert E})$ is non-trivial, $W$ is not
contained in the base locus of $\vert-K_{X} \vert$.

The following cases need to be considered:
\begin{enumerate}
\item If $\codim(W)=3$, $W$ is a closed point $\{x\}$. Then:
\[
H^0(\widetilde{X}, \mu^{\ast}(-K_X)+A) \to H^{0}(E, (\mu^{\ast}(-K_{X})+A)_{\vert E}) \simeq H^{0}(E,
A_{\vert E}) \to 0, 
\]
and the Leray spectral sequence shows that $h^0(\widetilde{X},A)= h^0(X,
\mathcal{O}_X)=1$. The divisor $A_{\vert E}$ is effective, therefore
$h^0(E, A_{\vert E})=1$. Hence, $\{x\}$ does not
belong to the base locus of $\vert -K_X \vert$.
\item If $\codim(W)=2$, $W$ is a curve $\Gamma$. By
  Propositions~\ref{pro:7} and \ref{pro:8}, $\Gamma$ is non-singular and
  there is an effective divisor $M_{\Gamma}$ on $\Gamma$ such that:
\[
(K_{X}+b(1-\epsilon)S+ \epsilon kD)_{\vert \Gamma}=(1-b)K_X= K_{\Gamma} +M_{\Gamma}.
\]
The divisor $K_X$ has non-positive degree on $\Gamma$ and
$K_{X}+b(1-\epsilon)S+ \epsilon kD$ is numerically
equivalent to $(1-b-\epsilon (k-b))K_X$, hence, as $0<\epsilon<\!\!<1$ can be
taken arbitrarily small, $\Gamma$ has arithmetic genus $p_a(\Gamma)=0$ or $1$. 
A flopping
curve on $X$ is rational \cite{MR986434}: if $\Gamma$ is elliptic,
$(-K_{X})_{\vert \Gamma}$ has positive degree. In both cases, $h^0(\Gamma,
{-K_{X}}_{\vert \Gamma})>0$ and 
\[
H^{0}(E,(\mu^{\ast}(-K_X)+A)_{\vert E}) \supset
H^{0}(E,(\mu^{\ast}(-K_X))_{\vert E}) \simeq H^{0}(\Gamma,
{-K_X}_{\vert\Gamma}) \neq (0)
\]
shows that $\Gamma$ is not contained in $\Bs \vert -K_{X} \vert$. 
\item If $\codim(W)=1$, $f_{\ast}(E)=W$ is a base component of the linear system
  $\vert -K_X \vert$.
 
Denote by $P(t)$ the polynomial of degree $2$ defined as the Euler
characteristic $\chi((t\mu^{\ast}(-K_X)+A)_{\vert E})$ of the
linear system $\vert t\mu^{\ast}(-K_X)+A \vert_{\vert E} $ on the
divisor $E$.

By definition of the divisor $N(t)$ \eqref{eq:46}, 
\[
(t\mu^{\ast}(-K_X)+A)_{\vert E}= (N(t)+K_{\widetilde{X}}+E+F)_{\vert E}.
\]
The divisor $E$ is Cartier, hence $(K_{\widetilde{X}}+E)_{\vert
  E}=K_E$; the support of $F\cup E$ has simple normal crossings
  because $\mu$ is a log resolution, the divisor $F_{\vert E}$
  therefore is a boundary divisor and $P(t)=
  h^0((t\mu^{\ast}(-K_X)+A)_{\vert E}))$ by the Kawamata-Viehweg
  vanishing for $t+1-b+\epsilon (b-k)>0$.
Recall that $W= \mu(E_0)$ is a minimal lc centre of codimension $1$ of the pair
  $(X,bS)$, hence $W$ is a fixed component of $\vert -K_X \vert$. 
From~(\ref{eq:28}), the multiplicity of the fixed component $W$ in
  $\vert -K_X \vert$ is $\frac{1}{r_0}$. 

If $b=1$, $r_0=1$ and $W$ is a fixed component of $\vert -K_X \vert$
with multiplicity $1$. Since the divisor $-K_X$ is
connected and movable, $S\in \vert -K_X \vert$ is singular along a
codimension $2$ set $\Gamma \subset W$. This contradicts $W$ being a
minimal lc centre of $(X, S)$. 

The multiplicity $r_0$ of the fixed component
$W$ is therefore strictly greater than $1$ and $b<1$. 
For $\epsilon$ sufficiently small, $1-b+\epsilon (b-k)>0$ and
$P(t)=h^0((t\mu^{\ast}(-K_X)+A))$ for all $t \geq 0$. 

The divisor $W$ is a base component of the linear system $\vert -K_X
\vert$, therefore $P(1)= h^0(E, (\mu^{\ast}(-K_X)+A)_{\vert E})=
h^0(W, -K_X)=0$. As $h^0(\widetilde{X},A)= h^0(X,
\mathcal{O}_X)=1$ and $A_{\vert E}$
is effective, $P(0)=1$. The leading coefficient of $P(t)$ is
$\frac{((\mu^{\ast}(-K_X))_{\vert E})^2}{2!}=\frac{(-K_X)^2\cdot
  W}{2!}$ and $P(t)$ can be written:
\[
P(t)= \frac{1}{2} \big((-K_X)^2\cdot W \cdot t^2+(-(-K_X)^2\cdot
W-2)\cdot t+2 \big).
\]
By Riemann-Roch, however,
\begin{eqnarray*}
P(t)= \frac{1}{2}((t\mu^{\ast}(-K_X)+A)_{\vert E})\cdot
((t\mu^{\ast}(-K_X)+A)_{\vert E}-K_E)+1 \\
= \frac{1}{2}(t\mu^{\ast}(-K_X)+A)_{\vert E})^2-\frac{1}{2}K_E \cdot((t\mu^{\ast}(-K_X)+A)_{\vert E}))+1.
\end{eqnarray*}
The leading coefficients in both expressions agree, equating the
coefficients of the terms of degree $1$ in $t$ shows that:
\begin{eqnarray}
  \label{eq:30}
  -(-K_X)^2\cdot W-2=(A_{\vert E}-\frac{1}{2}K_E)\mu^{\ast}(-K_X)_{\vert E}.
\end{eqnarray}
The anticanonical divisor $-K_X$ is nef and big, so that the
expression on the left hand side is strictly negative. 
The divisor $E$ is Cartier, therefore, by adjunction:
\begin{eqnarray*}
K_E=(K_{\widetilde{X}}+E)_{\vert E}=
(\mu^{\ast}(-K_X+b(1-\epsilon)S+ \epsilon k D)+A-F)_{\vert E}\\ 
\sim
\big(\mu^{\ast}(-(1-b+\epsilon(b-k))K_X)+A-F\big)_{\vert E}
\end{eqnarray*}
where $A$ is an effective exceptional divisor and $F$ is
effective and supported on $\cup_{1\leq i\leq d}E_i$. 

The left hand side of \eqref{eq:30} becomes (up to a factor of $2$):
\begin{eqnarray*}
(A+ F)_{\vert E} \cdot \mu^{\ast}(-K_X)_{\vert E}+\mu^{\ast}(-(1-b+\epsilon(b-k))K_X)_{\vert E} \cdot
\mu^{\ast}(-K_X)_{\vert E}\\
\sim (A+ F)_{\vert E} \cdot \mu^{\ast}(-K_X)_{\vert
  E}+(1-b+\epsilon(b-k))(\mu^{\ast}(-K_X)_{\vert E})^2 .
\end{eqnarray*}

 The divisor $(A+F)_{\vert E}$ is effective and the first term is
 non-negative because the divisor $\mu^{\ast}(-K_{X})_{\vert E}$ is
 nef. Recall that $\epsilon$ has be chosen so that $(1-b+\epsilon(b-k))>0$. 

The left hand side of \eqref{eq:30} is positive; this yields a
contradiction. The linear system $\vert -K_W \vert$ has no base
component.
\end{enumerate}

I have shown that the linear system $\vert -K_X \vert$ has no base
component. The free
linear system $\vert L \vert$ is not composed with a pencil \cite{re} and hence,
by Bertini,
the general
member $S \in \vert -K_X \vert$ is irreducible and reduced.  
For all indices $i$ of exceptional divisors $E_i$
whose centre is contained in the base locus of $\vert -K_X\vert$, $a_i-r_i \geq 0$.  
The surface $S$ is irreducible, reduced and has canonical singularities. 

\begin{step}
The general member $S \in \vert -K_X  \vert$ is a $K3$ surface with
rational double points.  
\end{step}
I first recall the definition of $K3$ surfaces:
\begin{dfn}
  \label{dfn:4}
  A $K3$ surface is a surface $S$ with no worse than Du Val singularities
  such that:
  \begin{enumerate}
\item $H^1(S, \mathcal{O}_S)=(0)$,
\item The canonical divisor class $K_S$ is trivial.
  \end{enumerate}
\end{dfn}

Consider the long exact sequence in
cohomology associated to 
\begin{equation}
\label{eq:10}
0 \to \mathcal{O}_X (-S) \to \mathcal{O}_X \to \mathcal{O}_S \to 0.
\end{equation}
By Kawamata-Viehweg vanishing, 
\[
H^1(X, \mathcal{O}_X)= H^2(X, \mathcal{O}_X(-S))=(0),
\]
and $H^1(S, \mathcal{O}_S)=(0)$.
By adjunction, $K_S=(K_X+S)_{\vert S} \sim 0$ and $S$ is a $K3$ surface with
rational double points.
\end{proof}
\setcounter{step}{0}

Tensor the exact sequence \eqref{eq:10} with $\mathcal{O}_X(-K_X)$. The
surface $S$ is a section of $\vert -K_X \vert$. By the
Kawamata-Viehweg vanishing theorem, the map \[H^0(X, -K_X) \to H^0(S,
{-K_X}_{\vert S})\] is surjective. It follows that the base points of
$\vert -K_X \vert$ are precisely the base points of $\vert
{-K_X}_{\vert S} \vert$.

I study the linear system $\vert {-K_X} _{\vert S}\vert $ on
the surface $S$ to complete the proof of Theorem~\ref{thm:1}. Linear systems
on non-singular K$3$ surfaces were studied by Saint-Donat
\cite{MR0364263}. I recall his --- now classical --- results below.

\begin{thm}{\cite{MR0364263}}
  \label{thm:7}
  Let X be a non-singular K$3$ surface and $D$ a nef and big
  effective divisor on $X$. 
One of the following is true:
  \begin{enumerate}
  \item The linear system $\vert D \vert$ is basepoint free, and:
    \begin{enumerate}
    \item Either: the morphism $\phi_{\vert D \vert}$ induced by
$\vert D \vert$ on $X$ is
birational on its image $\overline{X}$, $\overline{X}$ has rational
double points, and $\phi_{\vert D \vert}$ is the minimal resolution. 
    \item Or: $D$ is hyperelliptic and $\phi_{\vert D\vert}$ is
generically $2$-to-$1$.
    \end{enumerate}
  \item $D$ is monogonal, that is $D \sim kE + \Gamma $, with $\vert E
  \vert $ a free elliptic pencil and $\Gamma$ a $(-2)$-curve, which is
  the fixed component of
  $\vert D \vert$. In that case, $\Gamma^2=-2$, $\Gamma \cdot E=1$.
  \end{enumerate}
\end{thm}
I prove the following easy extension of Saint-Donat's result  to
linear systems on K$3$ surfaces with rational double points:   
\begin{thm}
  \label{thm:2}
Let $S$ be a $K3$ surface with rational double points and let $D$ be a
nef and big Cartier divisor on $S$. Then one of the following holds:

\begin{enumerate}
\item $|D|$ is basepoint free. \item $D$ is monogonal. In that
case:
\[D \sim kE+\Gamma,\] where $|E|$ is a free elliptic pencil,
$E\cdot \Gamma=1$ and $\Gamma$ is a $(-2)$-curve contained in the
non-singular locus of $S$, 
\item $S$ is birational to a (special)
complete intersection $X_{2,6} \subset \PS(1^3,2,3)$ with an ordinary
double point given by:
\[
\left\{ \begin{array}{c}
a_2=0\\
z^2+y^3+ya_4 +a_6=0\\
\end{array} \right.
\]
where $a_i$ is a homogeneous form of degree $i$ in the coordinates $x_1,x_2,x_3$.
\end{enumerate}  
\end{thm}

\begin{proof}
Let $f \colon T \to S$ be the minimal resolution of $S$. The surface
$T$ is a non-singular K$3$ surface. Indeed, the exceptional locus of
$f$ consists of a set of $(-2)$-curves $G_i$. By adjunction, $K_T \cdot
G_i=0$. Moreover, $K_T$ may be written as
\[ K_T=f^*(K_S)+ \sum_i a_i G_i= \sum_i a_i G_i,
\] with $a_i \geq 0$. The matrix $(G_i \cdot G_j)$ is negative definite,
hence $a_i=0$ for all $i$ and $K_T$ is trivial. 

The singularities of
$S$ are rational and, in particular,
$R^1f_{\ast}\mathcal{O}_T=0$. By the Leray spectral sequence,
$H^1(S, \mathcal{O}_S)=H^1(T, \mathcal{O}_T)=(0)$.

Consider the complete linear system $\vert D'\vert=f^* \vert D
\vert= f^* \vert {-K_X} _{\vert S}\vert$ on $T$. The divisor 
$D'$ is nef and big: the results of \cite{MR0364263} show that
$\vert D' \vert$ is basepoint free or $D'$ is monogonal.

In the first case, $\Bs \vert D' \vert= \Bs \vert D \vert =
\emptyset$. 

Now assume that $D'= kE + \Gamma$, the base locus
 of $\vert D'\vert$ is a $(-2)$-curve $\Gamma$ and $\vert E \vert$ is
 a free elliptic pencil. 

Let $G$ be an exceptional irreducible curve of the resolution $f$.
Then, \[G\cdot D'=0= kG \cdot E+ G\cdot \Gamma \] and $E$ is nef. There are
two possible configurations:
\begin{enumerate}
\item Either $G$ is disjoint from both $E$ and $\Gamma$, and $\Bs
  \vert {-K_X}_{\vert S} \vert= f_{\ast} (\Gamma)$ is a $(-2)$-curve
  contained in the non-singular locus of $S$, so we are in the second
  case of the theorem.
\item Or $G=\Gamma$ and, as $-kE\cdot G=k=G^2$, this implies $k=2$. Thus,
  $D' \sim 2E + \Gamma$, and the base locus $\Bs \vert {-K_X}_{\vert S}
  \vert= \{p\}$ is an ordinary double point on the surface $S$.   
\end{enumerate}

Assume that $G=\Gamma$. By the observation above, any section $C$ of
the linear system
$\vert {-K_X}_{\vert S} \vert$ is of the form $C= E_1+ E_2$ with
$E_1\cdot E_2 =1$. The intersection of $E_1$ and $E_2$ is a reduced
point $\{p\}$.

I determine explicitly the graded ring 
\[
R(S, {-K_X}_{\vert S})=
\bigoplus_{n \in \N} H^0(S, -n{K_X}_{\vert S})
\] and show that the
linear system $\vert -n {K_X}_{\vert S} \vert$ defines a birational
morphism from $S$ to $\proj R(S, {-K_X}_{\vert S})$.

Let $C$ be a section of the linear system $\vert {-K_X}_{\vert S} \vert$ on the K$3$
surface $S$. By Kawamata-Viehweg vanishing on $S$, the long exact
sequence in cohomology associated to
\[
0 \to \mathcal{O}_S((n-1)({-K_X}_{\vert S}))\to
\mathcal{O}_S(n({-K_X}_{\vert S}))\to \mathcal{O}_C(-n{K_X}_{\vert S})
\to 0
\]
shows that  
\[R(C, {-K_X}_{\vert S})=R(S, {-K_X}_{\vert S})/sR(S,
{-K_X}_{\vert S}), 
\]
where $s$ is a section of $D= {-K_X}_{\vert S}$.
The exact sequence on $S$
\[
0 \to \mathcal{O}_{E_1 \cup E_2} \to \mathcal{O}_{E_1} \oplus
\mathcal{O}_{E_2} \to \mathcal{O}_{E_1 \cap E_2} \to 0
\]
gives rise to the long exact sequence in cohomology:
\begin{gather*}
 \ldots \to H^i(C,\mathcal{O}_{C}(nD)) \to H^i(E_1,
 \mathcal{O}_{E_1}(nD)) \oplus  H^i(E_2,
 \mathcal{O}_{E_2}(nD)) \to \ldots \\
 \to H^i(C, \mathcal{O}_{\{ p \}}(nD)) \to
 H^{i+1}(C,\mathcal{O}_{C}(nD)) \to \ldots.
\end{gather*}
By Kawamata-Viehweg vanishing, $H^1(C, \mathcal{O}_C(nD))=(0)$, and therefore: 
\begin{gather*}
0 \to H^0(C,\mathcal{O}_C(nD)) \to \ldots \\
\ldots \to H^0(C,\mathcal{O}_{E_1}(nD)) \oplus
 H^0(C,\mathcal{O}_{E_2}( nD)) \stackrel{\alpha} \to  H^0(C,\mathcal{O}_{\{p\}}(nD)) \to  0,
\end{gather*}
in particular:
\[
R(C, D)=(R(E_1, D) \oplus R(E_2, D))/(\alpha R(p,D)).
\] 
By the above, 
\[
 H^0(C,\mathcal{O}_{E_i}(nD))=H^0(E_i,\mathcal{O}_{E_i}(nD))=
H^0(E_i,\mathcal{O}_{E_i}(np))
\]
 and the map $\alpha $ is a relation of degree $2$. 

Applying Riemann-Roch's theorem to the linear system $\vert np \vert$ on the
non-singular elliptic curve $E_i$ yields:
\[
R(E_i, p)= R(E_i,D_{\vert{E_i}}) \simeq \C[x_i,v,w]/(\phi), 
\]
where $x_i$, $v$, $w$ have degrees $1$, $2$ and $3$ respectively and
$\phi$ is a relation of degree $6$. 
Finally, the map $\alpha$ is a relation of
degree $2$ of the form $\{ b_2(x_1,x_2)=0 \}$, where $b_2$ is a
homogeneous map of degree $2$. This
determines $R(C,D)$ and $R(S,D)$.

Hence $\proj R(S, D_{\vert S})$ is a special complete intersection
\[
X_{2,6} \subset \PS(1,1,1,2,3)=\PS(x_1,x_2,x_3,v,w)
\]
with an ordinary double point. 

The map $S \to \proj R(S,D_{\vert S})$ is birational, induced by
$\vert 2D_{\vert S} \vert$ because $R(S,D_{\vert S})$ is generated in
degrees less or equal to $2$. The proof of Theorem~\ref{thm:2} is
complete.
\end{proof}

\begin{rem}
  \label{rem:11}
In the third case of Theorem~\ref{thm:2}, the map induced by $\vert
-K_{X} \vert$ is
generically $2$-to-$1$. 
The map
\[
S \to \{ a_2(x_1, x_2, x_3)=0 \} \subset \PS^2 \subset \PS(1,1,1,2,3)
\]  
is $2$-to-$1$.
\end{rem}

\begin{proof}[End of proof of \ref{thm:1}]
By \ref{lem:2}, the general section $S$ of $\vert -K_X \vert$ is a
K$3$ surface with rational double points. Consider the
long exact sequence in cohomology associated to:
\begin{equation}
\label{eq:11}
0 \to \mathcal{O}_X (-nK_X-S) \to \mathcal{O}_X(-nK_X) \to
\mathcal{O}_S({-nK_X}_{\vert S}) \to 0.
\end{equation}
By Kawamata-Viehweg vanishing, since $-nK_X-S \sim -(n-1)K_X$ with
$-K_X$ nef and big, $H^1(X,
\mathcal{O}_X(-nK_X-S))=(0)$ for $n \geq 0$. In particular, $H^0(X, -K_X) \to H^0(S,
{-K_X}_{\vert S})$ is surjective and $\Bs \vert -K_X \vert=\Bs {\vert
-K_{X\vert S}}\vert$. By Bertini's theorem, a general member $S$
can only have singularities at base points of $\vert -K_X \vert$. 

According to \ref{thm:7}, one of the following holds:
\begin{enumerate}
\item $\vert -K_X \vert $ is basepoint free,
\item $\Bs\vert -K_X \vert = \Gamma$, where $\Gamma$ is a $(-2)$-curve
  contained in the non
  singular locus of $X$, or
\item $\Bs\vert -K_X \vert =\{p\}$ an ordinary double point. 
\end{enumerate}
In the first case, $\vert -K_X \vert$ has a non-singular member. As is noted in Remark~\ref{rem:9}, if $X$ is a weak$\ast$ Fano
$3$-fold, this is the only possible case.

In the second case, let $S \in \vert -K_X \vert$ be a general $K3$
section. Theorem~\ref{thm:2} shows that $S$ is non-singular along
$\Bs \vert -K_X \vert$, hence (by Bertini) it is non-singular everywhere.

In the third case, the anticanonical map is a birational map
(determined by $\vert -2K_{X} \vert$) to
$\proj R(X, -K_X)$. The exact sequence \eqref{eq:11} implies the ``hyperplane section principle'': 
\[
R(S, {-K_X}_{\vert S}) \simeq R(X, -K_X)/(s R(X, -K_X)),
\] where $s$
is a variable of degree $1$. The anticanonical model of $X$ is
therefore birational to a
special $X_{2,6} \subset \PS(1^4, 2, 3)$ with a node, given by
equations of the following form (the variables $x_i$ have degree $1$,
$v$ and $w$ have degrees $2$ and $3$ respectively):
\[
\left\{ \begin{array}{c}
a_2=0\\
z^2+y^3+ya_4 +a_6=0\\
\end{array} \right.
\]
where each $a_j$ is a homogeneous form of degree $j$ in the variables
$x_i$, $i= 0, \cdots , 4$.

To conclude the proof of Theorem~\ref{thm:1}, notice that, as is mentioned
in Remark~\ref{rem:9}, if $X$ is weak$\ast$, the linear system $\vert
-K_{X} \vert$ is basepoint free.
\end{proof}
\begin{rem}
  \label{rem:3}
Weak Fano $3$-folds such that $\Bs \vert -K_X \vert \simeq \PS^1$ is
a $(-2)$-curve are called monogonal. One can show that the only
occurences are:
\begin{enumerate}
\item $X= \PS^1 \times S$, where $S$ is the weighted
hypersurface $S_6 \subset \PS(1,1,2,3)$; $X$ has Picard rank $10$,
\item
$X=\Bl_{\Gamma} V$, the blow up of $V$, the weighted hypersurface $V_6
\subset \PS(1^3,2,3)$ along the plane $\Pi=\{ x_0=x_1=0\}$; $X$ has
Picard rank $2$.
\end{enumerate}
If $X$ is monogonal, the linear
system $\vert -K_X\vert$ determines a rational map from $X$ to a surface.
\end{rem}

\subsection[MMP for weak$\ast$ Fano $3$-folds]{Minimal Model Program for weak$\ast$ Fano $3$-folds}
\label{sec:minim-model-progr}

In \cite{MR936328}, Cutkosky studies the contractions of extremal
rays on projective $3$-folds with normal Gorenstein $\Q$-factorial
terminal singularities. He proves the following results:

\begin{thm}[Birational operations of the minimal model program]
  \label{thm:8}
Let $f \colon X \to Y$ be the birational contraction of an extremal
ray. Assume that $f$ is not an isomorphism in codimension $1$. Then $Y$ is factorial and:
  \begin{description}
  \item[E$1$:] Suppose that $f \colon X \to Y$ contracts a surface $E$ to a
  curve $\Gamma$. Then $Y$ is non-singular near $\Gamma$, $\Gamma$ is
  locally a complete intersection and has planar singularities: in the local
  ring $\mathcal{O}_{Y, p}$ of any point $p \in \Gamma$, one of the
  local equations of $\Gamma$ is a smooth
  hypersurface near $p$. The contraction $f$ is
  the blow up of the ideal sheaf $I_{\Gamma}$. The $3$-fold $X$ has only $cA_n$ singularities on $E$.
  \item Suppose that $f \colon X \to Y$ contracts a surface $E$ to a
  point $p$. Then $f$ is one of the following:

\item[E$2$:] The $3$-fold $Y$ is non-singular, $E \simeq \PS^2$ and
  $\mathcal{O}_E(E)\simeq \mathcal{O}_{\PS^2}(-1)$: $f$ is the inverse
  of the blow up of a non-singular $3$-fold at a point.
\item[E$3$:] The local ring of $Y$ at $p$ is of the form
  $\mathcal{O}_{Y,p}= k[x,y,z,w]/(x^2+y^2+z^2+w^2)$, $E \simeq \PS^1 \times \PS^1 $ and $\mathcal{O}_E(E)=
  \mathcal{O}(-1,-1)$: $f$ is the inverse of the blow up of an
  ordinary double point.
\item[E$4$:] The local ring of $Y$ at $p$ is of the form
$\mathcal{O}_{Y,p}= k[x,y,z,w]/(x^2+y^2+z^2+w^n)$, with
  $n \geq 3$ and $E \simeq Q$, an irreducible reduced singular quadric surface
  in $\PS^3$, with normal bundle $\mathcal{O}_E(E) \simeq
  \mathcal{O}_E \otimes \mathcal{O}_{\PS^3}(-1)$: $f$ is the inverse
  of the blow up of a $cA_{n-1}$ singular point.
\item[E$5$:]  The local ring of $Y$ at $p$ is of the form
  $\mathcal{O}_{Y,p}= k[[x,y,z]]^{(2)}$, the ring of invariants for
  the $\Z_2$-action. $E \simeq \PS^2$ with
  $\mathcal{O}_E(E) \simeq \mathcal{O}_{\PS^2}(-2)$: $f$ is the inverse
  of the blow up of a non Gorenstein point of index $2$. 
\end{description}
\end{thm}
\begin{thm}[Conic bundles]
  \label{thm:9}
Suppose that $f \colon X \to Y$ is the contraction of an extremal ray
to a surface $Y$. Then $Y$ is non-singular and $X$ is a possibly
singular conic bundle over $Y$.  
\end{thm}

I want to determine whether I can run a Minimal Model Program in
the categories of $3$-folds I have introduced in
Section~\ref{sec:defin-basic-results}. The category of
weak Fano $3$-folds is not suitable, since a birational contraction
of type E$5$ would take us out of the category. I prove that 
the category of weak$\ast$ Fano $3$-folds, by contrast, is stable under the
operations of the Minimal Model Program.
 
\begin{thm}
  \label{thm:3}
The category of weak$\ast$ Fano $3$-folds is stable under the
birational operations of the Minimal Model Program.  
\end{thm}

\begin{proof}  
Let $X$ be a weak$\ast$ Fano $3$-fold and let $\phi:X
\to X'$ be a divisorial contraction. By
Theorem~\ref{thm:8}, the variety $X'$ is $\Q$-factorial and has
terminal singularities.
\setcounter{step}{0}
\begin{step}
The extremal contraction $\phi$ is not of type E$5$.  
\end{step}
By the description of extremal contractions given in
Theorem~\ref{thm:8}, if $\phi$ is of type E$5$, $X$ contains a Cartier
divisor
$E\simeq \PS^2$ with normal bundle $\mathcal{O}_{E}(E) \simeq
\mathcal{O}_{\PS^2}(-2)$. By adjunction, the restriction of the
anticanonical divisor to $E$ is ${-K_{X}}_{\vert \PS^2}=
\mathcal{O}_{\PS^2}(3-2)= \mathcal{O}_{\PS^2}(1)$. The rational map $\Phi_{\vert
-K_X \vert_{\vert E}}$ coincides with a projection from a possibly
empty linear subspace
\[
\xymatrix{\nu\colon \PS(H^0(E, \mathcal{O}_{\PS^2}(1)))\simeq \PS^2 \ar@{-->}[r]
& \PS(H^0(E, \vert -K_X\vert_{\vert E})}
\] 
associated to the inclusion  $\vert -K_X \vert_{{\vert E}}
\subset \vert \mathcal{O}_{\PS^2}(1) \vert$.
 By
Theorem~\ref{thm:1}, the linear system $\vert -K_X \vert_{\vert E}$ is
basepoint free, so that $\Phi_{\vert -K_X\vert_{\vert E}}$
is a morphism, and $\vert -K_X \vert_{\vert E}=
\vert \mathcal{O}_{\PS^2}(1) \vert$.
The anticanonical model $Y$ contains a plane
$\PS^2$ with ${-K_{Y}}_{\vert \PS^2}= \mathcal{O}_{\PS^2}(1)$: this
contradicts $X$ being weak$\ast$.

\begin{step}
The anticanonical divisor $-K_{X'}$ is nef and big.  
\end{step}
I first establish the following lemma:
\begin{lem}
  \label{lem:3}
Let $X$ be a weak Fano $3$-fold and let $\phi\colon X \to X'$ be a
divisorial extremal contraction. Then:
\begin{enumerate}
\item The anticanonical divisor $-K_{X'}$ is big. 
\item One of the following holds:
\begin{enumerate}
\item The anticanonical divisor $-K_{X'}$ is nef,
\item  The contraction $\phi$ is of type E$1$; its exceptional divisor
  $E$ is isomorphic to $\F_1$. The negative section $\sigma$
of $E$ is contracted by the anticanonical map, i.e. $-K_{X} \cdot
\sigma=0$.
\end{enumerate}
\end{enumerate}  
\end{lem}

\begin{proof}[Proof of \ref{lem:3}] 
The morphism $\phi$ is an extremal divisorial contraction and 
\begin{eqnarray}
  \label{eq:31}
   -K_X=\phi^{\ast}(-K_{X'})-aE
\end{eqnarray}
where $E$ is the effective prime exceptional divisor of $\phi$ and $a$
is a positive natural number. More
precisely, $a=2$ when $\phi$ is of type E$2$ and $1$ otherwise (recall
that $\phi$ is not of type E$5$).

In particular, the anticanonical rings of $X$ and $X'$ satisfy 
\[
R(X, -K_X) \subset R(X', -K_{X'}),
\]
so that the anticanonical model of $X'$ has dimension at least
$3$: $-K_{X'}$ is big.

Assume first that $\phi$ contracts a divisor $E$ to a point. 
Let $Z'$ be an effective irreducible curve lying on $X'$. 
Denote by
$Z= \phi_{\ast}^{-1}(Z')$ the proper transform of $Z'$. The curve $Z$
is effective and irreducible and it either instersects $E$ properly or
not at all. In particular, by the projection formula:
\begin{eqnarray*}
-K_{X'} \cdot Z'= \phi^{\ast}(-K_{X'}) \cdot Z= -K_X \cdot Z+aE\cdot Z
\geq -K_X\cdot Z \geq 0.
\end{eqnarray*}
The anticanonical divisor $-K_{X'}$ is nef.

Assume now that $\phi$ contracts a divisor $E$ to $\Gamma$, an irreducible
reduced curve. The anticanonical divisors of $X$ and $X'$ satisfy:
\begin{eqnarray*}
\phi^{\ast}(-K_{X'})= -K_X+E  .
\end{eqnarray*}
Let $Z'$ be an irreducible effective curve lying on $X'$ and denote by
$Z$ any effective irreducible curve on $X$ that maps $1$-to-$1$ to
$Z'$.
There are two cases to consider: either $Z'$ and $\Gamma$ intersect in
a $0$-dimensional (or empty) set or $Z'=\Gamma$. 

If $Z'$ and $\Gamma$ intersect in a $0$-dimensional or empty
  set, $E \cdot Z \geq 0$ and by the projection formula: 
  \begin{eqnarray*}
-K_{X'}\cdot Z'= -K_{X}\cdot Z+aE \cdot Z \geq -K_{X}\cdot Z \geq 0.    
  \end{eqnarray*}
 
If $Z'= \Gamma$, by contrast, 
  \begin{eqnarray*}
 -K_{X'}\cdot \Gamma= -K_{X}\cdot Z+ E \cdot Z.
  \end{eqnarray*}
As above $-K_X \cdot Z \geq 0$ but $E \cdot Z$ can be negative and
$-K_{X'}$ can fail to be nef.

If $-K_{X'}$ is not nef, $-E_{\vert E}$ is ample. The
divisor $E$ is Cartier and, by adjunction,
\begin{eqnarray*}
-K_{E}= (-K_{X}-E)_{\vert E}  
\end{eqnarray*}
is ample as the sum of a nef and an ample divisor. The surface $E$ is
a possibly nonnormal del Pezzo surface. 
\begin{cla}
The surface $E$ and the curve $\Gamma$ are normal.   
\end{cla}

Theorem~\ref{thm:8} states that the curve $\Gamma$ has planar
singularities and that at any point $P \in \Gamma$, one of the local
equations of $\Gamma$ is a smooth hypersurface near $P$. Locally, the
equation of $\Gamma$ is of the form $\{x=f(y,z)\}$. Near a singular
point $P \in \Gamma$, the blow up of $\Gamma$ (and hence $E$) is given
by the equation:
\begin{eqnarray*}
   \{ t_0 x-t_1
f(y,z)=0\} \subset \PS^1_{t_0, t_1} \times \C^3.
\end{eqnarray*}
Writing down the equations of $E$ on the affine pieces of the blow up
shows that $E$ is a $\PS^1$-bundle over $\Gamma$. 

Consider the normalisation maps $\Gamma^{\nu} \to \Gamma$ and $E^{\nu} \to E$.
The surface $E^{\nu}$ has a structure of $\PS^1$-bundle over $\Gamma^{\nu}$ which
makes the following diagram commutative:
\[
\xymatrix{ E^{\nu} \ar[r]^{n} \ar[d]_{g^{\nu}} & E \ar[d]^g\\
 \Gamma^{\nu} \ar[r]^n  & \Gamma }.
\]
Denote by $\Delta$ (resp. $\delta$) the divisor of $E^{\nu}$
(resp. $\Gamma^{\nu}$) defined by the conductor ideal of the
normalisation map \cite{MR1311389}, so that
\begin{align*}
n^{\ast}(K_{\Gamma})= K_{\Gamma^{\nu}}+ \delta\\
n^{\ast}(K_{E})= K_{E^{\nu}}+ \Delta
\end{align*} 
where $\Delta = g^{\nu \ast}(\delta)$ is an effective sum of fibres of $E^{\nu} \to
\Gamma^{\nu}$. The Weil divisor $\Delta$ is defined scheme
theoretically by the conductor ideal $\mathcal{I}_{\Delta, E^{\nu}}$,
i.e. the inverse image by $n$ of the ideal
$\Ann(n_{\ast}(\mathcal{O}_{E^{\nu}})/ \mathcal{O}_E)$. As a set,
$\Delta$ is the codimension $1$ double locus of $n$. 

The anticanonical divisor $-K_{E^{\nu}}$ is ample.
Indeed, $n^{\ast}(-K_{E})$ is ample as the pullback of an ample
divisor by a finite morphism and $\Delta$ is nef as an effective sum
of fibres of the $\PS^1$-bundle $E^{\nu}$  over a non-singular curve $\Gamma^{\nu}$.
The surface $E^{\nu}$ is a normal del Pezzo and has a structure of
$\PS^1$-bundle over a non-singular curve $\Gamma^{\nu}$. 
According to the classification of normal del Pezzo surfaces,
$E^{\nu}$ is isomorphic to $\F^1$. 

Let $\sigma$ be a negative section of $E^{\nu}$.  
Since $n^{\ast}(-K_{E})$ is ample,
\[
 (-K_{E^{\nu}} -\Delta) \cdot \sigma= 1- \Delta \cdot \sigma
>0.  
\]
Consequently, $\Delta$ is empty and $E$ is normal. Moreover, by adjunction:
\begin{eqnarray*}
-K_{E}\cdot \sigma =   (-K_{X}-E)_{\vert E} \cdot \sigma=1
\end{eqnarray*}
and as $-E_{\vert E}$ is ample and $-K_X$ is nef, $-E_{\vert E}\cdot
\sigma=1$ and $-K_X \cdot \sigma =0$.
\end{proof}

\begin{cla}
The anticanonical divisor $-K_{X'}$ is nef and big: $X'$ is a weak
Fano $3$-fold.
\end{cla}
Lemma~\ref{lem:3} shows that $-K_{X'}$ is always big and that it fails to be
nef only if $\phi$ is an E$1$ contraction and if the exceptional
divisor $E$ is a rational scroll $\F_1$ with negative section
$\sigma$ such that $-K_{X} \cdot \sigma=0$.

I show that this case does not occur when $X$ is a weak$\ast$ Fano
$3$-fold. 

The computations in the proof of Lemma~\ref{lem:3} show that $\vert
-K_{X} \vert_{ \vert E}$ is contained in the linear system $\vert
\sigma +f \vert$ on the scroll $\F_1 =E$ (the linear system that
contracts the negative section). Following the usual
notation conventions, I denote by $\sigma$ the negative section on the
scroll, by $f$ a fibre and the intersection numbers are $\sigma^2=-1$,
$\sigma f=1$ and $f^2=0$.

By Theorem~\ref{thm:1}, the anticanonical linear system is
basepoint free.
The basepoint free linear system $\vert -K_{X}\vert_{\vert E}$
determines a morphism onto a surface $\overline{E}$. More precisely, the morphism
$\Phi_{\vert -K_{X}\vert_{\vert E}}$ has a factorisation $\nu \circ \Phi_{\vert \sigma+f
  \vert}$ where $\nu$ is the projection from a (possibly empty) linear
subspace 
\[
\xymatrix{\PS(H^0(E, \vert \sigma+f \vert))\simeq \PS^2 \ar@{-->}[r]&
  \PS(H^0(E,\vert -K_{X}\vert_{\vert E}})
\] 
associated to the inclusion $\vert -K_{X}\vert_{\vert E} \subset \vert
\sigma +f \vert$. The surface $\overline{E}$ is the projection of
$\PS^2=\Phi_{\vert \sigma+f \vert}(E)$ from a possibly empty linear
subspace: it is equal to $\PS^2$. This shows that $\vert -K_X
\vert_{E}$ is the complete linear system $ \vert \sigma+f \vert$. The
anticanonical model $Y$ then contains a plane $\PS^2=\overline{E}$,
which is the image of the scroll $E =\F_1$ by the
anticanonical map. Moreover, ${-K_{Y}}_{\vert \PS^2}=
\mathcal{O}_{\PS^2}(1)$ and this contradicts $X$ being weak$\ast$.

\begin{step}
The anticanonical map $X' \to Y'$ contracts finitely many curves.
\end{step}

Let $Z'$ be an irreducible effective curve of $X'$ that is contracted
by the anticanonical map of $X'$. Denote by $E$ the exceptional
divisor of the contraction $\phi$ and recall that
\begin{eqnarray*}
\phi^{\ast}(-K_{X'})= -K_{X}+aE  
\end{eqnarray*}
for some positive integer $a$. Either $Z'$ meets 
the centre of the contraction $\phi(E)$ in a $0$-dimensional or empty set or $\phi$ is an E$1$
contraction with centre $Z'$. If $Z'$ meets $\phi(E)$ in a $0$-dimensional or empty set, denote by
$Z$ its proper transform. Then, by the projection formula,
$-K_{X'}\cdot Z'=0$  if and only if
$E \cdot Z=0$ and $-K_{X} \cdot Z=0$, that is if $Z$ is itself
contracted by the anticanonical map of $X$ and does not meet the
centre of the contraction $\phi$. The $3$-fold $X$ is weak$\ast$: there
are finitely many such curves. The anticanonical map $X' \to Y'$ is small.

\begin{step}
 The anticanonical ring $R(X', -K_{X'})$ is generated in degree $1$. 
\end{step}

The proof of Theorem~\ref{thm:1} shows that the anticanonical ring
$R(X',-K_{X'})$ of $X'$ is generated in degree $1$ if and only if the
rational map $\Phi_{\vert -K_{X'} \vert}$ determined by $\vert -K_{X'} \vert$ is birational onto its
image. More precisely, it shows that if $\Phi_{\vert -K_{X'}
  \vert}$ is not birational onto its image, either $\Phi_{\vert -K_{X'}
  \vert}$ is generically $2$-to-$1$ and $Y'$ is birational to a
special complete intersection $X_{2,6} \subset \PS(1^4, 2,3)$ with a
node, or $\Phi_{\vert -K_{X'}\vert}$ maps $X'$ to a
surface ($X'$ is then monogonal).

Recall that
\[
\phi^{\ast}(-K_{X'})= -K_{X}+aE,  
\]
for some positive integer $a$, and hence:
\[
0 \to H^{0}(X, -K_{X}) \to H^{0}(X, \phi^{\ast}(-K_{X'})) \simeq
H^{0}(X', -K_{X'}). 
\]
The linear system $\vert -K_{X} \vert$ is naturally a subsystem
of $\vert \phi^{\ast} (-K_{X'}) \vert$. This embedding defines a natural projection
from a linear subspace:
\[
\begin{xy}
\xymatrix{\nu \colon\PS(H^0(X,\phi^{\ast} (-K_{X'})))\ar@{-->}[r] &\PS((H^0(X, -K_{X}))).}  
\end{xy}
\]
The rational map $\Phi_{\vert -K_{X}\vert}$ determined by  $\vert -K_{X} \vert$
is  birational onto its image because $X$ is a weak$\ast$ Fano
$3$-fold.
The rational map $\Phi_{\vert -K_{X}\vert}$ factorises as $\nu \circ\Phi_{\vert \phi^{\ast}(-K_{X'})
  \vert}$, hence the dimension of $\Phi_{\vert -K_{X}\vert}(X)$ is less
than or equal to that of $\Phi_{\vert
  \phi^{\ast}(-K_{X'})\vert}(X)=\Phi_{\vert -K_{X'}\vert}(X')$, and
$X'$ is not monogonal. The map
$\Phi_{\vert -K_{X}\vert}$ cannot be birational onto its image if
$\Phi_{\vert -K_{X'}\vert}$ is generically $2$-to-$1$. Hence, the map
$\Phi_{\vert -K_{X'}\vert}$ is birational onto its image.

The anticanonical ring of $X'$ is generated in degree $1$. This
implies that $\vert -K_{X'}\vert$ is basepoint free (Remark~\ref{rem:9}).

\begin{step}
The anticanonical model $Y'$ of $X'$ does not contain a plane $\PS^2$
with ${-K_{Y'}}_{\vert \PS^2}= \mathcal{O}_{\PS^2}(1)$.  
\end{step}
Denote by $g$ (resp. $g'$) the anticanonical map of $X$ (resp. $X'$).
Assume that $X'$ contains a surface $S$ that is mapped by $g'$ to a plane $\PS^2$ with
${-K_{Y}}_{\vert \PS^2}= \mathcal{O}_{\PS^2}(1)$. The linear system
$\vert -K_{X'} \vert$ is basepoint free  and determines the
anticanonical map of $X'$, therefore, as above, $\vert -K_{X'} \vert_{\vert S}=
\vert (g')^{\ast}\mathcal{O}_{\PS^2}(1) \vert$. Denote by
$\widetilde{S}$ the proper transform of $S$. The linear system
\begin{eqnarray*}
\vert -K_{X} \vert_{\vert \widetilde{S}} = \vert \phi^{\ast}(-K_{X'})-aE
\vert _{\vert \widetilde{S}}   
\end{eqnarray*}
is a subsystem of $\vert \phi^{\ast}(-K_{X'}) \vert_{\vert \widetilde{S}}=\vert
(g'\circ\phi)^{\ast}\mathcal{O}_{\PS^2}(1)\vert_{\vert \widetilde{S}}$. The linear system $\vert -K_{X}
\vert_{\widetilde{S}}$ is strictly contained in $\vert \phi^{\ast}(-K_{X'})
\vert_{\vert \widetilde{S}}$ if the surface $S$ intersects the centre of the
contraction $\phi$. By assumption on $X$, $\vert -K_{X}
\vert_{\vert \widetilde{S}}$ is basepoint free and determines a
morphism $\Phi_{\vert -K_{X} \vert_{\vert \widetilde{S}}}$
from $\widetilde{S}$ to a surface. This morphism factorises as
$\Phi_{\vert -K_{X} \vert_{\vert \widetilde{S}}}= \nu \circ\Phi_{\vert
(g'\circ\phi)^{\ast}\mathcal{O}_{\PS^2}(1)\vert_{\vert
  \widetilde{S}}}$ where $\nu$ is the projection from a possibly empty
linear subspace 
\[
\xymatrix{\PS(H^0(\widetilde{S}, 
{(g'\circ\phi)^{\ast}\mathcal{O}_{\PS^2}(1)}_{\vert
  \widetilde{S}}))\simeq \PS^2 \ar@{-->}[r] & \PS(H^0(\widetilde{S},\vert -K_{X}
\vert_{\widetilde{S}}))} 
\] 
associated to the inclusion $\vert -K_{X}
\vert_{\widetilde{S}}\subset \vert(g'\circ\phi)^{\ast}\mathcal{O}_{\PS^2}(1)\vert_{\vert
  \widetilde{S}}$. The image of $\widetilde{S}$ by $\Phi_{\vert -K_{X}
  \vert_{\vert \widetilde{S}}}$ is therefore $\nu(\PS^2)=
\nu(\Phi_{\vert(g'\circ\phi)^{\ast}\mathcal{O}_{\PS^2}(1)\vert_{\vert
  \widetilde{S}}}(\widetilde{S}))$. The linear system $\vert -K_{X}
\vert_{\vert \widetilde{S}}$ determines a morphism onto a surface: $\nu$
is the identity, 
$\vert -K_X\vert_{\vert \widetilde{S}}$ is complete and equal to $
\vert \mathcal{O}_{\PS^2}(1) \vert$. This contradicts $X$
being a weak$\ast$ Fano $3$-fold.
\end{proof}
\setcounter{step}{0}
The Minimal Model Program can therefore be run in the category of
weak$\ast$ Fano $3$-folds. 

\begin{lem}
\label{lem:12}
If $X := X_0$ is a weak$\ast$ Fano $3$-fold whose anticanonical
divisor $Y_0$ has Picard rank $1$, there is a sequence of extremal contractions:
\[
\xymatrix{
X_0 \ar[r]^{\phi_1} \ar[d]& X_1 \ar[r]^{\phi_2}\ar[d] & \cdots
& X_{n-1}\ar[r]^{\phi_n} \ar[d]& X_n \ar[d] \\
Y_0 & Y_1 & \cdots & Y_{n-1} & Y_n
}
\]
where for each $i$, $X_i$ is a weak$\ast$ Fano $3$-fold, $Y_i$ is its anticanonical model, and
$\phi_i$ is a birational contraction of an extremal ray. Moreover, for
each $i$, the Picard rank of $Y_i$, $\rho(Y_i)$ is equal to $1$. The
$3$-fold $X_n$ either has Picard rank $1$ or is an extremal Mori fibre space. 
\end{lem}

\begin{proof}
Assume that $\rho(X_i)>1$. By Lemma~\ref{lem:1}, there is an
extremal ray $R$ on $X_i$ and $R$ can be contracted. If the contraction
associated to $R$ is not birational, $X_i=X_n$ is an extremal Mori fibre space
and there is nothing to prove. 

If the
extremal contraction 
\[ \phi_i\colon X_i\to X_{i+1}
\] is birational,
Theorem~\ref{thm:3} shows that $X_{i+1}$ is a weak$\ast$ Fano
$3$-fold. I prove that the anticanonical model $Y_{i+1}$ of $X_{i+1}$ has Picard rank $1$.

If $\phi_i$ is a flopping contraction, there is nothing to prove as $Y_i=Y_{i+1}$.
Let us assume that $\phi_i$ is a divisorial contraction of type
E$1$-E$4$. Recall that the proof of Theorem~\ref{thm:3} shows that
extremal contractions of weak$\ast$ Fano $3$-folds are not of type E$5$. 

Let $E$ be the exceptional divisor of the contraction $\phi_i$ and let $f_i$
be the anticanonical map of $X_i$. 

\begin{step}
The image of $E$ by the anticanonical map $\overline{E}=f_i(E)$ is a
Weil non $\Q$-Cartier divisor.  
\end{step}
The divisor $E$ is covered by $K_{X_i}$-negative rational curves $\Gamma$ that have
strictly negative intersection with $E$, $E \cdot \Gamma<0$. If the
divisor $\overline{E}$ is Cartier, $\overline{E}$ is covered by curves
$\overline{\Gamma}=(f_i)_{\ast}(\Gamma)$ and by the projection formula:
\[
\overline{E}\cdot \overline{\Gamma}= \overline{E}\cdot
(f_i)_{\ast}(\Gamma)= (f_i)^{\ast}(\overline{E})\cdot \Gamma <0.
\]
The $3$-fold $Y_i$ has Picard rank $1$: $\overline{E}$ is not
$\Q$-Cartier because it is not ample.

Denote by 
$Z_i=\PProj \bigoplus_{n\geq 0}
{f_i}_{\ast}\mathcal{O}_{X_i}(nE)$ a small partial
$\Q$-factorialisation of $Y_i$; $Z_i$ is the symbolic blow up of $Y_i$
along the Weil non $\Q$-Cartier divisor $f_i(E)$. 
\[
\xymatrix{E \subset X_i \ar[r]^{\phi_i} \ar[d]_{g_i} & X_{i+1}\\
E' \subset Z_i \ar[d]_{h_i} & \\
\overline{E} \subset Y_i &
}
\]  
In the above diagram, $h_i \circ g_i= f_i$, $E = f_i^{\ast}(\overline{E})$. 
By construction, the divisor $E'$ is 
Cartier on $Z_i$.
\begin{step}
The $3$-fold $Z_i$ has Picard rank $2$. There is an extremal
contraction on $Z_i$ that contracts $E'$.  
\end{step}
The $3$-fold $Z_i$ is the crepant blow up of $Y_i$ along a single Weil
non $\Q$-Cartier divisor $\overline{E}$: $\rho(Z_i/Y_i)=1$. 
The ($\Q$-Cartier) divisor $E'$ is covered by curves
$\Gamma'$ such that $E'\cdot \Gamma'<0$ and $-K_{Z_i}\cdot \Gamma'>0$. Indeed, the map $g_i$ is small and
$E= g_i^{\ast}(E'), -K_{X_i}= g_i^{\ast}(-K_{Z_i})$. The divisor $E$ is
covered by curves $\Gamma$ such that $E\cdot \Gamma
<0$ and $-K_{X_i}\cdot \Gamma<0$; by the projection formula, so is $E'$. The cone
$\overline{NE}(Z_i)$ is $2$-dimensional and rational polyhedral; by the contraction theorem
\cite[Theorem 3-2-1]{MR946243}, there exists an extremal ray
$R$ on which $E'$ is negative and this extremal ray may be contracted. 

Denote by $\psi_i$ the contraction of the extremal face $E'$; $\psi_i$
fits in the following diagram:
\[
\xymatrix{E \subset X_i \ar[r]^{\phi_i} \ar[d]_{g_i} & X_{i+1} \\
E'\subset Z_i \ar[d]_{h_i} \ar[r]^{\psi_i} & Z_{i+1} \\
Y_i &
}
\]  
and $\rho(Z_i/Z_{i+1})=\rho(X_i/X_{i+1})=1$. 

\begin{step}
  There is a small map $g_{i+1}\colon X_{i+1}\to Z_{i+1}$ such that
  $\psi_i\circ g_i= g_{i+1}\circ \phi_i$.  
\end{step}
Consider the projective and surjective morphism
\[
\psi_i \circ g_i \colon X_i \to Z_{i+1}
\]
and run a relative Minimal Model Program on $X_i$ over $Z_{i+1}$. The
divisor $E \in \overline{NE}(X_i/Z_{i+1})$ is covered by curves $\Gamma \in N_1(X_i/Z_{i+1})$ with $E\cdot
\Gamma<0$ and $-K_{X_i}\cdot \Gamma>0$, because $E$ is covered by such curves $\Gamma\in
N_1(X_i)$ and these are contracted by $\psi_i \circ g_i$ by definition
of $\psi_i$. The contraction theorem shows that the contraction $\phi_i \colon
X_i\to X_{i+1}$ of the
extremal face $E$ factorises $\psi_i\circ g_i$ and makes the diagram
\[
\xymatrix{E \subset X_i \ar[r]^{\phi_i} \ar[d]_{g_i} & X_{i+1}\ar[d]^{g_{i+1}} \\
E'\subset Z_i \ar[d]_{h_i} \ar[r]^{\psi_i} & Z_{i+1} \\
Y_i &
}
\]  
commutative.

The map $g_{i+1}$ is crepant because $\rho(X_i/X_{i+1})=
\rho(Z_i/Z_{i+1})$ and $g_i$ maps the exceptional divisor of $\phi_i$ to
that of $\psi_i$.

\begin{step}
The anticanonical model of $Z_{i+1}$ is $Y_{i+1}$. The Picard rank of
$Y_{i+1}$ is $1$.  
\end{step}

Since $K_{X_{i+1}}= g_{i+1}^{\ast}(K_{Z_{i+1}})$, $X_{i+1}$ and $Z_{i+1}$ have the same anticanonical
model. The Picard rank of $Z_{i+1}$ is $\rho(Z_i)-1=1$ because
$\overline{E}$ is $\Q$-Cartier, hence the
anticanonical model of $Z_{i+1}$ has Picard rank $1$. 
\setcounter{step}{0}
\end{proof}

I now study strict Mori fibrations. Let $\phi_{n+1}$ be an extremal
contraction that is a strict Mori fibration, i.e. a del Pezzo
fibration over a curve $\Gamma$ or a conic bundle over
a surface $S$. 

If $\phi_{n+1}$ is a del Pezzo fibration over a curve $\Gamma$,
$\Gamma$ has arithmetic genus $0$ by the Leray spectral sequence
because $h^1(X_n,\mathcal{O}_{X_n})=0$. The
curve $\Gamma$ is isomorphic to $\PS^1$.

\begin{lem}
\label{lem:35}
 If $\phi \colon X \to \PS^1$ is a weak$\ast$ Fano $3$-fold that is an
 extremal del Pezzo
 fibration of degree $k$, $k$ is not equal to $1$ or $2$. 
\end{lem}
\begin{proof}
A general fibre $F$ of $\psi$ is a non-singular del Pezzo surface
of degree $k$. As $-K_{X\vert F}= -K_F$, the linear system $\vert -K_X
\vert_{\vert F}$ is naturally a subsystem of $\vert -K_F
\vert$. Theorem~\ref{thm:1} shows that $\vert -K_X \vert$ is basepoint
free, hence the degree $k$ cannot be equal to $1$ as the anticanonical
linear system of a non-singular del Pezzo surface of degree $1$ has
base points.   

Since $X$ is a weak$\ast$ Fano, the anticanonical map is birational,
contracts finitely many curves, and is determined by the linear system
$\vert -K_X \vert$. The restriction of the anticanonical map to $F$
factorises as $\Phi_{\vert -K_X \vert_{\vert F}}=\nu \circ \Phi_{\vert
  -K_F \vert}$, where $\nu$ is the projection from a possibly empty
linear subspace
\[
\xymatrix{\nu \colon \PS(H^0(F, -K_F)) \ar@{-->}[r] & \PS(H^0(F, -K_{X\vert F}))}
\] 
naturally associated to the inclusion of linear systems $\vert -K_X
\vert_{\vert F} \subset \vert -K_F \vert$.
 
This is impossible if $k=2$, as $\Phi_{\vert -K_F\vert}$ is generically
$2$-to-$1$.  
\end{proof}
 
Assume that $\phi_{n+1}$ is a conic bundle, then Cutkosky
shows (Theorem~\ref{thm:9}) that the
surface $S$ is non-singular. 

\begin{dfn}
  \label{dfn:5}
 A conic bundle $\phi \colon X \to S$ is a \emph{weak Fano conic
  bundle} if $X$ is a weak Fano $3$-fold and $\rho(X/S)=1$.
\end{dfn}

\begin{lem}
  \label{lem:17}
If $\phi \colon X\to S$ is a weak Fano conic bundle, $-K_S$ is nef.
\end{lem}

\begin{proof}
By definition of the discriminant curve $\Delta$ of a conic bundle, 
\[
-4K_{S}= \phi_{\ast}(-K_X)^2 +\Delta.
\]

Assume that $-K_S$ is not nef and let $C \subset S$ be an irreducible
curve such
that $-K_S \cdot C <0 $. The curve $C$ is necessarily contained in
$\Delta$, as $\phi_{\ast}(-K_X)^2$ is nef and $C^2\leq C\cdot \Delta <0$.
By adjunction:
\begin{eqnarray*}
-2 \leq 2p_a(C)-2 \leq (K_S+C)\cdot C\\
  \leq (K_S+\Delta)\cdot C \leq (-3K_S-\phi_{\ast}(-K_X)^2)\cdot C \leq
  -3K_S\cdot C.
\end{eqnarray*}
This is impossible because $-K_S \cdot C$ is an integer.
\end{proof}

If the surface $S$ is not minimal, one can run a Minimal Model Program
on $S$. There is
a chain of contractions of $(-1)$-curves $ S \to S_1 \to \ldots \to
S_N$ such that either $K_{S_N}$ is nef, or $S_N$ is $\PS^2$ or a
$\PS^1$-bundle over a curve $\Gamma$.

I show the following basic lemma:
\begin{lem}
  \label{lem:6}
Let $\phi \colon X \to S$ be a weak Fano conic bundle. Then, there exists a weak
  Fano conic bundle $\phi' \colon X' \to S'$, with $S'= \PS^2, \F_0$ or $\F_2$,
  such that the following diagram is commutative:
\[
\xymatrix{ X \ar[r] \ar[d]& X' \ar[d]\\
 S \ar[r]  & S' }.
\]
\end{lem}
\setcounter{step}{0}
\begin{proof} I show that steps of the Minimal Model Program on
  $S$ are dominated by extremal contractions of $X$. 
\begin{step}
 Let $X \to S$ be a weak Fano conic bundle and $S \to S'$ the
 contraction of a $(-1)$-curve. Then, there exists a weak Fano
 conic bundle $X' \to S'$ and an extremal
 contraction $X \to X'$ making the following diagram
 commutative:
\[
\xymatrix{ X \ar[r] \ar[d]& X' \ar[d]\\
 S \ar[r]  & S' }.
\]   
\end{step}
 The relative Picard rank $\rho(X/S')$ is equal to $2$, so that by
 Lemma~\ref{lem:1}, $\overline{NE}(X/S')$ has exactly $2$ extremal
 rays that can be contracted. A $2$-ray game gives the
 following two possible situations:
\[
\xymatrix{ X \ar[r]^{\psi} \ar[d]& X' \ar[d] & \mbox{or} & X \ar[r]^{f} \ar[d]& T \ar[d]\\
 S \ar[r]  & S'  & \quad  &  S \ar[r]  & S'}
\]
where $\psi$ is birational and $f$ is a conic bundle. 

In the first case, denote by $\Gamma$ the $(-1)$-curve on $S$
contracted to a point $P \in S'$. Then, by definition $X_{\vert S \smallsetminus
  \Gamma} \simeq X'_{\vert S'\smallsetminus\{P\}}$, and $X' \to S'$
naturally has a structure of weak Fano conic bundle.

In the second case, all maps are
proper. Considering analytic neighbourhoods of the contracted $(-1)$-curves
in $S$ and in $T$
shows that $S$ and $T$ are isomorphic. The two maps $X \to S$
and $X \to T$ correspond to the contraction of the same extremal ray.  

There is a sequence of contractions of $(-1)$-curves $S \to S_1 \to
\ldots \to S_N$ such that $K_{S_N}$ is nef or $S_N$ is $\PS^2$ or
$\PS^1 \times \PS^1$.
\begin{step}
The minimal surface $S_N$ is
birational to $\PS^2$ or a $\PS^1$-bundle over a rational curve.  
\end{step}
As $-K_{S_N}$ is nef, the surface $S_N$ is isomorphic to $\PS^2$ or is a
$\PS^1$-bundle over a curve $\Gamma$. 

If $S_N$ is a $\PS^1$-bundle over a curve $\Gamma$, as $h^1(S_N, \mathcal{O}_{S_N})=h^1(X_N,
\mathcal{O}_{X_N})=0$, $\Gamma$ has arithmetic genus
$h^1(\Gamma,\mathcal{O}_\Gamma)=0$ and $\Gamma$ is
isomorphic to $\PS^1$. 

Suppose that $S_N$ is
a $\PS^1$-bundle over $\PS^1$. The surface $S_N$ is a Hirzebruch
surface of the form $\F_a= \PS(\mathcal{O}_{\PS^1}
\oplus\mathcal{O}_{\PS^1}(a))$. Since $-K_{S_N}$ is nef, $S_N$ is either
$\F_0$ (i.e. $\PS^1\times \PS^1$) or $\F_2$. 
\end{proof}

\setcounter{step}{0}

\begin{thm}[End product of the MMP]
  \label{thm:4}
By the above, the end product of the Minimal Model Program $X_n$ is a
weak$\ast$ Fano $3$-fold. More precisely, $X_n$ is one of the
following:
\begin{enumerate}
\item $X_n$ is a terminal Gorenstein $\Q$-factorial Fano $3$-fold of
   rank $1$. $X_n$ is the deformation of a non-singular Fano $3$-fold of
   rank $1$.
\item $X_n \to \PS^1$ is a del Pezzo fibration of degree $k$, with $3 \leq k
  \leq 9$. $X_n$ has Picard rank
$2$. 
\item $X_n \to S$ is a conic bundle and $S$ is either $\PS^2$, $\PS^1
  \times \PS^1$ or $\F_2$. $X_n$ has Picard rank $2$ or $3$.
\end{enumerate}
\end{thm}
\begin{proof}
 Theorem~\ref{thm:3} shows that the category of weak$\ast$ Fano
 $3$-folds is stable under the birational operations of the Minimal
 Model Program.
Running the Minimal Model Program on a weak$\ast$ Fano $3$-fold $X_0$
 yields a sequence
 of extremal contractions and weak$\ast$ Fano $3$-folds:  
\[
X_0 \stackrel{\phi_0}\to X_1 \stackrel{\phi_1}\to \ldots
\stackrel{\phi_n}\to X_n.
\]
The Program terminates if $\rho(X_n)=1$ or if $X_n$ is a Mori fibre
space over a minimal base.
 
If $\rho(X_n)=1$, the anticanonical divisor is
ample. Indeed, assume $-K_X$ is nef but not ample. There exists an
effective curve $Z$ such that $-K_X \cdot Z=0$. Since $\rho(X_n)=1$,
$-K_{X_n}$ is trivial and therefore not big: this contadicts $X_n$
being weak$\ast$. 
In this case, Mukai shows that $X_n$
can be deformed to a non-singular Fano $3$-fold with Picard rank $1$ 
\cite{MR1944132}. Such Fano $3$-folds have been classified \cite{MR463151,MR503430}.

If $X_n \to \PS^1$ is a del Pezzo fibration, then
$\rho(X_n)=\rho(\PS^1)+1=2$.

If $X_n \to S$ is a conic bundle, by the discussion above, either $S$ is
$\PS^2$ and $\rho(X_n)=2$, or $S$ is $\PS^1 \times \PS^1$ or $\F_2$ and
$\rho(X_n)=3$.
\end{proof}

\section[A bound on the defect]{A bound on the defect of some Fano $3$-folds}
\label{sec:bound-defect-some-1}

The results obtained in Section~\ref{sec:categ-weak-fano} show that it
is possible to run a Minimal
Model Program (MMP) in the category of weak$\ast$ Fano
$3$-folds. If $X$ is a weak$\ast$ Fano $3$-fold and $Y$ is its
anticanonical model, the end product of
the MMP on $X$ is either a Fano $3$-fold, a del Pezzo fibration over $\PS^1$, or a conic
bundle over $\PS^2, \PS^1\times\PS^1$ or $\F_2$. Since the
defect of $Y$ is equal to $\rk \Pic(X)-\rk \Pic(Y)$, the number of
divisorial contractions needed to reach an
end product of the Minimal Model Program on a weak$\ast$ Fano $3$-fold $X$ 
determines the defect of its anticanonical model $Y$.

In this Section, I show that if the anticanonical model $Y$ of a weak$\ast$ Fano $3$-fold $X$ does
not contain a quadric $Q$ with $-K_{Y\vert Q}= \mathcal{O}_{Q}(1)$, a divisorial contraction increases the
anticanonical degree by at least $4$. I then prove that this condition
on $Y$ is preserved by the operations of the MMP. Let the weak$\ast$
Fano $3$-folds $X_i$, for $1\leq i\leq n$, be the intermediate steps of
the MMP on $X= X_0$ and let $Y_i$ be the anticanonical model of
$X_i$. If $Y=Y_0$ does not contain an irreducible reduced quadric $Q$
with $-K_{Y\vert Q}= \mathcal{O}_Q(1)$, then no $Y_i$ contains an
irreducible reduced quadric $Q$
with $-K_{Y_i\vert Q}= \mathcal{O}_Q(1)$. 
If $Y$ does contain an irreducible reduced quadric $Q$ with
$-K_{Y\vert Q}= \mathcal{O}_Q(1)$ however, there exists an extremal
divisorial contraction $\phi \colon X \to X_1$ whose exceptional
divisor is $\widetilde{Q}$, the pull back of $Q$ by the anticanonical
map. In this case, I show that there are at
most $10-g$ contractions of quadrics when running the MMP on $X$.  

Following these observations, I determine a bound on the defect of
terminal Gorenstein Fano $3$-folds $Y$ that contain neither a plane,
nor a quadric. I 
then give a general bound on the defect of terminal Gorenstein Fano
$3$-folds that do not contain a plane.

\subsection{Further study of the Minimal Model Program of weak$\ast$ Fano $3$-folds}
\label{sec:further-study-mmp}
 
\begin{lem} 
  \label{lem:5}
Let $X$ be a weak$\ast$ Fano $3$-fold. Assume that $Y$, the
anticanonical model of $X$,
does not contain an irreducible reduced quadric
$Q=\PS^1 \times \PS^1$ or $Q\subset \PS^3$ with ${-K_Y}_{\vert Q}= \mathcal{O}_Q(1)$.
Every divisorial extremal contraction $\phi$ increases the anticanonical degree
$(-K)^3$ by at least $4$.  
\end{lem}

\begin{proof}
  Let $\phi \colon X \to X'$ be a divisorial contraction. 
\setcounter{step}{0}
\begin{step}
The contraction $\phi$ is of type E$1$ or E$2$.  
\end{step}

Recall from the proof of Theorem~\ref{thm:3} that the contraction
$\phi$ is not of type E$5$.

If $\phi$ is of type E$3$ (resp. E$4$), $X$ contains a quadric $Q= \PS^1\times \PS^1$
  (resp. $Q \subset \PS^3$) with normal bundle
  $\mathcal{O}_{Q}(Q)= \mathcal{O}_{Q}(-1)= \mathcal{O}_{\PS^1 \times \PS^1}(-1,-1)$
  (resp. $\mathcal{O}_{Q}(Q)= \mathcal{O}_Q(-1)$). 
The divisor $Q$ is Cartier, hence, by adjunction ${-K_X}_{\vert Q}=
  \mathcal{O}_{\PS^1 \times \PS^1}(1,1)= \mathcal{O}_{Q}(1)$ (resp. ${-K_X}_{\vert Q}=
  \mathcal{O}_{Q}(1)$). 

 Since $X$ is a
  weak$\ast$ Fano $3$-fold, the linear system $\vert -K_X\vert_{\vert Q}$ is basepoint free
  (Theorem~\ref{thm:1}) and determines a morphism $\Phi_{\vert
  -K_X\vert_{\vert Q}}$ onto a surface.
 
As $\vert -K_X\vert_{\vert Q} \subset \vert \mathcal{O}_{Q}(1)\vert$, the morphism
  $\Phi_{\vert -K_X \vert_{\vert Q}}$ factorises as $\nu \circ \Phi_{\vert
  \mathcal{O}_Q(1) \vert}$, where $\nu$ is the projection from a
  possibly empty linear subspace
\[
\xymatrix{\PS(H^0(Q,\mathcal{O}_Q(1)))\simeq \PS^3 \ar@{-->}[r]&
  \PS(H^0(Q,{-K_X}_{\vert Q})) }
\] 
naturally associated to the inclusion of linear systems.
The image of $Q$ by $\Phi_{\vert -K_X\vert_{\vert Q}}$ is $\overline{Q}=\nu(Q)$. The map
$\nu$ is the identity or a projection from a point, line or plane. As
$\overline{Q}$ is a surface, $\nu$ can only be  the identity or a projection
from a point and $\overline{Q}$ either is the quadric $Q$ or is a plane $\PS^2$
and, by construction, ${-K_Y}_{\vert \overline{Q}}= \mathcal{O}_{\overline{Q}}(1)$.

In both cases, this contradicts the
assumption on $X$ since $Y$ contains neither a plane $\PS^2$ with
$-K_{Y \vert \PS^2}= \mathcal{O}_{\PS^2}(1)$ nor a quadric $Q$ with
$-K_{Y \vert Q}= \mathcal{O}_{Q}(1)$. The contraction $\phi$ cannot
therefore be of type E$3$ or E$4$; $\phi$ is of type E$1$ or E$2$.

\begin{step}
The anticanonical degree increases by at least $4$: $-K_{X'}^3
\geq -K_X^3+4$.  
\end{step}

If $\phi$ is of type E$2$, then $-K_X= \phi^{\ast}(-K_{X'})-2E$ and
the anticanonical degrees of $X$ and $X'$ satisfy $-K_X^3=-K_{X'}^3-8E^3=-K_{X'}^3-8$. The degree
increases by precisely $8$.
  
Assume that $\phi \colon X \to X'$
is an E$1$-contraction and denote by $\Gamma$ the centre of $\phi$.
Lemma~\ref{lem:32} shows that the following intersection table holds: 
\begin{align}
\label{eq:16}
(-K_X)^3= (-K_{X'})^3 - 2(-K_{X'} \cdot \Gamma -p_a(\Gamma)+1)\\
\label{eq:12}
(-K_X)^2 \cdot E=-K_{X'} \cdot \Gamma +2-2p_a(\Gamma)\\
\label{eq:13}
(-K_X)\cdot E^2= 2 p_a(\Gamma)-2\\
\label{eq:14}
E^3= -\deg\mathcal{N}_{\Gamma/X'}= 2-2p_a(\Gamma)- (-K_{X'}\cdot \Gamma) . 
\end{align}    

The anticanonical divisor $-K_X$ is  Cartier and nef. The
anticanonical map is small and determined by $\vert -K_X \vert$. The
anticanonical map $\Phi_{\vert -K_X \vert}$ contracts no divisor,
hence in particular, $(-K_{X})^2 \cdot E >0$. 
Equation \eqref{eq:12} implies that:
\begin{eqnarray*}
 -K_{X'}\cdot \Gamma \geq
2p_a(\Gamma)-1 \quad \mbox{that is:}\\
-K_{X'} \cdot \Gamma -p_a(\Gamma)+1\geq
p_a(\Gamma). 
\end{eqnarray*}
Equation \eqref{eq:16} shows that the required result holds for
$p_a(\Gamma)\geq 2$.

Assume that the centre of $\phi$ is a curve $\Gamma$ with arithmetic
genus $p_a(\Gamma) \leq 1$.
Since $\Gamma$ has planar singularities (Theorem~\ref{thm:8}), if its
arithmetic genus is $1$, its degree $-K_{X'}\cdot \Gamma$ is at least
$3$ and the desired inequality holds. 

Finally, assume that the curve $\Gamma$ has arithmetic genus $p_a(\Gamma)=0$,
that is that $\Gamma$ is rational and non-singular. Equation \eqref{eq:16} shows that
the inequality holds for $-K_{X'}\cdot \Gamma \geq 1$. 

\begin{cla}
The centre of $\phi$ is not a flopping curve. 
\end{cla}
Suppose that the centre of $\phi$ is a curve $\Gamma$ with $-K_{X'} \cdot
\Gamma =0$. 
From \eqref{eq:12}, $-K_X^2\cdot E>0$, and the curve $\Gamma$ is
rational and non-singular.  

The degree of $\mathcal{N}_{\Gamma/X'}$, the normal bundle of $\Gamma$ in $X'$,
 is $-2$. The curve $\Gamma \simeq \PS^1$ is locally a complete
 intersection, so that its normal bundle is 
\[
\mathcal{N}_{\Gamma/X'}=\mathcal{O}_{\PS^1}(n)\oplus
\mathcal{O}_{\PS^1}(-2-n)
\]
 for some integer $n\geq -1$.

Let $s$ be a section of $\mathcal{N}_{\Gamma/X'}$ corresponding to the
exact sequence: 
\[
0 \to  \mathcal{O}_{\PS^1}(2+n) \to \mathcal{N}_{\Gamma/X'} \to
\mathcal{O}_{\PS^1}(-n) \to 0.
\]  
The intersection of $s$ with the anticanonical divisor is $-K_{X'}
\cdot s=-n$. Since $-K_{X'}$ is nef, $n$ is either $0$ or
$-1$. The exceptional divisor $E$ of the contraction $\phi$ is
$\PS(\mathcal{N}_{\Gamma/X'})$, that is either $E\simeq\PS(\mathcal{O}_{\PS^1} \oplus
\mathcal{O}_{\PS^1}(-2))\simeq \F_2$ or $E\simeq\PS(\mathcal{O}_{\PS^1}(-1) \oplus
\mathcal{O}_{\PS^1}(-1))\simeq \PS^1\times\PS^1$. 

I show that these
cases are impossible unless $Y$ contains a quadric $Q$
with ${-K_Y}_{\vert Q}= \mathcal{O}_{Q}(1)$.

If $E\simeq \F_2$, the linear system $\vert -K_X \vert_{\vert E}$ is basepoint
free and determines a morphism $\Phi_{\vert -K_X \vert_{\vert E}}$ onto a
surface. The morphism $\Phi_{\vert -K_X \vert_{\vert E}}$ is determined by a linear system
$\vert -K_X \vert_{\vert E} \subset \vert af+b\sigma \vert$ on the
scroll. By convention, $f$ denotes the fibre of the scroll,
$\sigma$ its negative section, and the intersection numbers are $f^2=0$, $\sigma^2=-2$,
and $\sigma \cdot f=0$. Since $f$ is a fibre of the contraction $\phi$, 
$-K_X\cdot f=1$ and $a=1$. From \eqref{eq:12}, the degree is $(-K_X)^2\cdot E=
(\sigma+bf)^2=2$, so that $b=2$. Since $(\sigma+2f)\cdot \sigma=0$
and $(\sigma+2f)\cdot f=1$, $\vert \sigma +2f \vert$ determines a
morphism that contracts the negative section and that maps the scroll
onto an irreducible reduced singular quadric $\overline{Q}= \proj(\F_2, \vert\sigma+2f
\vert)\subset \PS^3=\PS(H^0(\F_2,\vert \sigma+2f
\vert ))$. The image of the morphism $\Phi_{\vert -K_X \vert_{\vert
    E}}$ is $\nu(\overline{Q})$, where $\nu$ is
the projection from a (possibly empty) linear subspace
\[
\xymatrix{\PS^3=\PS(H^0(\F_2,\vert \sigma+2f \vert ))\ar@{-->}[r]&
 \PS(H^0(E, {-K_X}_{\vert E}))} 
\]
associated to the inclusion $\vert -K_X \vert_{\vert E} \subset \vert \sigma+2f \vert$.

Similarly, if $E\simeq \PS^1\times\PS^1$, the linear system $\vert
-K_X\vert_{\vert E}$ is contained in the linear system
$\vert \mathcal{O}_{\PS^1\times \PS^1}(1,1)\vert$, because $-K_X \cdot l=1$ for any
ruling $l$ of $E$. The linear system $\vert \mathcal{O}_{\PS^1\times
  \PS^1}(1,1)\vert$ determines a morphism of $E$ onto an irreducible
reduced non-singular quadric
$\overline{Q}=\PS^1 \times \PS^1\subset \PS^3=\PS(H^0(E,\mathcal{O}_{\PS^1\times
  \PS^1}(1,1)))$.  The image of the morphism $\Phi_{\vert -K_X \vert_{\vert
    E}}$ is $\nu(\overline{Q})$, where $\nu$ is
the projection from a (possibly empty) linear subspace
\[
\xymatrix{\PS^3=\PS(H^0(E,\mathcal{O}_{\PS^1\times
  \PS^1}(1,1)))\ar@{-->}[r]&
  \PS(H^0(E, {-K_X}_{\vert E})} 
\]
associated to the inclusion $\vert -K_X \vert_{\vert E} \subset \vert
\mathcal{O}_{\PS^1\times\PS^1}(1,1)\vert$.

In both cases, the image of $\Phi_{\vert -K_X \vert_{\vert E}}$ is a
surface that is the image of an irreducible reduced quadric $\overline{Q}\subset \PS^3$
under a projection of $\PS^3$ from a linear subspace. As above,
$\Phi_{\vert -K_X \vert_{\vert E}}(E)$ is a plane
$\PS^2$ with $-K_{Y \vert \PS^2}= \mathcal{O}_{\PS^2}(1)$ or the
quadric $\overline{Q}$ with $-K_{Y \vert \overline{Q}}= \mathcal{O}_{\overline{Q}}(1)$. .
This yields a contradiction if $X$ is weak$\ast$ and $Y$ does not
contain a quadric $\overline{Q}$ with ${-K_{Y}}_{\vert \overline{Q}}= \mathcal{O}_{\overline{Q}}(1)$.
\end{proof}
\setcounter{step}{0}
I show that the hypotheses of Lemma~\ref{lem:5} are
preserved by the birational operations of the Minimal Model Program.
\begin{lem}
  \label{lem:26}
Let $X$ be a weak$\ast$ Fano $3$-fold and $\phi \colon X \to X'$ a
birational extremal contraction. Denote by $Y$ and $Y'$ the
anticanonical models of $X$ and $X'$. If $Y'$ contains an irreducible reduced quadric $Q'$ with
${-K_{Y'}}_{\vert {Q'}}= \mathcal{O}_{Q'}(1)$,
$Y$ also contains an irreducible reduced quadric $Q$ with ${-K_{Y}}_{\vert Q}= \mathcal{O}_{Q}(1)$.  
\end{lem}
\begin{proof} 
Denote by $g$ (resp. $g'$) the anticanonical map of $X$ (resp. $X'$). 
Assume that $Y'$ contains an
irreducible quadric $Q'$ with ${-K_{Y'}}_{\vert Q'}=
\mathcal{O}_{Q'}(1)$. Denote by $\widetilde{Q'}$ the proper transform of
$Q'$ by $g'$ and by $\widetilde{Q}$ the proper transform of $\widetilde{Q'}$ on $X$. 

The linear system $\vert -K_{X}\vert_{\vert \widetilde{Q}}= \vert
\phi^{\ast}(-K_{X'})-E \vert_{\vert \widetilde{Q}}$ is naturally a
subsystem of
$\vert (\phi)^{\ast}(-K_{X'}) \vert_{\vert \widetilde{Q'}}= \vert
(g'\circ \phi)^{\ast} \mathcal{O}_{Q'}(1)\vert$. 
The inclusion is strict if the centre of $\phi$ intersects
$\widetilde{Q}$. The $3$-fold $X$ is a weak$\ast$ Fano hence, by
Theorem~\ref{thm:1}, $\vert -K_{X}\vert_{\vert \widetilde{Q}}$ is
basepoint free and determines a morphism $\Phi_{\vert
  -K_{X}\vert_{\vert \widetilde{Q}}}$ that maps $\widetilde{Q}$ onto a
surface. 

Let $Q$ be the image of $\Phi_{\vert
  -K_{X}\vert_{\vert \widetilde{Q}}}$.
The morphism $\Phi_{\vert -K_{X}\vert_{\vert \widetilde{Q}}}$
is the composition $\nu \circ \Phi_{\vert (g'\circ
  \phi)^{\ast}\mathcal{O}_{Q'}(1)\vert}$ where $\nu$ is the projection
from a (possibly empty) linear subspace 
\[\xymatrix{
\PS(H^0(\widetilde{Q},(g'\circ \phi)^{\ast}\mathcal{O}_{Q'}(1)))\simeq
\PS^3 \ar@{-->}[r] & \ \PS(H^0(\widetilde{Q},{-K_X}_{\vert \widetilde{Q}}))
}\] 
associated to the inclusion $\vert {-K_X}_{\vert \widetilde{Q}}\vert
\subset \vert (g'\circ \phi)^{\ast}\mathcal{O}_{Q'}(1)\vert$. 

The surface $Q$ therefore is the image of $Q'\subset \PS^3$ under the
projection $\nu$ of $\PS^3$ from a linear subspace. As $Q$ is a
surface, the rational map
$\nu$ is either the identity, or the projection of $\PS^3$ from a
point. The surface $Q$ is either $Q'$ or a plane $\PS^2$ and
${-K_Y}_{\vert Q}= \mathcal{O}_Q(1)$ by construction. 
The surface $Q$ is not a plane because $X$ is a weak$\ast$ Fano
$3$-fold, hence $Q$ is an irreducible reduced quadric $Q\simeq Q'$. 
The proof shows in addition that $\widetilde{Q'} \subset X'$ does not intersect the centre of $\phi$. 
\end{proof}

\subsection{A bound on the defect of some Fano $3$-folds}
\label{sec:bound-defect-some}
 
Let $Y$ be a non $\Q$-factorial normal Gorenstein
terminal Fano $3$-fold of genus $g \geq 3$. The $3$-fold $Y$ is not
$\Q$-factorial, hence:
\[ 
\rk (\W (Y))= \dim H_4(Y,\Z) \neq \dim H^2(Y,\Z)=
\rk(\Pic(Y)).
\]
If $Y$ is a quartic
$3$-fold with terminal (Gorenstein) singularities, the
Grothendieck-Lefschetz hyperplane theorem shows that $\rho(Y)= \rk
\Pic(Y)=1$. In this section, I always assume that the Picard rank of $Y$ is $1$.

Recall the statement of Kawamata's
$\Q$-factorialisation theorem (Proposition~\ref{pro:4}).
\begin{pro}[$\Q$-factorialisation]
\label{pro:1}
Let $Y$ be an algebraic threefold with only terminal
singularities. Then, there is a birational morphism $f \colon X \to Y$
such that $X$ is terminal and $\Q$-factorial, $f$ is an
isomorphism in codimension $1$ and $f$ is projective.
\end{pro}
\begin{rem}
  \label{rem:7}
  The map $f \colon X \to Y$ is a chain of `symbolic blow-ups' of
  Weil non $\Q$-Cartier divisors on $Y$. The anticanonical map
  introduced in Section~\ref{sec:categ-weak-fano} is a
  $\Q$-factorialisation map.  
\end{rem}
\begin{rem}
  \label{rem:12}
As mentioned in Remark~\ref{rem:8}, the Fano indices of $X$ and $Y$
are equal. The map $f$ is crepant, and so the $3$-fold $X$ is also
Gorenstein.  
\end{rem}

The morphism $f$ is small and $X$ is $\Q$-factorial, hence:
\[
\rk\W(Y)=\rk\W(X)=\rk \Pic(X),
\]
and the defect of $Y$ is:
\[
\sigma(Y)= \rk \Pic(X)-1.
\]

The $\Q$-factorialisation $X$ of $Y$ is Gorenstein, terminal,
$\Q$-factorial and its
anticanonical divisor $-K_X= f^{\ast}(-K_Y)$ is nef and big. The
$3$-fold $X$ is a weak Fano $3$-fold.

The proof of
Theorem~\ref{thm:1} shows that if $\vert -K_X \vert$ is not basepoint
free, then either $X$ is birational to a special complete intersection
$X_{2,6}\subset \PS(1,1,1,2,3)$ and $X$ has genus $2$, or $X$ is
monogonal and $X$ has Picard
rank at least $2$. If $Y$ has Picard rank $1$ and genus $g \geq 3$,
$R(Y,-K_Y)=R(X, -K_X)$ is generated in degree $1$.
 
Unless $Y$ contains a plane $\PS^2$ with
${-K_Y}_{\vert \PS^2}=\mathcal{O}_{\PS^2}(1)$, its $\Q$-factorialisation $X$
is a weak$\ast$ Fano $3$-fold. 

I state here a weak version of Theorem~\ref{thm:16} \cite{MR1489117},
which is established in Section~\ref{sec:deformation-theory}. This
theorem is used in the rest of this section.
\begin{thm}
\label{thm:12}
  Let $Y$ be a Fano $3$-fold with terminal Gorenstein
  singularities. There is a $1$-parameter flat deformation of $Y$
  $f\colon \mathcal{Y}\to \Delta$, where $\mathcal{Y}_t$ is a terminal
  Gorenstein Fano $3$-fold for all $t\in \Delta$ and
  $\mathcal{Y}_{t_0}$ is
  non-singular for some $t_0 \in \Delta\smallsetminus
  \{0\}$.  
\end{thm}
\begin{rem}
If $\mathcal{Y} \to \Delta$ is such a
$1$-parameter flat deformation, $\rk \Pic(Y)=
\rk \Pic(\mathcal{Y}_t)$ for all $t\in\Delta$ (Lemma~\ref{lem:27}). 
The Fano index and the plurigenera are constant in a
$1$-parameter flat deformation, so that $-K_Y^3$, the degree of $Y$,
is equal to the degree of $\mathcal{Y}_{t_0}$ (see Section~\ref{sec:deformation-theory}).  
\end{rem}
\begin{cor}
  \label{cor:1} 
  {(of all the above):} Let $Y$ be a terminal Gorenstein Fano $3$-fold
  with Fano index $1$, Picard rank $1$ and genus $g \geq 3$.
If $Y$ does not contain a quadric or a plane, the defect of $Y$ is
  bounded by $\big[\frac{12-g}{2}\big]+4$, where $g$ is the genus of $Y$. 
\end{cor}

\begin{proof}
Denote by $X$ a small $\Q$-factorialisation of $Y$, $X$ is a
weak$\ast$ Fano $3$-fold. I prove that the Picard rank of $X$ is at
most $\big[\frac{12-g}{2}\big]+5$.

The anticanonical model $Y$ of $X$ is a terminal Gorenstein Fano
$3$-fold. The anticanonical degree of
$X$ is of the form $-K_Y^3=-K_X^3=2g-2$ with $2 \leq g \leq 10$ or $g=12$
by Iskovskih's classification \cite{MR463151,MR503430} and Theorem~\ref{thm:12}.

Lemma~\ref{lem:5} states that, when running the Minimal Model Program on $X$, each
divisorial contraction increases the anticanonical degree by at least $4$. By
Theorem~\ref{thm:12}, at each step of the MMP, the anticanonical model
$Y_i$ of $X_i$ is a terminal Gorenstein Fano $3$-fold and it can be
smoothed to a non-singular Fano $3$-fold with Picard rank $1$. 

 Iskovskih's classification of non-singular Fano $3$-folds shows that at each step,
$-K_{X_i}^3$ can only take one of finitely many values. If the Fano index
of $X_i$ is $1$, it is of the form $2g_i-2$, with $g+2 \leq g_i \leq
10$ or $g_i=12$. If $X_i$ has Fano index $2$, it is of the form $8d$, with $1 \leq d \leq 5$. If
$X_i$ is a possibly singular quadric, it is equal to $54$, and if $X_i
\simeq \PS^3$, it is equal to $64$ (see Theorem~\ref{thm:17}). 

In each case, the Picard rank of the weak$\ast$ Fano $3$-fold $X$
is equal to the sum of the number of divisorial contractions and the
Picard rank of the end product of the Minimal Model Program. The
end product of the MMP is one of: a terminal Gorenstein Fano
$3$-fold, a del Pezzo fibration over $\PS^1$,
or a conic bundle over $\PS^2$, $\PS^1 \times \PS^1$ or $\F_2$. 

If, at any intermediate step, $X_i$ is a
weak$\ast$ Fano $3$-fold of Fano index $2$, $\phi_i$ can only be an E$2$ contraction,
and $X_{i+1}$ also has Fano index $2$ (Lemma~\ref{lem:25}). If the end product $X_n$ has index
$2$, it is one of: a Fano index $2$ Fano
$3$-fold, an \'etale conic bundle, or it admits a fibration
by del Pezzo surfaces of degree $4$.
More precisely, depending on the end product of the Minimal Model
Program, one of the folllowing cases occurs:
\begin{enumerate}
\item{$Y_n$ is a Fano $3$-fold of Fano index $1$.} The genus decreases by at
  least $2$ with each divisorial contraction. The Picard rank of $X$
  is bounded above by $\big[ \frac{12-g}{2} \big]+1$.
\item{$Y_n$ is a Fano $3$-fold of Fano index $2$.} Similarly, by inspection,
  the Picard rank of $X$ is bounded by
  $\big[ \frac{12-g}{2} \big]+3$. 
\item{$Y_n$ is a possibly singular quadric.}
In this case, no intermediate $3$-fold can have index $2$ and the
Picard rank is bounded by $\big[\frac{12-g}{2}\big]+2$.
\item{$Y_n \simeq \PS^3$.}
No intermediate $3$-fold can have Fano index $2$, hence the Picard rank of
$X$ is bounded by $\big[\frac{12-g}{2}\big]+3$. 
\item{$Y_n$ is a strict Mori fibre space.}
In this case, as $Y_n$ is a terminal Gorenstein Fano $3$-fold with
Picard rank $1$, I can apply the same results to bound the degree of
$Y_n$ (refining the argument to take account of possible indices of Mori fibre spaces does
not improve the bound). This yields an upper bound on the 
Picard rank equal to $\big[\frac{12-g}{2}\big]+5$.
\end{enumerate}
\end{proof} 
\begin{cor} 
\label{cor:4} 
 Let $Y$ be a terminal Gorenstein Fano $3$-fold with Fano index $1$,
  Picard rank $1$ and genus $g\geq 3$. If $Y$ contains a
  quadric but does not contain a plane, the defect of $Y$ is at most $14-g$.
\end{cor}
\begin{proof}
Denote by $g \colon X \to Y$ a small $\Q$-factorialisation of $Y$; $X$ is a
weak$\ast$ Fano $3$-fold. I prove that the
Picard rank of $X$ is at most $15-g$.

I can run a Minimal Model Program on $X$. 
If $Q$ is a quadric lying on $Y$ with ${-K_Y}_{\vert Q}=
\mathcal{O}_Q(1)$, denote by $\widetilde{Q}$ its proper transform on
$X$. We claim that there is a $K$-negative extremal ray of $X$ on which
$\widetilde{Q}$ is negative, implying that there is an extremal contraction that
contracts $\widetilde{Q}$. Indeed, let $\overline{\Gamma}$ be any ruling of
$Q$ and let $\Gamma$ be the proper transform of $\Gamma$. As
$\widetilde{Q}$ is Cartier,
\[
\widetilde{Q}\cdot \Gamma= -K_X \cdot \Gamma+ K_{\Gamma}-(\Gamma \cdot
\Gamma)_{\widetilde{Q}}.
\]
In addition,  $-K_X \cdot \Gamma= -K_Y \cdot \overline{\Gamma}=1$ and
$-K_{\Gamma}= -K_{\overline{\Gamma}}=-2$ by the Leray spectral
sequence. The divisor $\overline{\Gamma}\subset Q$ is nef: $(\Gamma \cdot
\Gamma)_{\widetilde{Q}} \geq 0$.  Hence, $\widetilde{Q} \cdot \Gamma
<0$ and there exists a $K$-negative extremal ray $R$, on which $\widetilde{Q}$ is negative.
The contraction theorem \cite{MR946243} shows that $R$ can be
contracted. Denote by $\phi\colon X \to X'$ the extremal divisorial
contraction associated to $R$. 
\begin{cla}
The degree increases by at least $2$: $-K_X^3 \leq -K_{X'}^3 +2$.
\end{cla}
If $\phi$
contracts a quadric $\widetilde{Q}$ to a point,  $\phi$ is of type E$3$ or E$4$ and
the degree increases by $2$. 
If $\widetilde{Q}$ is contracted to a curve $\Gamma$,
\eqref{eq:12} shows that $-K_{X'} \cdot \Gamma= 2p_a(\Gamma)$ and
therefore by \eqref{eq:16}, $(-K_X)^3= (-K_{X'})^3 - 2(p_a(\Gamma)+1)$. 
The degree increases by at least $2$.

The classification of non-singular Fano $3$-folds of Picard rank $1$
shows that the number of quadrics lying on $Y$ is bounded by
$11-g=10-g+1$. 
As long as quadrics are contracted, the indices of the
intermediate steps of the Minimal Model Program are equal to $1$. The
bound on the Picard rank then follows from
Corollary~\ref{cor:4}.   
\end{proof}
\begin{rem}
I have not managed to construct a terminal Gorenstein quartic $3$-fold
$Y=Y_4^3\subset \PS^4$ that contains a
quadric  but no plane for which the bound
on the defect would be attained. My guess is that this bound is not
optimal.
First, if $Y$ contains two quadrics $Q_0=\{x_0=q_0=0\}$ and
$Q_1=\{x_1=q_1=0\}$, $Y$ naturally has a structure of del Pezzo
fibration of degree $4$ over $\PS^1$. One could try to
build an example by considering such del Pezzo fibrations with $3$
reducible fibres.
Second, if $Y$ does not contain two quadrics lying in distinct
hyperplane sections, $Y$ contains two
conjugate quadrics of the form $Q= \{x_0=q=0\}$ and $Q'=\{x_0=q'=0\}$
and $Y$ can be ``unprojected'' to a terminal Gorenstein
$Y_{2;3}\subset \PS^5$ with a node at $P= (0:0:0:0:0:1)$. One then
needs to consider quadrics in $Y_{2,3}$ that do not contain the point
$P$ (Lemma~\ref{lem:26}).   
\end{rem}
\begin{rem}
 \label{rem:34}
This analysis need not be limited to Gorenstein
terminal Fano $3$-folds of Picard rank $1$.
Let $Y$ be a non $\Q$-factorial terminal Gorenstein Fano $3$-fold of
Picard rank $\rho(Y)$ and denote by $X$ a small $\Q$-factorialisation of $Y$.
It is possible to run a Minimal Model
Program on $X$. Whereas Lemma~\ref{lem:12} shows that the Picard rank
of the anticanonical models of the intermediate steps of the MMP is
always $1$ when $\rho(Y)=1$, in the general case, at each step of the MMP,
the Picard rank of the anticanonical model
$Y_i$ is at most $\rho(Y)$.
Following the strategy of Corollaries~\ref{cor:1}
and~\ref{cor:4}, it is possible to bound the number of divisorial
contractions. This time, the end product is either a terminal Gorenstein Fano
$3$-fold of Picard rank at most $\rho(Y)$ or a strict Mori fibre space. The anticanonical
model of each intermediate
step can be deformed to a non-singular Fano $3$-fold; these have been
classified by Mori and Mukai
\cite{MR1102273,MR641971}.  
\end{rem}

\section{Fano $3$-folds containing a plane}
\label{sec:fano-3-folds}

In Sections~\ref{sec:categ-weak-fano} and \ref{sec:bound-defect-some-1},
I have determined an upper bound on the defect
of terminal Gorenstein Fano $3$-folds with Picard rank $1$ that do not
contain a plane.
The simplest examples of non $\Q$-factorial Fano $3$-folds are,
however, singular Fano $3$-folds $Y$ with Picard rank $1$ that contain
a plane. In this section, I study terminal Gorenstein Fano $3$-folds
with Picard rank $1$ that contain a plane.

Let $Y$ be a non $\Q$-factorial Fano $3$-fold with terminal Gorenstein
singularities and Picard rank $1$. Let $E\simeq \PS^2 \subset Y$ be a
plane contained in $Y$. The plane $E$ is a Weil non $\Q$-Cartier divisor. The 
anticanonical divisor of $Y$ is ample and $Y$ has Picard rank
$1$, hence $-K_Y= \mathcal{O}_Y(i(Y))$, where $i(Y)$ is the Fano index
of $Y$. In particular, if $E$ is any Weil non $\Q$-Cartier plane lying
on $Y$, $-K_{Y\vert E}= \mathcal{O}_E(i(Y))$. If $Y$ is
a non $\Q$-factorial Fano $3$-fold with terminal Gorenstein
singularities and Picard rank $1$, then the restriction of the
anticanonical divisor to any plane lying on $Y$ is determined by the Fano index of $Y$. 

In Section~\ref{sec:non-q-factorial-1}, I study terminal Gorenstein
Fano $3$-folds of Fano index strictly greater than $1$ that contain a
plane. The methods used in sections~\ref{sec:categ-weak-fano} and
\ref{sec:bound-defect-some-1} can be applied and yield a bound on the
defect.

In Section~\ref{sec:quartic-3-fold}, I focus on the case of terminal
quartic $3$-folds $Y= Y_4^3 \subset \PS^4$ that contain a
plane. Notice that such a quartic $3$-fold $Y= Y_4^3 \subset \PS^4$
has Fano index $1$ and $Y$ is Gorenstein because it is a hypersurface in
$\PS^4$. Section~\ref{sec:mixed-hodge-theory} explains that the defect
of a terminal Gorenstein Fano $3$-fold depends on the number of
singular points and on the position of its singularities. If $Y= Y_4^3
\subset \PS^4$ is nodal, that is if $Y$ has no worse than ordinary
double points, $Y$ has at most $45$ nodes
\cite{MR712934,MR848512}. Further, it is known that, up to projective
equivalence, there is a unique quartic hypersurface with $45$ nodes
\cite{MR1032936}: the Burkhardt quartic (Example~\ref{exa:1}).    
This quartic contains forty planes and has
defect $15$. I show in Section~\ref{sec:quartic-3-fold} that this is
an upper bound for the defect of terminal quartic $3$-folds. 
Finally, I make a conjecture (\ref{con:1}) on the defect of terminal Gorenstein
Fano $3$-folds with Picard rank $1$ and Fano index $1$ that contain a plane.

\subsection[Higher index Fano $3$-folds]{Non $\Q$-factorial Fano $3$-folds of higher index}
\label{sec:non-q-factorial-1}

Higher Fano index canonical Gorenstein Fano $3$-folds are classified by
Fujita and Shin. Shin proves the following: 

\begin{thm}[\cite{MR1030501}]
\label{thm:17}
Let $Y$ be a Gorenstein Fano $3$-fold with at most canonical
singularities. Denote by $i(Y)$ the Fano index of $Y$, then:
\begin{enumerate}
\item $i(Y) \leq 4$ with equality if and only if $Y \simeq \PS^3$.
\item $i(Y)=3$ if and only if $Y$ is a possibly singular quadric in $\PS^4$.
\end{enumerate}
\end{thm}

Let $Y$ be a terminal Gorenstein Fano $3$ fold with Fano index
$i(Y)>1$  and Picard rank $1$, and let $X$ be a small
$\Q$-factorialisation of $Y$. Recall that $i(X)=i(Y)$ by
Remark~\ref{rem:8}. The $3$-fold $X$ is a weak$\ast$
Fano $3$-fold. Indeed, $Y$ cannot
contain a plane with $-K_{Y \vert \PS^2}=
\mathcal{O}_{\PS^2}(1)$, the anticanonical ring of $Y$ is generated in
degree $1$ (Theorem~\ref{thm:1}) and the anticanonical map is small.
The results obtained in Section~\ref{sec:categ-weak-fano} show that it is
possible to run a Minimal Model Program on $X$.

Let $Y$ be a terminal Gorenstein Fano $3$-fold of Fano index
$2$ and let $X$ be a small $\Q$-factorialisation of $Y$. I state in
the following Lemma some straightforward results on
contractions of extremal rays on $X$. 
 
\begin{lem}
\label{lem:25}
Let $X$ be a Fano index $2$ weak$\ast$ Fano $3$-fold and $f\colon
X \to X'$ the contraction of an extremal ray. Then, one of the
following holds:
\begin{enumerate}
\item $f$ is birational, and $f$ is either a flop or an E$2$ contraction. 
\item $f\colon X\to S$ is a conic bundle, in fact $f$ is an \'etale $\PS^1$-bundle.
\item $f\colon X \to \PS^1$ is a del Pezzo fibration of degree $8$,
  and $f$ is a quadric bundle.  
\end{enumerate}
\end{lem}

\begin{rem}
\label{rem:32}
 Similarly, if $X$ has Fano index $3$ and $\rho(X)\geq 2$, then $\rho(X)=2$.
 Any contraction of an extremal ray is either a flop or a del
 Pezzo fibration of degree $9$, that is, $X$ is a $\PS^2$-bundle over
 $\PS^1$.  
\end{rem}
\begin{proof}
Let $f \colon X \to X'$ be the birational contraction of an extremal
ray. Since $X$ is weak$\ast$, $f$ cannot be of type E$5$. If $f$ is of
type E$1$, E$3$ or E$4$, an irreducible reduced curve $C$ contracted by $f$
has anticanonical degree $-K_X \cdot C=1$: this contradicts $i(X)=2$. 

If $f$ is a birational contraction, it is either a flopping contraction or the inverse of the blow up of a
smooth point. The exceptional divisor of $f$ is then a plane $E \simeq
\PS^2$ with $\mathcal{O}_E(-E)= \mathcal{O}_E(1)$, that is, by
adjunction, with ${-K_{X}}_{\vert E}= \mathcal{O}_{\PS^2}(2)$. 

Let $f$ be a fibering contraction. Unless $f$ is
an \'etale $\PS^1$-bundle over $S=\PS^2,\PS^1\times \PS^1$ or $\F_2$,
or a del Pezzo fibration of degree $8$, there exists a curve $C$
contracted by $f$ with anticanonical degree $-K_X\cdot C=1$ or $3$.   
\end{proof}

\begin{rem}
  \label{rem:14}
In particular, if $\rk \Pic(X)> 2$, $X$ has to contain a
plane.
The problem of bounding the defect of a terminal Gorenstein
Fano $3$-fold $Y$ of Fano index $2$ is equivalent to that of determining the
maximal number of disjoint planes that can lie on $X$ (see the proof
of Corollary~\ref{cor:3}).    
\end{rem}

\begin{lem}
  \label{lem:11}
  Let $f \colon X \to X'$ be an extremal contraction of type E$2$. Assume
  that $X$ is a weak$\ast$ Fano $3$-fold of index $2$ and degree $h^3=
  \frac{-K_{X}^3}{8}$. Then $X'$ is a weak$\ast$ Fano $3$-fold of
 Fano index $2$ and degree $h'^3=h^3-1$.
\end{lem}
\begin{proof}
By Theorem~\ref{thm:3}, $X'$ is a weak$\ast$ Fano $3$-fold. The
anticanonical divisors of $X$ and $X'$ satisfy the relation:
\[
-K_{X}= f^{\ast}(-K_{X'})-2E,
\] 
and $-K_{X}^3=-K_{X'}^3-8$. This also shows that $2$ divides the
Fano index of $X'$, so that $i(X')=2$ or $i(X')=4$.

 The $3$-fold $X'$ cannot have Fano index $4$. If it did,
its anticanonical model $Y'$ would be $\PS^3$
 (Theorem~\ref{thm:17}).  By Theorem~\ref{thm:12}, $Y$ is a one-parameter deformation of a
 non-singular Fano $3$-fold $Y_t$ with Fano index $2$ and Picard rank $1$. 
The degree is invariant
 in a one-parameter flat deformation, hence $-K_Y^3= -K_X^3
 \leq 40$.
In particular, $X'$ cannot have degree $64=-K_{\PS^3}^3$; thus
$X'$ has Fano index $2$.
\end{proof}

Let $Y$ be a non $\Q$-factorial Gorenstein Fano $3$-fold with Picard
rank $1$ and Fano index $2$. Let $X$ be a $\Q$-factorialisation of $Y$. The $3$-fold $X$
is a weak$\ast$ Fano $3$-fold. I run a Minimal Model Program on $X$. 
\[
X=X_0 \stackrel{f_1}\to X_1 \stackrel{f_2} \to X_2 \to \ldots
\stackrel{f_n}\to X_n
\]
By Lemma~\ref{lem:25}, each birational contraction is either a flop or
of type E$2$. Theorem~\ref{thm:4} gives the following possible end
products of the Minimal Model Program on
$X$:
\begin{enumerate}
\item $X_n$ is a terminal Gorenstein $\Q$-factorial Fano $3$-fold of
  rank $1$ and Fano index $2$. 
\item  $X_n \to \PS^1$ has a structure of del Pezzo
  fibration of degree $8$ (the general fibre is isomorphic to $\PS^1
  \times \PS^1$).
\item $X_n \to S$ has a structure of \'etale $\PS^1$-bundle over $S= \PS^2,
  \PS^1 \times \PS^1$ or $S=\F_2$.
\end{enumerate}
The Picard number of $X$ is determined by the number of planes
contracted while running the Minimal Model Program, that is by the number of E$2$ contractions.
\begin{cor}
  \label{cor:2}
 Let $Y$ be a non $\Q$-factorial terminal Gorenstein Fano $3$-fold of Picard
rank $1$ and Fano index $2$, and  let $X$ be a small $\Q$-factorialisation of
$Y$. Denote by $h^3= \frac{-K_X^3}{8}$ the degree of $X$ and $Y$. 
The Picard rank of $X$ is at most $8-h^3$.
\end{cor}
\begin{proof}
By Lemma~\ref{lem:11}, each divisorial contraction increases the
degree by $1$. Denote by $l$ the number of divisorial contractions
encountered in running
the Minimal Model Program on $X$.
\begin{enumerate}
\item $X_n$ is a terminal Gorenstein
  $\Q$-factorial $3$-fold of Fano index $2$. By Theorem~\ref{thm:12}, $X_{n}$ has degree  $1 \leq
  h_{n}^3 \leq 5$. As
  $h_{n}^3=h^3+l$ and $\rho (X)=l+1$, $\rho (X) \leq 6-h^3$.
\item $X_n \to \PS^1$ is a del Pezzo
  fibration of degree $8$. By Theorem~\ref{thm:12}, $Y_{n}$ can be deformed to a
  non-singular Fano $3$-fold of Fano index $2$ and Picard rank $1$. In
  particular, the degree $h_{n}^3 \leq 5$. Since $h_{n}^3=h^3+l $ and $\rho
  (X)=l+2$, $\rho (X) \leq 7-h^3$.
\item $X_n \to \PS^2$ (resp. $X_n  \to \PS^1 \times \PS^1$ or $\F_2$ ) is an \'etale
  $\PS^1$-bundle, and as above $X_{n}$ has degree
  $h_{n}^3 \leq 5$. Since $h_{n}^3=h^3+l $ and $\rho (X)=l+2$
  (resp. $\rho (X)=l+3$), $\rho (X) \leq 7-h^3$ (resp. $\rho(X) \leq 8-h^3$).
\end{enumerate}
\end{proof}
\begin{cor}
  \label{cor:3}
   Let $Y$ be a non $\Q$-factorial terminal Gorenstein Fano $3$-fold of Picard
rank $1$ and Fano index $2$ and let $X$ be a small
$\Q$-factorialisation of $Y$. Then $X$ contains at most $7-h^3$ disjoint planes.
\end{cor}
\begin{proof}
 Let $\overline{E}= \PS^2$ be a plane contained in $X$, then $-K_{X \vert \PS^2}=
 \mathcal{O}_{\PS^2}(2)$. 

By the contraction theorem, there is an extremal contraction which
contracts $E$ to a point. Indeed, $E$ is covered by curves $\Gamma$
with $-K_X \cdot \Gamma=2$ and $E \cdot \Gamma=-1$. 
 
\begin{cla}
The planes contracted when running the Minimal Model
 Program on $X$ are disjoint non $\Q$-Cartier planes on $Y$.    
\end{cla}
 
I show that I may assume that the flopping
contractions are performed first, that is, that the Minimal Model Program on $X$ is of
the form:
\[
X \to X_1 \to \cdots \to X_{n-l-1} \to \cdots \to X_{n}
\]
where all the extremal contractions $f_i$ for $i \leq n-l-1$ are
flopping contractions and the  extremal contractions $f_i$ for $i \geq
n-l$ are divisorial contractions of type E$2$. 

A flopping contraction and a divisorial contraction of type
E$2$ always commute. Indeed, the exceptional locus of a flopping
contraction is $K$-trivial,
while that of an E$2$ contraction is $K$-negative. 
More precisely, assume that an E$2$ contraction $f_{i-1}$ is followed
by a flopping contraction $f_i$:
\[ X_{i-1}\stackrel{f_{i-1}}\to X_i \stackrel{f_i}\to X_{i+1} \]
Then the centre $P$ of $f_{i-1}$
does not belong to a flopping curve $C$.
Otherwise, let $\widetilde{C}$ be a curve, which maps $1$-to-$1$ to
$C$; $-K_{X_{i-1}} \cdot \widetilde{C}<0$ and this contradicts
$X_{i-1}$ being weak$\ast$.

Let $E_i$ be the exceptional divisor of $f_i$ for $i \geq n-l$.
Without loss of generality, I want to prove that for $i \geq n-l+1$,
the centre $P_i$ of $f_{i}$ on $X_i$ does not belong to
$E_{i+1}$. The proper transform of any line through
$P_i$ is a flopping curve in $X_{i-1}$. The morphism determined by
$\vert -K_{X_{i-1}}\vert$ contracts finitely many curves, in
particular, there are finitely many curves $C$ through $P_i$ such that $-K_{X_i}
\cdot C=2$. The point $P_i$ cannot lie on $E_{i+1}$.

In particular, the proper transforms on $X$ of the 
divisors contracted when running the Minimal Model Program on $X$ are
disjoint planes. There is no other plane disjoint from these on $X$,
as this would give rise to another extremal contraction.
\end{proof}
\begin{rem}
The image of a plane $E \subset X$ on $Y$ is $\nu(S)$, where $S \simeq
\PS^2\subset \PS^5= \PS(H^0(\PS^2, \mathcal{O}_{\PS^2}(2)))$ is the
Veronese surface and $\nu$ is the projection of $\PS^5$ from a
possibly empty linear subspace
\[
\xymatrix{\PS^5= \PS(H^0(\PS^2, \mathcal{O}_{\PS^2}(2))) \ar@{-->}[r]
  &  \PS(H^0(\PS^2,-K_{X \vert \PS^2 })} 
\]
associated to the inclusion of linear subspaces $\vert -K_{X \vert
  \PS^2} \vert \subset  \vert \mathcal{O}_{\PS^2}(2) \vert$. As $X$ is
weak$\ast$, $\vert -K_X \vert_{\vert E}$ is basepoint free and $\nu$
is the identity or the projection of $S$ from a point $P\in
  \PS^5\smallsetminus S$.   
\end{rem}

\subsection{Quartic $3$-fold containing a plane}
\label{sec:quartic-3-fold}
Assume that $Y=Y_4^3 \subset \PS^4$ is a quartic
$3$-fold with terminal singularities that contains a plane $\Pi =
\{x_0=x_1=0\}$. As $Y$ is a complete intersection in $\PS^4$, $Y$ is
Gorenstein. The Grothendieck-Lefschetz hyperplane theorem shows that
$\rho(Y)=1$.
 I aim at
determining a bound on the defect of $Y$.
The equation of $Y$ is of the form: 
\begin{equation}
  \label{eq:35}
Y= \{x_0 \cdot a_3(x_0, \cdots , x_4)+
x_1 \cdot b_3(x_0,\cdots, x_4)=0 \}.  
\end{equation}

Let $X$ be the blow up of $Y$ along the plane $\Pi$. Locally, the
equation of $X$ can be written as:
\begin{multline}
  \label{eq:33}
\{ t_0 \cdot a_3(t_0x, t_1x, x_2, x_3, x_4)+ t_1 \cdot b_3(t_0x, t_1x , x_2, x_3, x_4)=0 
\} \\\subset \PS (t_0, t_1) \times \PS (x, x_2, x_3, x_4),  
\end{multline}
where the variable $x$ is defined by $x_0=t_0 \cdot x$, $x_1=t_1 \cdot x$. 

I fix some notation for rational scrolls over
$\PS^1$. Let $t_0, t_1$ be coordinates on $\PS^1$ and consider
$\C^{n+1}$ with coordinates $x_0, \ldots , x_n$. Fix a set of non negative
integers $a_0, \ldots , a_n$ in increasing order.

 I consider actions of the group $G=\C^{\times} \times \C^{\times}$
on the affine space $\C^2 \times \C^{n+1}$, where the two factors of
$G$ act by:
\begin{eqnarray*}
  \label{eq:3}
  \lambda : (t_0, t_1,x_0, \ldots , x_n ) \mapsto  (\lambda
  t_0,\lambda t_1,\lambda^{-a_0}x_0, \ldots ,\lambda^{-a_n} x_n )\\
\mu : (t_0, t_1,x_0, \ldots , x_n ) \mapsto (t_0, t_1,\mu x_0, \ldots ,\mu x_n )
\end{eqnarray*}
I use the following matrix notation to summarise this action: 
\[
\left( 
  \begin{array}{ccccc}
1 & 1 & a_0 & \ldots & a_n \\
0 & 0 & 1 & \ldots & 1
  \end{array}
\right)
\]
The scroll $\F(a_0, \ldots , a_n)$ is the quotient
$(\C^2 \smallsetminus\{0\})\times (\C^{n+1} \smallsetminus \{0\})/G$. The scroll $\F$ is a
$\PS^n$-bundle over $\PS^1$.
In the case of $\PS^3$-bundles over $\PS^1$, I denote by
$\F_{a_1,a_2,a_3}$ the scroll $\F(0,a_1,a_2, a_3)$.
Line bundles on $\F$ are in $1$-to-$1$ correspondence with characters
$\chi \colon G \to \C^{\ast}$ of the group $G$. Indeed, associate
to a character $\chi$ the line bundle $L_{\chi}$ such that the space of
sections of $L_{\chi}$ is: 
\[
H^0(\F, L_{\chi})= \{ f:(\C^2-\{0\})\times (\C^{n+1}-\{0\}) \to \C
\mid f(gx)= \chi(g)f(x) \}.
\]
 Let $L$ and $M$ be the line bundles associated to the characters
 $(\lambda, \mu) \mapsto \lambda$ and $(\lambda, \mu) \mapsto \mu$
 respectively. Denoting by $\pi \colon \F \to \PS^1$ the natural
 projection, $L$ is the pull back $\pi^{\ast} \mathcal{O}_{\PS^1}(1)$.

The $3$-fold $X$ is naturally
embedded in the scroll $\F_{0,0,1}$ over $\PS^1$. It is an element of
the linear system $\vert 3M+L\vert $ on $\F$, and it has a natural
structure of cubic del Pezzo fibration
over $\PS^1$. The generic fibre $X_{\eta}$ is reduced and irreducible
because $Y$ is. However,
special fibres might be reducible. Moreover, $X$ is a weak Fano $3$-fold,
$\rho(X/Y)=1$ and the map $X \to Y$ is small
($X$ is a small partial $\Q$-factorialisation of $Y$).

Bounding the rank of the group of Weil divisors of
$X$ suffices to determine the defect of $Y$.
 
\begin{lem}\cite{MR1426888}
\label{lem:15}
Let $X_i$, $0 \leq i \leq n$, be the reducible fibres of the cubic
fibration $X$. Let $D_{\eta, k}$ be irreducible generators of $NE(X_{\eta})$, $D_k$
the closure of $D_{\eta, k}$ in $X$ (each $D_k$ is
irreducible) and let the $E_{i,j}$ be the irreducible components of $X_i$.
The cone of effective divisors $N^1(X)$ is generated by the
$E_{i,j}$ and by the divisors $D_k$.
\end{lem}
\begin{rem}
  \label{rem:18}
A similar description holds for any del Pezzo fibration. 
\end{rem}
Corti proves in \cite{MR1426888} that, when projecting a terminal
weak Fano
del Pezzo fibration away from a plane contained in a reducible fibre, the
number of irreducible components of this reducible fibre
decreases. Let me recall his results more precisely.

Let $\mathcal{O}$ be a discrete valuation ring, with fraction field
$K$, parameter $t$ and residue field $k$. In the cases I consider,
$\mathcal{O}$ is the local ring of $\PS^1$ at a point $P$ that
corresponds to a reducible fibre and $k$ is
$\C$. Let $S = \Spec(\mathcal{O})$, $\eta= \Spec K$ be the generic
point and $0= \Spec k$.  
Let $\PS= \PS^n_{\mathcal{O}}$ be an $n$-dimensional projective
space over $S$ and $L=L_d$ a $d$-dimensional projective subspace,
defined over $k$. If $d \leq n-1 $, there is a birational
transformation $\Phi=\Phi_L \colon \PS \to \PS$, centred at $L$. The
map $\Phi$ is a projection from $L$. In homogeneous coordinates, it is
given by:
\[
\Phi \colon (x_0 : \ldots :x_n) \mapsto (tx_0: \ldots :tx_d:x_{d+1}:
\ldots :x_n).
\] 
Now, let $X \subset \PS$ be a cubic surface defined over $S$, such
that $X_{\eta}=X_K$ is smooth. I assume that $X$ is a
weak Fano $3$-fold in the sense of Definition~\ref{dfn:1} and that the
anticanonical map of $X$ is small. Let $n(X)$ denote the
number of $k$-irreducible components of the central fibre $X_0$.  
\begin{lem}\cite[Lemma 2.17]{MR1426888}
\label{lem:13}
 Assume that $X_0$ contains a $2$-plane $\Pi \subset \PS_0$, defined
 over $k$. Let $\Phi \colon \PS \dasharrow \PS^+$ be the projection
 away from $\Pi$ and $X^+= \Phi_{\ast} X$. Then:
 \begin{enumerate}
\item{}$X^+$ has terminal singularities, $-K_{X^+}$ is nef and big over $\Spec(\mathcal{O})$. 
\item{} $n(X^+) < n(X)$.
 \end{enumerate}
\end{lem}
\begin{rem}
  \label{rem:20}
Notice that in such an operation of projection away from a plane, the
Weil rank decreases by at most two. Indeed, this follows from
Lemma~\ref{lem:15}, as such a transformation only
affects $X$ in the central fibre.
\end{rem}

\begin{cor}[\cite{MR1426888}, Flowchart 2.13 and Corollary 2.20]
\label{cor:5}
If $X \subset \PS_{\mathcal{O}}^3$ is a weak Fano cubic fibration with small anticanonical map, then
projecting away from planes contained in the central fibre yields a
standard integral model $X'$ for $X_{\eta}$, that is a flat subscheme
$ X' \subset \PS_{\mathcal{O}}^3$ with isolated cDV singularities and
reduced and irreducible central fibre.   
\end{cor}

\begin{lem}
\label{lem:8}
The rank of the Weil group of $X$ is bounded by $8+2\cdot N + M$,
where $N$ is the number of reducible fibres with $3$ irreducible
components and $M$ the number of reducible fibres with $2$ irreducible
components. The defect of $Y$ is hence bounded above by $7 + 2\cdot N
+ M$. 
\end{lem}
\begin{proof}
If a fibre of $X \to \PS^1$ is reducible, it has at most $3$ irreducible
components (one of which has to be a plane). Project away from
these reducible components in order to obtain a standard integral
model. The result follows from the following exact sequence, valid
when $\pi \colon X' \to \PS^1$ has reduced and irreducible fibres:
\[
0 \to \Z[\pi^{\ast} \mathcal{O}_{\PS^1}(1)] \to \W(X')\to \Pic
(X'_{\eta})\to 0.
\]
\end{proof}
The problem is now to bound the possible number of reducible
fibres of $X$. 

\begin{exa}[The Burkhardt quartic]
\label{exa:1}  
It is known that a nodal quartic hypersurface in $\PS^4$ has at most $45$ nodes
\cite{MR712934,MR848512}. Moreover, up to projective equivalence, there
is only one such $Y= Y^3_4 \subset \PS^4$ \cite{MR1032936} and it has the following equation:
\[ 
 \{x_0^4 -x_0 (x_1^3+x_2^3+x_3^3+x_4^3)+ 3x_1x_2x_3x_4 = 0 \}. 
\]
Let $\widetilde{Y}$ be the blow up of $Y$ at the $45$ nodes.
The Hodge theoretic approach (Section~\ref{sec:mixed-hodge-theory},(\ref{eq:47}))
leads to $\sigma(Y)=b_2(\widetilde{Y})-45-b_2(Y)=
b_2(\widetilde{Y})-46$.
The cohomology of $\widetilde{Y}$ is
determined in \cite{MR1828606} and $b_2(\widetilde{Y})=61$, so that the defect of the Burkhardt
quartic is $15$.

Alternatively, the plane $\Pi=\{x_0=x_1=0\}$ is contained in $Y$. Let
$X$ be the $3$-fold obtained by blowing up $Y$ along the plane $\Pi$;
$X$ has the following expression:
\[
t_0(t_0^3-t_1^3)x^3-t_0(x_2^3+x_3^3+x_4^3)+ 3t_1x_2x_3x_4=0,
\] 
which on the affine piece $t_1 \neq 0$ reads
\[
t_0(t_0^3-1)x^3-t_0(x_2^3+x_3^3+x_4^3)+ 3x_2x_3x_4=0.
\] 
The central fibre, which corresponds to  $(t_0:t_1)=(0:1)$, has three
irreducible components: the planes $\{ t_0=x_i=0\}$.
Each fibre reduced and irreducible for $(t_0:t_1) \neq (0:1)$ and
$(t_0:t_1) \neq (1:\omega^i)$ for $0 \leq i \leq 2$, where $\omega$ is a cube root
of unity (notice that considering the affine piece $t \neq 1$ yields
the same results). The generic fibre $X_{\eta}$ is a non-singular
cubic surface in $\PS^3$. 

Consider for instance the fibre $X_1$ over $(1:1)$. The fibre $X_1$ is the
union of three planes in $\PS(x, x_2, x_3, x_4)$, namely $\Pi_1=\{
x_2+x_3+x_4=0\}$, $\Pi_2=\{\omega x_2+x_3+\omega^2 x_4=0\}$ and
$\Pi^3=\{\omega^2 x_2+x_3+\omega x_4=0\}$. The situation is
analogous for the other two reducible fibres.
There are $27$ closed subschemes in
$X \smallsetminus\bigcup_{i=0,1,2}
X_i$ isomorphic to $\PS^1_{\PS^1 \smallsetminus
  (\cup_i\{(\omega^i:1)\})}$. In other words, the $27$
lines on the generic fibre may be completed to divisors on $X$, they
are rational over $\PS^1 \smallsetminus (\cup_i\{(\omega^i:1)\})$. The
Picard rank of $X_{\eta}$ is $7$ and the generators of
$\Pic(X_{\eta})$ complete to independent Cartier divisors on $X$. 

This shows that the rank of the Weil group of $X$ is $1+7+8$.
The rank of the Weil group of the  Burkhardt quartic is $16$ and its
defect is $15$. 
\end{exa}

I show that if $X$ is a cubic fibration obtained by blowing up a
plane lying on a quartic $3$-fold $Y=Y_4^3 \subset \PS^4$ with
terminal singularities, then $X$ has at most $4$ reducible
fibres. By construction, $Y$ is the anticanonical model of $X$.

\begin{lem}
  \label{lem:14}
Let $Y =Y_4^3\subset \PS^4$ be a terminal quartic $3$-fold that contains a plane
$\Pi$ and let $X$ be the cubic del Pezzo fibration obtained by
blowing up $Y$ along $\Pi$. The cubic fibration $X$ has at most $4$
reducible fibres. 
\end{lem}
\begin{cor}
\label{cor:6}
The defect of a terminal quartic $3$-fold $Y=Y_4^3 \subset \PS^4$
 that contains a plane is at most $15$.  
\end{cor}
\begin{proof}[Proof of Lemma~\ref{lem:14}]
 \setcounter{step}{0}
Recall that $X$ is the blow up of a terminal quartic
$3$-fold $Y \subset \PS^4$ of the form  
\begin{equation}
  \label{eq:37}
 Y=\{x_0a_3(x_0, \cdots , x_4)+ x_1b_3(x_0, \cdots , x_4)=0\} \subset \PS^4
\end{equation}
along the plane $\Pi= \{x_0=x_1=0\}$. 

From \eqref{eq:33}, the equation of $X$ reads:
\[
X=\{t_0 a_3(t_0x, t_1x, x_2,x_3, x_4) +t_1 b_3(t_0x, t_1x, x_2,x_3, x_4)=0\},  \]
where $a_3$ and $b_3$ are sections of the linear system $\vert 3M \vert$ on
the scroll $\F_{0,0,1}$. Denote by $f$ the morphism $f \colon X\to
Y$. 

The terminal Gorenstein $3$-fold $X$ is a section of the linear system $\vert
 3M+L \vert $ on the scroll $\F_{0,0,1}$.

A change of coordinates on the scroll $\F_{0,0,1}$ is
given by a usual change of coordinates on $\PS^1$ and by the action
of a matrix
of the following form on $(x_2:x_3:x_4:x)$:
\[
\left( 
  \begin{array}{cccc}
 a & b & c & l_1(t_0, t_1) \\
d&e&f&l_2(t_0, t_1)\\
g&h&i& l_3(t_0, t_1)\\
0&0&0& j
  \end{array}
\right),
\]
where the upper $3$ by $3$ minor is an invertible complex matrix,
$l_1, l_2$ and $l_3$ are linear forms and $j$ is a non-zero complex number.   

The fibre of the scroll $\F_{0,0,1}$ over $(\lambda:\mu)\in \PS^1$ is
isomorphic to the hyperplane $H_{(\lambda:\mu)}= \{\mu x_0-\lambda
x_1=0\}\subset \PS^4$. There is a $1$-to-$1$ correspondence
between fibres of the scroll $\F_{0,0,1}$ and hyperplanes of $\PS^4$
that contain the plane $\Pi$.

The hyperplane section of $Y\subset \PS^4$ corresponding to
$H_{(\lambda: \mu)}$ is a reducible quartic surface that contains the
plane $\Pi$. The intersection $Y \cap H_{(\lambda:\mu)}$ is $\Pi \cup
Y'_{(\lambda: \mu)}$, where $Y'_{(\lambda: \mu)}$ is a possibly
reducible cubic surface, naturally isomorphic to the fibre
$X_{(\lambda: \mu)}$.  

If $X$ has a reducible fibre over $(\lambda:\mu)\in \PS^1$, $X$ contains a plane either of the form: 
\[
\Pi_{(\lambda:\mu)}=\{\lambda t_0 +\mu t_1= l(x_2, x_3, x_4)+l'(t_0x, t_1x)=0 \},
\]
where $l$ and $l'$ are linear, or of the form:
\[
\Pi_{(\lambda:\mu)}=\{ \lambda t_0 +\mu t_1=x=0 \}.
\]
The plane $\Pi_{(\lambda: \mu)}\subset X_{(\lambda:\mu)}$ is isomorphic to a plane $\Pi'_{(\lambda:
  \mu)}$ lying on $Y\cap H_{(\lambda:\mu)}$. In the first case,
  $\Pi'_{(\lambda:\mu)}$ is of the form:
\[
\Pi'_{(\lambda:\mu)}=\{\lambda x_0 +\mu x_1= l(x_2, x_3, x_4)+l'(x_0, x_1)=0 \}
\] 
and meets $\Pi$ in a line $L_{(\lambda:\mu)}=\{x_0=x_1= l(x_2, x_3,
x_4)=0\}$, while in the second case $\Pi'_{(\lambda:\mu)}= \Pi$.
\begin{step}
The cubic fibration $X$ does not contain two planes of the form:  
\[
\{ \alpha t_0 +\beta t_1= x=0 \}.
\]
\end{step}
I show that if $X$ contains two distinct planes of this form, every
fibre of $X$ is reducible. This contradicts $Y$ having terminal singularities.

Indeed, we may change coordinates on $\PS^1$ so that the two
reducible fibres lie over $0$ and $\infty$. Let us assume that $X$ contains the planes $\{t_0=x=0 \}$ and
$\{t_1=x=0\}$. Then the equation of $X$ reads:
\[
t_0t_1 A_{3M-L} + x B_{2M+2L}=0, 
\] 
where $A_{3M-L}$ is a section of $\vert 3M-L \vert $  and $B_{2M+2L}$
is a section of $\vert 2M+2L \vert$. Every monomial in
$H^0(\F_{0,0,1}, 3M-L)$ is divisible by $x$, hence $X$ is reducible.

From the point of view of $Y\subset \PS^4$, if $X$ contains
two such planes, the hyperplane sections
of $Y\cap H_{(0:1)}$ and $Y \cap H_{(1:0)}$ both contain the plane
$\Pi=\{x_0=x_1=0\}$ with multiplicity $2$, so that in the equation
\eqref{eq:37}, both $a_3$ and $b_3$ lie in the ideal of $\Pi$. The
quartic $Y$ has generic multiplicity at least $2$
along $\Pi$.

\begin{step}
Assume that $X$ has at least two reducible fibres.
Let $\Pi_1= \Pi_{(\lambda_1:\mu_1)}$ and
$\Pi_2=\Pi_{(\lambda_2:\mu_2)}$ be two planes lying in distinct fibres of $X$
and denote by $\Pi'_1= \Pi'_{(\lambda_1: \mu_1)}$ and
$\Pi'_2=\Pi'_{(\lambda_2: \mu_2)}$ the planes of $\PS^4$ to which
they are isomorphic. If $\Pi'_1$ and $\Pi'_2$ are both distinct
from $\Pi$, they meet at a point. 
\end{step}

I show that if $X$ contains two such planes $\Pi_1$ and $\Pi_2$,
associated to planes $\Pi'_1$ and $\Pi'_2$,  
that intersect $\Pi$ along
the same line $L_1=L_2=L$, $Y$ has multiplicity $2$ along the line
$L$. This contradicts $Y$ having isolated singularities.

The planes $\Pi'_1$ and $\Pi'_2$ of $\PS^4$ lie in distinct
hyperplanes of $\PS^4$ because they correspond to distinct fibres of
the cubic; $\Pi'_1$ and $\Pi'_2$ intersect $\Pi$ in lines $L_1$ and $L_2$
respectively. If $L_1=L_2$, without loss of generality, I may assume
that $\Pi'_1=\{x_0= x_2=0\}$ and $\Pi'_2=\{x_1=x_2=0\}$ after
coordinate change on $\PS^4$. The planes $\Pi'_1$ and $\Pi'_2$ both
intersect the plane $\Pi$ along the line $L_1=L_2=L=\{x_0=x_1=x_2=0\}$. 

The plane $\Pi'_1$ (resp. $\Pi'_2$) lies on $Y$ if and only if, in \eqref{eq:37}, the
homogeneous form $b_3$ is in the ideal $\langle x_0,x_2 \rangle$ of $\Pi'_1$
(resp. the homogeneous form $a_3$ is in the ideal $\langle x_1, x_2 \rangle$ of
$\Pi'_2$). The quartic $Y$ then has multiplicity $2$ along the line $L$.   

\begin{step}
If $X$ contains three planes $\Pi_1, \Pi_2$ and $\Pi_3$ which lie in
distinct fibres and correspond to planes $\Pi'_1, \Pi'_2$ and $\Pi'_3$
lying on $Y$ that are all distinct from $\Pi$, the three planes do not
meet at a single point. 
\end{step}

I show that if the planes $\Pi'_1, \Pi'_2$ and
$\Pi'_3$ meet at a point $P$, the quartic $Y$ does not have isolated
cDV singularities. More precisely, I show that $X$, a small partial
$\Q$-factorialisation of $Y$, does not have isolated singularities. The cubic
fibration $X$ is hence not terminal: $P \in
Y$ is not a cDV point.

Up to coordinate change on $\PS^4$, we may assume that
\begin{eqnarray*}
\Pi'_1=\{x_0=x_2=0\}\\ 
\Pi'_2= \{x_0=x_2=0\}\\
\Pi'_3= \{x_0+x_1= x_2+x_3+l(x_0,x_1)=0\}, 
\end{eqnarray*}
where $l$ is a linear form.

Equivalently, up to global coordinate change on the scroll, $X$
contains the following $3$ planes:
\begin{eqnarray*}
\Pi_1=\{t_0=x_2=0\}\\
\Pi_2=\{t_1=x_3=0\}\\
\Pi_3=\{t_0 + t_1=x_2+x_3+l(t_0,t_1)x=0\}.  
\end{eqnarray*}

The equation of $X$ therefore has to be in the ideal spanned by the monomials:
\begin{eqnarray*}
I =\{ (t_0+t_1)t_0t_1, t_0t_1(x_2+x_3+l(t_0,t_1)x), (t_0+t_1)t_1x_2,
(t_0+t_1)t_0x_3, \\t_1x_2(x_2+x_3+l(t_0,t_1)x),
 t_0x_3(x_2+x_3+l(t_0,t_1)x), (t_0+t_1)x_2x_3,\\ x_2x_3(x_2+x_3+l(t_0,t_1)x)\}, 
\end{eqnarray*}

The cubic fibration $X$  has multiplicity $2$ along
$\Gamma=\{x=x_2=x_3=0\}$. 

Indeed, the only monomials in $I$ which do
not have multiplicity $2$ along $\Gamma$ are $(t_0+t_1)t_0t_1$,
$t_0t_1(x_2+x_3+l(t_0,t_1)x)$, $(t_0+t_1)t_1x_2$ and  $(t_0+t_1)t_0x_3$. However,
$(t_0+t_1)t_0t_1$ is a section of $\vert 3L \vert$, so that if it
appears in the equation of $X \in \vert 3M+L \vert$, it is multiplied
by $A_{3M-2L}$, a section of
$\vert 3M-2L \vert$. Any monomial in $\vert 3M-2L \vert$ is divisible
by $x^2$, hence the term $(t_0+t_1)t_0t_1 \times A_{3M-2L}$ which
appears in the equation of $X$ (if non-zero) has multiplicity $2$
along $\Gamma$. 

Similarly, $t_0t_1(x_2+x_3+l(t_0,t_1)x)$, $(t_0+t_1)t_1x_2$ and
$(t_0+t_1)t_0x_3$ are sections of $\vert M+2L \vert$ and have
multiplicity $1$ along $\Gamma$. If they appear in the equation of
$X$, they are multiplied by a section of $\vert 2M-L \vert$. Any monomial
of $\vert 2M-L \vert$ is divisible by $x$, so that these terms, if
non-zero, also have multiplicity $2$ along $\Gamma$.

As $Y$ has terminal Gorenstein singularities, $X$, a small partial
resolution of $Y$, also has isolated
singularities: this yields a contradiction. 

In terms of $Y$, if $X$ contains such a combination of planes, the singular point
$(0:0:0:0:1)\in Y$ is not a cDV point. 


If $X$ contains three planes $\Pi_1$,
$\Pi_2$ and $\Pi_3$ lying in distinct fibres, which correspond to planes $\Pi'_1$, $\Pi'_2$ and
$\Pi'_3 \subset Y$ all distinct from $\Pi$, without loss of generality, up to
coordinate change on $\PS^4$, we may
assume that:
\begin{eqnarray*}
\Pi'_1=\{x_0= x_2=0 \} \subset H_{(0:1)}\\
\Pi'_2=\{x_1= x_3=0 \}\subset H_{(1:0)}\\
\Pi'_3=\{x_0+x_1= x_4=0 \} \subset H_{(1:1)}.  
\end{eqnarray*}
Equivalently, up to global change of coordinates on the scroll, we may assume that $X$ contains
the following $3$ planes:
\begin{eqnarray*}
\Pi_1=\{t_0= x_2=0 \} \subset X_{(0:1)}\\
\Pi_2=\{t_1= x_3=0 \}\subset X_{(1:0)} \\
\Pi_3=\{t_0+t_1= x_4=0 \} \subset X_{(1:1)}
\end{eqnarray*}
\begin{step}
The cubic fibration $X$ has at most $4$ reducible fibres.  
\end{step}

I argue by contradiction. Assume that $X$ has at least $5$ reducible
fibres and denote by $\Pi_1, \cdots , \Pi_5$ planes that lie in
distinct reducible fibres. Denote by $\Pi'_1, \cdots, \Pi'_5 \subset \PS^4$ the
planes lying on $Y$ naturally associated to $\Pi_1, \cdots ,
\Pi_5$. Recall that Step 1 shows that at most $1$ of the planes
$\Pi'_1, \cdots , \Pi'_5$ is isomorphic to $\Pi$. Without loss of
generality, I may assume that the four planes $\Pi'_1,\Pi'_2, \Pi'_3$
and $\Pi'_4$ are distinct from $\Pi$ and denote by $L_i= \Pi\cap
\Pi'_i$ for $1\leq i\leq 4$. Steps $2$ and $3$ show that the lines
$L_i \subset \Pi$ are distinct and that no three of these lines meet
at a point. Up to change of coordinates, we may therefore assume that:
\begin{eqnarray*}
\Pi_1=\{t_0= x_2=0 \} \subset X_{(0:1)}\\
\Pi_2=\{t_1= x_3=0 \}\subset X_{(1:0)} \\
\Pi_3=\{t_0+t_1= x_4=0 \} \subset X_{(1:1)}
\end{eqnarray*}
and that
\begin{equation*}
\Pi_4=\{t_0+\lambda t_1= ax_2+bx_3+cx_4+ l(t_0, t_1)x=0 \} \subset X_{(1:\lambda)},
\end{equation*}
where $a, b$ and $c$ are constants and $l$ is a linear form. The
constant $\lambda \neq 0,1$ or $\infty$ because no two of the planes
$\Pi_i$ lie in the same fibre of $X$. The
constants $a, b$ and $c$ are all non-zero,  as otherwise either two of
the lines $L_i$ coincide, or three of the lines $L_i$ meet at a
point. Up to rescaling, we may assume that:
\begin{equation*}
\Pi_4=\{t_0+\lambda t_1= x_2+x_3+x_4+ l(t_0, t_1)x=0 \} \subset X_{(1:\lambda)}.
\end{equation*}
\begin{cla}
The plane $\Pi'_5 \subset Y$ is not $\Pi$.   
\end{cla}
If $X$ has a fifth reducible fibre containing a plane $\Pi_5$ such
that $\Pi'_5=\Pi$, $X$ contains a plane of the form:
\[
\Pi_5=\{t_0+\mu t_1= x=0 \}.
\] 
with $\mu \neq 1,0,\infty, \lambda$.

The equation of $X$ may be written uniquely in the form:
\begin{equation}
\label{eq:4}
 A(t_0, t_1)x^3 + B(t_0, t_1, x_2, x_3, x_4)x^2+ C(t_0, t_1, x_2, x_3,
x_4)x+D(t_0, t_1, x_2, x_3, x_4)=0,  
\end{equation}
with $D$ in the ideal   
\begin{multline}
  \label{eq:20}
\{ t_0x_3x_4(x_2+x_3+x_4),
t_1x_2x_4(x_2+x_3+x_4),\\ (t_0+t_1)x_2x_3(x_2+x_3+x_4), (t_0+\lambda t_1)
x_2x_3x_4\},
\end{multline}
i.e. $D$ is a linear combination of the above monomials.

In the expression \eqref{eq:4}, $D$ has to be of the form \[D(t_0, t_1,
x_2, x_3, x_4)= (t_0+\mu t_1)P(x_2, x_3, x_4),\] with $P$ a
polynomial of degree $3$. 
Writing $D$ as a linear combination of the $4$ monomials in \eqref{eq:20} shows
that the only possibility is for $D$ to be equal to a scalar
multiple of one of these monomials. In this case, $\mu$ is equal to
one of $0,1, \infty$ or $\lambda$, so that $\Pi_5$ lies in one of the
fibres above $(0:1),(1:0),(1:1)$ or $(1, \lambda)$. This is a contradiction.

\begin{cla}
 $X$ cannot have a fifth reducible fibre that contains a
plane $\Pi_5$ such that $\Pi'_5$ is distinct from $\Pi$.
\end{cla}

If $X$ contains such a plane $\Pi_5$, the equation of $\Pi_5$ is of
the form:
\[
\Pi_5=\{t_0+\mu t_1= \alpha x_2+ \beta x_3 + \gamma x_4+
l'(t_0,t_1)x=0 \},
\] 
where $l'$ is a linear form and $\alpha, \beta$ and $\gamma$ are
constants. As any three lines of $L_1,\cdots , L_5$ have to satisfy
the conditions of Steps $2$ and $3$, $\alpha, \beta$ and $\gamma$ are
all non-zero. I may assume that $\alpha=1$. Considering triples of
lines $(L_4, L_5, L_i)$ for $1\leq i\leq 3 $ shows that $\beta\neq 1$,
$\gamma\neq 1$, and $(\beta: \gamma)\neq (1:1)$.

In the expression \eqref{eq:4}, $D(t_0, t_1)$
may be written as a linear combination of the monomials
\begin{multline*}
\{ t_0x_3x_4(x_2+x_3+x_4),
t_1x_2x_4(x_2+x_3+x_4), \\(t_0+t_1)x_2x_3(x_2+x_3+x_4),(t_0+\lambda t_1)
x_2x_3x_4\}\end{multline*}
 or as a linear combination of the monomials 
\begin{multline*}
\{ t_0x_3x_4(x_2+\beta x_3+\gamma x_4),
t_1x_2x_4(x_2+\beta x_3+\gamma x_4),\\(t_0+t_1)x_2x_3(x_2+\beta
x_3+\gamma x_4),(t_0+\mu t_1)
x_2x_3x_4 \}.  
\end{multline*}

Equating these two expressions, one finds that either $D=0$ or
$\lambda=\mu$. If $D=0$, every fibre of $X$ is reducible and this yields a
contradiction. If $\lambda= \mu$, $\Pi_5$ is contained in the
reducible fibre over $(1:\lambda)$. 

This show that $X$ has at most $4$ reducible fibres. 
\end{proof}
I have proved the following:
\begin{thm}[Main Theorem $1$]
\label{thm:10}
Let $X_4^3 \subset \PS^4$ be a quartic $3$-fold. Then the defect of
$X$ is at most:
\begin{enumerate}
\item $8$ if $X$ does not contain a plane or a quadric,
\item $11$ if $X$ contains a quadric but no plane,
\item $15$ if $X$ contains a plane.
\end{enumerate}
\end{thm}
\begin{rem}
  \label{rem:13}
Let $Y$ be a terminal Gorenstein Fano $3$-fold of Fano index $1$. Denote by
$g$ its genus. It is known that if $g \geq 6$, $Y$ does not contain a
plane. If $g$ is equal to $4$ or $5$ and if $Y$ contains a plane
$\Pi$, a similar analysis can be carried out. Let $X$ be the crepant blow up
of $\Pi$. Projecting $Y$ away from $\Pi$ shows that $X$ is a conic bundle over $\PS^2$ if $g=4$, and is birational to
$\PS^3$ if $g=5$. I make the following conjecture: 
\end{rem}

\begin{con}
\label{con:1}
  \begin{enumerate}
  \item
 The defect of a non $\Q$-factorial terminal Gorenstein
 $Y_{2,3}\subset \PS^5$ with Picard rank $1$ is at most $8$,  
  \item
 The defect of a non $\Q$-factorial terminal Gorenstein
 $Y_{2,2,2}\subset \PS^6$ with Picard rank $1$ is at most $8$,
\item
 The defect of a non $\Q$-factorial terminal Gorenstein
 Fano $3$-fold with Picard rank $1$ of genus $g \geq 6$ is at most $\big[\frac{12-g}{2}\big]+5$.
\end{enumerate}
\end{con}

\begin{rem}
  \label{rem:47}
I would like to conjecture that the defect of a terminal Gorenstein
double sextic with Picard rank $1$ ($g=2$) is at most $18$. I
unfortunately have no evidence to support that conjecture.
\end{rem}

\section{Deformation Theory}
\label{sec:deformation-theory}

Let $X$ be a weak$\ast$ Fano $3$-fold and $Y$ its anticanonical model.
In this section, I show that deformation theory can be used to
determine invariants of the intermediate weak$\ast$ Fano $3$-folds
$X_i$ and to describe the extremal contractions encountered when running the
Minimal Model Program on $X$. It
follows from Namikawa's results (Theorem~\ref{thm:16}) and from
Lemma~\ref{lem:27} that the degrees of the intermediate weak$\ast$
Fano $3$-folds
$X_i$ can only take values among the list of degrees of non-singular
Fano $3$-folds with Picard rank less than or equal to $\rho(Y)$, the
Picard rank of $Y$. 
This result has been used in
sections~\ref{sec:bound-defect-some-1} and \ref{sec:fano-3-folds} to
bound the Picard rank of $X$, and hence the defect of $Y$.
I show that each step $\phi_i \colon X_i \to X_{i+1}$ of the MMP on $X$ can
be deformed--- in a suitable sense--- to the contraction of an extremal
ray on a non-singular weak
Fano $3$-fold with Picard rank $2$.  
  
First, I recall the set-up of the theory of deformation functors
of Schlessinger and Lichtenbaum \cite{MR0209339}. Second, I survey
results obtained by Friedman and Namikawa on the deformation theory
of \emph{generalised Fano $3$-folds} \cite{MR1489117,MR848512}. 
Generalised Fano
$3$-folds are algebraic $3$-folds
that admit a small birational proper map to a Fano $3$-fold $Y$ with
terminal Gorenstein singularities. Weak$\ast$ Fano $3$-folds are
generalised Fano $3$-folds, but the category of
generalised Fano $3$-folds is larger, as
 it does not require
$\Q$-factoriality. For instance, any small partial
$\Q$-factorialisation of a terminal Gorenstein Fano $3$-fold $Y$ is a
generalised Fano $3$-fold. Third, I recall and extend some results
of Koll\'ar and Mori \cite{MR1149195} on deformations of extremal
contractions.

\subsection{Theoretical setup}
\label{sec:theoretical-setup}
Let $k$ be a field and let $\Lambda$ be a complete local noetherian
$k$-algebra with residue field $k$ and maximal ideal
$\mathfrak{m}_{\Lambda}$.
Consider $\mathcal{C}_{\Lambda}$, the category of local
$\Lambda$-algebras with residue field $k$ and local homomorphisms.
\begin{dfn}
  \label{dfn:7}
Consider a covariant functor $F \colon \mathcal{C}_{\Lambda}\to Sets$.
  \begin{enumerate}
\item Let $k[\epsilon]$ with $\epsilon^2=0$ be the ring of dual
  numbers over $k$; the \emph{tangent space} to $F$ is
  $t_F= F(k[\epsilon])$.
\item
The functor $F$ \emph{has a hull} if there exists a complete
local $\Lambda$-algebra $R$ and a morphism  
\[
\phi \colon \Hom_{\text{local $\Lambda$-algebras}}(R,\cdot) \to F, 
\]
 such that for all $A \in \mathcal{C}_{\Lambda}$, $\phi(A)$ is
 surjective, and such that $\phi$ induces an isomorphism on tangent spaces.
\item  The functor $F$ is
  \emph{pro-representable} if it is isomorphic to a functor of the form:
\[
F \simeq \Hom_{\text{local $\Lambda$-algebras}}(R,\cdot), 
\]
with $R$ a complete local $\Lambda$-algebra such that $R/
\mathfrak{m}_R^n \in \mathcal{C}_{\Lambda}$ for all $n$.
\item
Let $A,A' \in \mathcal{C}_{\Lambda}$ and $p\colon A'\to A$ be a
morphism. The morphism $p$ is a \emph{small extension} if $\ker p$ is
a principal ideal annihilated by the maximal ideal of $A'$.  
\end{enumerate}
\end{dfn}
Schlessinger states the following conditions for a functor of Artin
rings to have a hull or to be pro-representable.

\begin{thm}[\cite{MR0217093}]
\label{thm:21}
Let $F\colon \mathcal{C}_{\Lambda} \to Sets$ be a covariant functor
such that $F(k)$ is restricted to one element $(e)$. Let $A, A'$ and
$A''$ be objects of $\mathcal{C}_{\Lambda}$, and let $A'\to A$ and $A''\to A$
be morphisms in $\mathcal{C}_{\Lambda}$. Consider the natural map:
\[
\psi_{A,A',A''} \colon F(A'\times_{A} A'')\to F(A')\times_{F(A)} F(A'')
\]
\begin{enumerate}
\item $F$ has a hull if and only if the following conditions $H1$ to
  $H3$ hold:
  \subitem[H$1$] $\psi_{A,A',A''}$ is a surjection whenever $A''\to A$
  is a small extension,
  \subitem[H$2$] $\psi_{k,A',k[\epsilon]}$ is a bijection,
  so that $F(A'\times_k k[\epsilon])\simeq F(A')\times_{(e)}t_F$,
  \subitem[H$3$] $\dim_k t_F$ is finite.
\item $F$ is pro-representable if moreover:
  \subitem[H$4$] For any small extension $A'\to A$, $\psi_{A, A', A'}$
  is a bijection.
\end{enumerate}
\end{thm}
\begin{rem}
\label{rem:1}
\begin{enumerate}
\item The condition H$2$ applied to $A'=k[W]$ for $W$ a $k$-vector
  space endows $t_F$ with a canonical $k$-vector space structure.
\item Let $A$ be an element of $\mathcal{C}_{\Lambda}$ and let $I$ be
  an ideal annihilated by
  the maximal ideal of $A$; $A\times_{A/I} A$ is naturally isomorphic
  to $A\times_{k}k[I]$ by the map sending $(x,y)$ to $(x,
  x_0-x+y)$, where $x_0$ denotes the $k$-residue of $x$. 

Let
  $p\colon A'\to A$ be a small extension with kernel $I$. If
  condition H$2$ holds, $\psi_{A,A',A'}$ may be written:
  \[\psi_{A, A', A'}\colon F(A')\times_{F(A)}(t_F \otimes I) \stackrel{\sim}\to
  F(A')\times_{F(A)} F(A').\] 

For any $\eta \in F(A)$, $\psi$
  defines a group action of $t_F \otimes I$ on $F(p^{-1}(\eta))\subset
  F(A')$. Condition H$1$ ensures that this action is transitive,
  while condition H$4$ makes $F(p^{-1}(\eta))$ a formally principal homogeneous
  space (i.e. a torsor) under $t_F \otimes I$. 

If the first three conditions are
  satisfied, pro-representability is obstructed by fixed points of the
  action of $t_F\otimes I$. In the geometric
  setting, a point $\eta'\in F(p^{-1}(\eta))$ is fixed under the
  action of $t_F\otimes I$, if there is an automorphism of an
  object $y\in [\eta]$ that cannot be extended to an automorphism of
  some object of $[\eta']$. 
\end{enumerate}
\end{rem}
Let $\hat{\mathcal{C}}$ be the category of complete local
$\Lambda$-algebras $A$ such that $A/
\mathfrak{m}_A^n \in \mathcal{C}_{\Lambda}$ for all $n$. The functor
$F$ can be extended to $\hat{\mathcal{C}}$ by defining $\hat{F}(A)$ as
$\projL F(A/\mathfrak{m}_A^n)$. For any such ring $R\in \hat{\mathcal{C}}$,
$h_R=\Hom(R,\cdot)$ defines a functor on $\mathcal{C}$. For any
functor on $\mathcal{C}$, there is a canonical isomorphism $\hat{F}(R)
\simeq \Hom (h_R, F)$.

\begin{dfn}
\label{dfn:9}
\begin{enumerate}
\item 
$(A,\xi)$ is a \emph{pro-couple} for $F$ if $A\in \mathcal{C}$ and $\xi
\in F(A)$. A morphism of pro-couples $u \colon (A,\xi)\to (A',\xi')$
is a morphism in $\mathcal{C}$ such that $F(u)(\xi)=\xi'$. 
\item
A pro-couple $(R, \xi)$ for $\hat{F}$ is \emph{pro-representing} if $h_R \to
F$ is an isomorphism induced by $\xi$ in the following sense. Since by
definition $\xi= \projL \xi_n \in \hat{F}(R)$, it is possible to
associate to $\xi$ a morphism
of functors $\Phi_{\xi}$ defined as follows: for any  morphism $u\colon R\to
A$ factoring through $u_n\colon R/\mathfrak{m}_R^n \to A$ for
some $n$, let
\[
\Phi_{\xi}(A):u \in h_R(A)\to F(u_n)(\xi_n) \in F(A).
\]
\end{enumerate}
\end{dfn}

Let $X$ be a pre-scheme over $k$. An \emph{infinitesimal deformation} of
$X/k$ to a local $\Lambda$-algebra $A$ is a diagram:
\[
\xymatrix{
X \ar[r]^{i} \ar[d] & Y \ar[d]\\
\Spec k \ar[r] & \Spec A,
}
\]
with $X \simeq Y \times_{\Spec A}\Spec k$. The pre-scheme $Y$ is
required to be flat over $\Spec A$ and
$i$ is necessarily a closed immersion. 

If $Y'/A$ is another infinitesimal
deformation of $X/k$, it is equivalent to $Y/A$ if there exists a
morphism of pre-schemes $f\colon Y\to Y'$ defined over $A$ that
induces the identity on the closed (or central) fibre $X$. 

Given a deformation $Y/A$ and a morphism $A \to B$, there is an induced
deformation $(Y\otimes_A B)/B$. The notion of
morphism of deformations is thus well defined. 

Define a deformation functor $\mathcal{D}$ of $X/k$ by associating to each
local $\Lambda$-algebra $A$ the equivalence class of infinitesimal
deformations of $X/k$ to $A$. Notice that
$\mathcal{D}(k)=\{X/k\}$. 
Let the obstruction space $T$ of
$\mathcal{D}$ be the space making, for any $A, A'\in \mathcal{C}$ and $p\colon A'\to A$ small, the
  following sequence exact:
\[
\mathcal{D}(A')\to \mathcal{D}(A) \to  T \otimes \ker p \to 0.  
\]
The definition of the obstruction space is functorial and is
consistent with Remark~\ref{rem:1}.

\begin{thm}[\cite{MR0217093}]
If $X$ is proper over $k$, then $\mathcal{D}$ has a hull $(R, \xi)$. 
The pro-couple $(R, \xi)$ pro-represents $\mathcal{D}$ if and only if, for each small
extension $A'\to A$ and for each deformation $Y'/A'$ of $X/k$, every
automorphism of the deformation $Y'\otimes_{A'}A$ is induced by an
automorphism of $A$.
\end{thm}

\begin{dfn}
\label{dfn:10}
Let $X$ be a projective variety. The \emph{Kuranishi space} $\Def(X)$ of $X$ is
the hull of the functor $\mathcal{D}$, that is, the semi-universal (or
miniversal) space
of flat deformations of $X/k$. If $\mathcal{D}$ is
pro-representable, the pro-representing couple $(\Def(X),
\mathcal{X})$ is a universal deformation object and $\mathcal{X}$ is
called the \emph{Kuranishi family} of $X$.  
\end{dfn}

Lichtenbaum, Schlessinger \cite{MR0209339} and Illusie
\cite{MR0491680} make explicit the relation between the problem of infinitesimal deformations of
$X/k$ and that of extensions of local
$\Lambda$-algebras. In that spirit, they relate infinitesimal
deformations of $X/k$ to the cotangent complex of $X$.
\begin{pro}[\cite{MR0209339, MR0491680}]
  \label{pro:5}
Let $X$ be a proper locally complete intersection (lci) variety.
Let $A\in \mathcal{C}_{\Lambda}$. The obstruction to the existence of
a flat deformation of $X/k$ to $A$ is a class $\omega$ lying in
$\Ext^2(\Omega^1_X,\mathcal{O}_X)$. If $\omega$ is zero, the
isomorphism classes of deformations of $X/k$ to $A$ is a torsor under
$\Ext^1(\Omega^1_X,\mathcal{O}_X)$ and the group of automorphisms of a given
deformation is canonically identified with
$\Ext^0(\Omega^1_X,\mathcal{O}_X)= \Hom (\Omega^1_X, \mathcal{O}_X)$. 
\end{pro}

\begin{rem}
\label{rem:23}
In particular, if $\Ext^2(\Omega^1_X,\mathcal{O}_X)=(0)$, the Kuranishi
space $\Def(X)$ of $X$ is a smooth analytic complex space. The tangent space
of the deformation functor (and of $\Def(X)$ at the origin) is $\Ext^1(\Omega^1_X, \mathcal{O}_X)$. If
$\Hom(\Omega^1_X, \mathcal{O}_X)=(0)$, the deformation functor of $X/k$
is pro-representable (or, equivalently, the Kuranishi space is universal). 
The functor of first order \emph{local
  deformations} is defined similarly, and is controlled by the sheaves $\EExt^i(\Omega^1_X,
  \mathcal{O}_X)$.
\end{rem}

Let $(X,D)$ be a pair of a normal $3$-fold over $k$ and an effective Cartier
divisor. An \emph{infinitesimal deformation of the
pair} $(X,D)/k$ to a local $\Lambda$-algebra $A$ is a proper flat
morphism $X_A \to \Spec A$ together with an effective divisor $D_A$ on
$X_A$, such that
$X_A\times_{\Spec A}\Spec k= X$ and $D_{A}\times_{\Spec A}\Spec k=
D$.  

Define, as above, a deformation functor $\mathcal{D}'$ of $(X, D)$ by
associating to each local $\Lambda$-algebra $A$ the equivalence class
of infinitesimal deformations of the pair $(X,D)$ to $A$; in
particular $\mathcal{D}'(k)=\{(X, D)/k\}$. 

I recall Kawamata's construction of the sheaf $\Omega^1_X(\log D)$ of logarithmic differential
forms \cite{MR814013} when $D$ is a not necessarily non-singular
effective Cartier divisor. 

Let $X$ be a complete algebraic variety,  $(D_j)_{j\in J}$ a finite
number of effective divisors, and let $D$ be the Cartier divisor
$D=\sum _{j\in J} D_j$. The Cartier divisor $D$ does not necessarily have normal
crossings.

Let $U= \Spec A \subset X$ be an affine open subset with a closed
embedding $U \to R$ into a non-singular affine algebraic variety $R
=\Spec B$ such that:
\begin{enumerate}
\item For each $j \in J$, $D_j \cap U= (h_j)$ with
$h_j \in A$ is extended to a non-singular prime divisor $D'_j=(h'_j)$
with $h'_j \in B$ on $R$;
\item the divisors $D'_j$ intersect
transversally.
\end{enumerate}
 Denote by $D'=(h')$ the normal crossings divisor $\sum_{j
  \in J}D'_j$ on R.

The usual residue exact sequence for a normal crossing divisor $D'$ on
a non-singular affine algebraic variety $R= \Spec B$ is:
\begin{equation*}
0 \to \Omega^1_B \to \Omega^1_B(\log (h')) \to \bigoplus_{j\in J}
B/(h'_j)\to 0.
\end{equation*}
 The sequence 
\begin{equation}
   \label{eq:42}
0 \to \Omega^1_B \otimes_B A\stackrel{\alpha} \to \Omega^1_B(\log (h'))\otimes_B A \to \bigoplus_{j\in J}
B/(h'_j)\otimes_B A=\bigoplus_{j\in J}A/(h_j) \to 0
 \end{equation}
is also exact, because $\Tor_1^B(B/(h'_j), A)$ is trivial.
Let $I$ be the ideal of $B$ such that $A=B/I$, so the conormal sequence is:
\begin{equation}
  \label{eq:43}
I/I^2 \to \Omega^1_B \otimes_B A \stackrel{\beta}\to \Omega^1_A \to 0.
\end{equation}
The exact sequence \eqref{eq:42} pushed forward along $\beta$ yields
an extension of modules:
\begin{equation}
  \label{eq:45}
 0 \to \Omega^1_A\to \mathcal{M}\to \bigoplus_{j\in J} A/(h_j)\to 0, 
\end{equation}
where the module $\mathcal{M}$ is
\[
\mathcal{M}= \frac{(\Omega^1_B(\log (h'))\otimes_B A)\oplus
  \Omega^1_A}{(\alpha\oplus -\beta)(\Omega^1_B \otimes_B A)}.
\]
If $e_j, j\in J,$ are the images in $\mathcal{M}$ of the elements
$\frac{dh'_j}{h'_j}\in \Omega^1_B(\log (h'))$, $\mathcal{M}$ is
generated by $\Omega^1_A$ and the elements $e_j$. 
\begin{dfn}
\label{dfn:15}
The sheaf $\Omega^1(\log D)$ is defined locally by setting
\[\Omega^1_A(\log D)= \mathcal{M}.\]  
\end{dfn}
\begin{rem}
  \label{rem:44}
If $X$ is non-singular and $D$ has normal crossings, the definition of
$\Omega^1_X(\log D)$ agrees with Deligne's definition. Note however that
the sheaf $\Omega^1_X(\log D)$ may have torsion, even when $X$ is non-singular.
\end{rem}
\begin{rem}
 \label{rem:16}
The residue sequence \eqref{eq:45} is exact. 
\end{rem}
\begin{lem}\cite{MR814013,MR0146040}
  \label{lem:34}
Let $(X,D)$ be a pair comprising a normal complete $3$-fold and an effective
Cartier divisor. Assume that $X$ and $D$ have no worse than isolated
l.c.i. singularities. The infinitesimal deformations of the pair
$(X,D)/k$ are controlled by $\Omega^1_X(\log D)$. 
Let $A\in \mathcal{C}_{\Lambda}$. The obstruction to the existence of
a flat deformation of $(X,D)/k$ to $A$ is a class $\omega$ lying in
$\Ext^2(\Omega^1_X( \log D),\mathcal{O}_X)$. If $\omega$ is zero, the
isomorphism classes of deformations of $(X,D)/k$ to $A$ is a torsor under
$\Ext^1(\Omega^1_X(\log D),\mathcal{O}_X)$ and the group of automorphisms of a given
deformation is canonically identified with
$\Ext^0(\Omega^1_X(\log D),\mathcal{O}_X)= \Hom (\Omega^1_X (\log D), \mathcal{O}_X)$.
\end{lem}

\begin{rem}
\label{rem:45}
Let $\mu\colon \widetilde{X}\to X$ be a good resolution of the pair
$(X,D)$. By definition of $\mu$, the divisor $\Exc(\mu)\cup \widetilde{D}$,
where $\widetilde{D}$ is the proper transform of $D$, has simple
normal crossings.
Denote by $E$ the union of $\mu^{\ast}(D)$ and of some irreducible
components of the exceptional locus of $\mu$. Kawamata proves
\cite{MR814013} that there is a natural transformation of functors
$\mathcal{D}''\to \mathcal{D}'$, where $\mathcal{D}'$ is the functor
of deformations of $(X, D)$ and $\mathcal{D}''$ the functor of
deformations of $(\widetilde{X}, E)$. This is consistent with the
construction of the sheaf $\Omega^1_X(\log D)$, which can be
interpreted as relating $\Omega^1_X(\log D)$ to
$\Omega^1_{\widetilde{X}}(\log \widetilde{D})$ for a good resolution $\mu$
of the pair $(X,D)$. 
\end{rem}

\subsection{Deformations of generalised Fano $3$-folds}
\label{sec:deform-theory-gener}
 
\begin{dfn}
\label{dfn:6}
A $3$-fold $X$ is a \emph{generalised Fano $3$-fold} if there exists a
  small birational proper morphism $X \to Y$ to a Fano $3$-fold with
  terminal Gorenstein singularities.
\end{dfn}
\begin{rem}
  \label{rem:25}
  \begin{enumerate}
  \item
Note that $X$ is a weak Fano $3$-fold; in particular $X$ has
terminal Gorenstein singularities.
\item
A weak$\ast$ Fano $3$-fold is a generalised Fano $3$-fold, but the
notion of generalised Fano $3$-fold does not
require $\Q$-factoriality. Any small partial
$\Q$-factorialisation of a terminal Gorenstein Fano $3$-fold is a
generalised Fano $3$-fold.  
  \end{enumerate}
\end{rem}

Let $X$ be a generalised Fano $3$-fold and let $D$ be a general
section of the anticanonical linear system
$\vert -K_X \vert$. Let $Y$ be the anticanonical model of $X$.

Recall that Theorem~\ref{thm:1} shows that a
general member $D$ of $\vert -K_X \vert$ is a K$3$-surface with no
worse than Du Val singularities. The pair $(X,D)$ is therefore log canonical. 
Let $\Sing(X) \cup \Sing(D)=\Sigma= \{P_1, \cdots , P_n\}$. 

If $Y$ is not birational
to a special complete intersection with a node $X_{2,6}\subset
\PS(1^4, 2, 3)$, there always exists a non-singular section of $\vert -K_X
\vert$. If $Y$ is, in addition, not monogonal, the general section $D
\in \vert -K_X \vert$ is non-singular. In particular, if $Y$ has
Picard rank $1$ and genus $g$ greater than or equal to $3$, the
general section $D$ is non-singular.

\begin{rem}
  \label{rem:43}
I use in this section some results on the Hodge theory of
surfaces with rational double points.
If $Z$ has isolated quotient
singularities and if $j \colon Z \smallsetminus \Sing(Z)\to Z$ is the
natural inclusion, the complex of coherent sheaves
$\widehat{\Omega}^{\bullet}_Z=
j_{\ast}(\Omega^{\bullet}_{Z\smallsetminus \Sing(Z)})$ is a
resolution of the constant sheaf $\C$  \cite{MR0485870}. For each $p$,
$\widehat{\Omega}^p_Z$ coincides with the double dual of
$\Omega^p_Z$. The spectral sequence of
hypercohomology $E_{1}^{p,q} = H^q(Z, \widehat{\Omega}^p_Z)$ abuts to
$H^{p+q}(Z, \C)$, degenerates at $E_1$ and the induced filtration on
$H^{p+q}(Z, \C)$ coincides with the canonical Hodge filtration. If
$\nu \colon \widetilde{Z}\to Z$ is a resolution of singularities, the
complexes $\widehat{\Omega}^{\bullet}_Z$ and
$\nu_{\ast}\Omega^{\bullet}_{\widetilde{Z}}$ are equal. 
\end{rem}

\begin{lem}\cite{MR1489117,MR1144434} 
\label{lem:19}
Let $X$ be a generalised Fano $3$-fold and $D$ a general
anticanonical section. Then:
  \begin{enumerate} 
\item $\Def(X)$ and $\Def(X,D)$ are smooth.
\item $\Def(X,D)$ is universal.
\item The natural map $\phi \colon  \Def(X,D) \to \Def(X)$
    is smooth.
  \end{enumerate}
\end{lem}
\begin{proof}
\setcounter{step}{0}

Denote by $\Omega^1_X(\log D)$ the sheaf of logarithmic differential
forms as constructed in Definition~\ref{dfn:15}; $\Omega^1_{X}(\log
D)$ is not, in general, locally free at points $P_i\in
\Sigma=\Sing(X)\cup \Sing (D)$.  

By Proposition~\ref{pro:5}, the Kuranishi spaces $\Def(X)$ and
$\Def(X,D)$ are smooth analytic spaces if: 
\begin{eqnarray*}
\Ext^2(\Omega^1_{X}(\log D), \mathcal{O}_X)= 
\Ext^2(\Omega^1_{X},\mathcal{O}_X)= (0).  
\end{eqnarray*}

The residue exact sequence \eqref{eq:45}
\begin{eqnarray*}
0 \to \Omega^1_{X} \to \Omega^1_{X}(\log D) \to \mathcal{O}_D \to 0  
\end{eqnarray*}
yields the exact sequence of $\Ext$ groups:
\begin{multline}
\label{eq:15}
0 \to \Hom(\mathcal{O}_D, \mathcal{O}_X) \to \Hom(\Omega^1_{X}(\log
D), \mathcal{O}_X) \to \Hom( \Omega^1_{X}, \mathcal{O}_X) \to\\ 
\to \Ext^1(\mathcal{O}_D, \mathcal{O}_X) \to \Ext^1(\Omega^1_X(\log
D),\mathcal{O}_X)\to \Ext^1(\Omega^1_X,\mathcal{O}_X)\to \\ \to 
 \Ext^2(\mathcal{O}_D, \mathcal{O}_X) \to \Ext^2(\Omega^1_X(\log
D),\mathcal{O}_X)\to \Ext^2(\Omega^1_X,\mathcal{O}_X)\to\cdots 
\end{multline}

\begin{cla}
The $\Ext$ groups $\Ext^i(\mathcal{O}_D, \mathcal{O}_X)$ for $i=0,1,2$
are as follows:
\begin{eqnarray*}
\Ext^0(\mathcal{O}_D, \mathcal{O}_X)=\Ext^2(\mathcal{O}_D, \mathcal{O}_X)=(0)  \\
\Ext^1(\mathcal{O}_D,\mathcal{O}_X)= H^0(X, \mathcal{O}_D(D)).  
\end{eqnarray*}
\end{cla}
The standard resolution of the sheaf $\mathcal{O}_D$
\begin{eqnarray*}
0 \to \mathcal{O}_X(-D)\to \mathcal{O}_X \to \mathcal{O}_D \to 0
\end{eqnarray*}
yields the long exact sequence of $\EExt$ sheaves:
\begin{multline*}
0 \to \HHom(\mathcal{O}_D, \mathcal{O}_X) \to \HHom(\mathcal{O}_X,
\mathcal{O}_X) \to \HHom(\mathcal{O}_X(-D),\mathcal{O}_X) \to \\
\to \EExt^1(\mathcal{O}_D, \mathcal{O}_X) \to \EExt^1(\mathcal{O}_X,
\mathcal{O}_X)\to \EExt^1(\mathcal{O}_X(-D),\mathcal{O}_X) \to \cdots 
\end{multline*}
As $\mathcal{O}_X $ and $\mathcal{O}_X(-D)$ are locally free,
\[\EExt^i(\mathcal{O}_X(-D), \mathcal{O}_X)=\EExt^i(\mathcal{O}_X,
\mathcal{O}_X)=0\] for all $i>0$. The sheaf
$\HHom(\mathcal{O}_D, \mathcal{O}_X)$ is trivial and the long exact sequence
reduces to
\begin{equation*}
0 \to \mathcal{O}_X\to \mathcal{O}_X(D)\to \EExt^1(\mathcal{O}_D,
\mathcal{O}_X)\to 0.  
\end{equation*}
This shows that $\EExt^i(\mathcal{O}_D,\mathcal{O}_X)=0$ for $i\neq 1$
and $\EExt^1(\mathcal{O}_D,\mathcal{O}_X)= \mathcal{O}_D(D)$.
The local-to-global spectral sequence of $\Ext$ sheaves and groups has $E_1$ term
\[E_1^{p,q}= H^q(X, \EExt^p(\mathcal{O}_D, \mathcal{O}_X)),\] and abuts
to $\Ext^{p+q}(\mathcal{O}_D, \mathcal{O}_X)$. This
shows that $\Ext^0(\mathcal{O}_D, \mathcal{O}_X)=(0)$ and that the
following sequence is exact \cite{MR0345092}:
\begin{multline*}
 0 \to H^1(X, \HHom(\mathcal{O}_D, \mathcal{O}_X)) \to
 \Ext^1(\mathcal{O}_D,\mathcal{O}_X)\to  
H^0(X,\EExt^1(\mathcal{O}_D,
 \mathcal{O}_X))\\ \to H^2(X,\HHom(\mathcal{O}_D, \mathcal{O}_X))) \to
 \Ext^2(\mathcal{O}_D, \mathcal{O}_X)\to 0
\end{multline*}
In particular, $\Ext^2(\mathcal{O}_D, \mathcal{O}_X)=(0)$ and
  $\Ext^1(\mathcal{O}_D,\mathcal{O}_X)= H^0(X, \mathcal{O}_D(D))$. 
\begin{step}
    The Kuranishi spaces $\Def(X)$ and $\Def(X,D)$ are smooth.
\end{step}
As $\Ext^2(\mathcal{O}_D,\mathcal{O}_X)=(0)$, if $\Ext^2(\Omega^1_X,
\mathcal{O}_X)=(0)$, both $\Def(X)$ and $\Def(X,D)$ are smooth by \eqref{eq:15}.

The group $\Ext^2(\Omega^1_X, \mathcal{O}_X)$ is Serre dual to
$H^1(X, \Omega^1_X\otimes \omega_X)$. 
The exact sequence  
\begin{eqnarray*}
0 \to \Omega^1_X \otimes \omega_X \to \Omega^1_X \to {\Omega^1_X}_{\vert D}
\to 0  
\end{eqnarray*}
yields the long exact sequence of cohomology groups:
\begin{multline}
\label{eq:36}
0 \to H^0 (X,\Omega^1_X \otimes \omega_X) \to H^0(X, \Omega^1_X) \to H^0(X,
\Omega^1_{X \vert D})\to \\
\to H^1(X, \Omega^1_X \otimes \omega_X)\to H^1(X,
\Omega^1_X) \stackrel{\beta} \to H^1(X,\Omega^1_{X \vert D}) \to \cdots  
\end{multline}
The conormal sequence is exact because $D$ is a Cartier divisor, hence: 
\[
0 \to \mathcal{O}_D(-D) \to \Omega^1_{X\vert D} \to \Omega^1_D \to 0.
\] 
This shows that  $H^0(X,\Omega^1_{X\vert D})=(0)$. Indeed, $D$ is a K$3$
surface with Du Val singularities; by Remark~\ref{rem:43},
$H^0(D,\widehat{\Omega}^1_ D)= H^1(D, \mathcal{O}_D)=(0)$. The divisor
$D$ is nef and big, hence by the Kawamata-Viehweg vanishing theorem $H^0(D,
\mathcal{O}_D(-D))=(0)$. The map $\Omega^1_D \to\widehat{\Omega}^1_D$ is
an injection because $\Omega^1_D$ is torsion free; thus
$H^0(X,\widehat{\Omega}^1_ D)= H^0(X,\Omega^1_ D)=(0)$ and $H^0(X,
\Omega^1_{X \vert D})=(0)$.

Consequently, it is sufficient to show that the map $\beta$ in
\eqref{eq:36} is injective in order to complete the proof.

The surface $D$ is a K$3$ with no worse than Du Val
singularities, hence $H^1(D, \mathcal{O}_D)$ is trivial. The map 
\[
\frac{1}{2i\pi}dlog \colon H^1(D, \mathcal{O}_D^{\ast})\otimes_{\Z} \C
\to H^1(D, \widehat{\Omega}^1_D)
\]
is injective and factors through $H^1(D, \Omega^1_D)$:
\[
\frac{1}{2i\pi}dlog \colon H^1(D, \mathcal{O}_D^{\ast})\otimes_{\Z} \C
\to H^1(D, \Omega^1_D)
\]
is also an injection.

The following diagram is commutative:
\[
\xymatrix{
H^1(X, \mathcal{O}_X^{\ast})\otimes_{\Z} \C \ar[r]^{j} \ar[d]^{\simeq}_{\frac{1}{2i\pi}dlog}
& H^1(D, \mathcal{O}_D^{\ast})\otimes_{\Z} \C  \ar[d]_{\frac{1}{2i\pi}dlog} \\ 
H^1(X, \Omega^1_X)\ar[r]^{k} & H^1(D, \Omega^1_D). 
}
\]
The conormal sequence shows that $k$ factors through $\beta$, hence the
morphism  $H^1(X,
\Omega^1_X)\to H^1(D, {\Omega^1_X}_{\vert D})$ is injective. 

\begin{cla}
The natural restriction map $\Pic X \to \Pic D$ is injective.  
\end{cla}
 
Following \cite{MR0282977}, one
shows that the pair $(X,D)$ satisfies the Lefschetz condition, i.e.
that for any coherent sheaf $\mathcal{F}$, the cohomology groups
$H^i(X-D, \mathcal{F})$ are trivial for $i>1$. 
Denoting by $\widehat{X}$ a formal completion of $X$ along
$D$, $\Pic X \to \Pic \widehat{X}$ is then an injection. The natural map $
\Pic \widehat{X}\to \Pic D$ is also an injection because
$\mathcal{O}_D(-D)$ is nef and big.
In the above diagram, $j$ is injective.

\begin{cla}
If $X$ is a generalised Fano $3$-fold, the map \[\frac{1}{2i \pi}dlog
\colon H^1(X, \mathcal{O}_X^{\ast})\otimes_{\Z} \C \to H^1(X,
\Omega^1_X)\] is an isomorphism.  
\end{cla}

Let $\mu \colon \widetilde{X}\to X$ be a good resolution of $X$. The
Leray spectral sequence shows that the cohomology groups
$H^i(\widetilde{X}, \mathcal{O}_{\widetilde{X}})= H^i(X,
\mathcal{O}_X)$ because the singularities of $X$ are rational. As $X$ is a
generalised Fano $3$-fold, by the Kawamata-Viehweg vanishing
theorem, $H^i(\widetilde{X}, \mathcal{O}_{\widetilde{X}})=(0)$ for $i>0$. The map 
 \[
\frac{1}{2i \pi}dlog
\colon H^1(\widetilde{X}, \mathcal{O}_{\widetilde{X}}^{\ast})\otimes_{\Z} \C \to H^1(\widetilde{X},
\Omega^1_{\widetilde{X}})\simeq H^2(\widetilde{X}, \C)
\]
 is an isomorphism.
Note first that when $X$ has isolated hypersuface singularities and
$P\in \Sing(X)$, the local cohomology groups
$H^i_{\{P\}}(X,\Omega^1_X)=(0)$ for $i=0,1$ \cite{MR1286924}. This implies that the
natural map $\Omega^1_X\to \mu_{\ast}\Omega^1_{\widetilde{X}}$ is an isomorphism.
The diagram
\[
\xymatrix{0 \ar[r] & H^1(X,
  \mu_{\ast}\Omega^1_{\widetilde{X}})=H^1(X,\Omega^1_X) \ar[r] &
  H^1(\widetilde{X},\Omega^1_{\widetilde{X}}) \ar[r] & H^0(X, R^1 \mu_{\ast}\Omega^1_{\widetilde{X}})\\
0 \ar[r] & H^1(X,
  \mu_{\ast}\mathcal{O}_{\widetilde{X}}^{\ast})=H^1(X,\mathcal{O}_X^{\ast}) \ar[r]\ar[u] &
  H^1(\widetilde{X},\mathcal{O}_{\widetilde{X}}^{\ast}) \ar[r] \ar[u]&
H^0(X, R^1 \mu_{\ast}\mathcal{O}_{\widetilde{X}}^{\ast})\ar[u]
}\]
is commutative.
The second vertical map is surjective by the Hodge decomposition of
  the non-singular $3$-fold $\widetilde{X}$. It is sufficient to
  prove that $H^0(X, R^1
  \mu_{\ast}\mathcal{O}_{\widetilde{X}}^{\ast}) \to  H^0(X, R^1
  \mu_{\ast}\Omega^1_{\widetilde{X}})$ is injective. This can be done
  exactly as in \cite{MR1286924}.
Denote by $E$ the exceptional divisor of the resolution $\mu$, by
  $V$ a neighbourhood of $E$ in $\widetilde{X}$, and by $L$ a line bundle
  in $V$. If $\frac{1}{2i \pi}dlog(L)=0$ in $H^0(X, R^1
  \mu_{\ast}\Omega^1_{\widetilde{X}})$, then $L$ is $\mu$-numerically
  trivial. The $3$-fold $X$ has rational singularities: it follows
  that $L$ is a
  torsion line bundle on $V$. 

The map $k$, and hence $\beta$, is an injection; thus $\Def(X)$ and
$\Def(X, D)$ are smooth analytic spaces.

\begin{step}
The group $\Ext^0(\Omega^1_{X}(\log D))$ is trivial; the Kuranishi
space $\Def(X)$ is universal. 
\end{step}

By Serre duality, it is sufficient to prove that:
\[
H^3(X,\Omega^1_{X}(\log D)\otimes \omega_X )=H^3(X,\Omega^1_{X}(\log
D)(-D))=(0). 
\]
The sequence
\[
0\to \Omega^1_X(\log D)(-D)\to \Omega^1_X\to \Omega^1_{D}\to 0
\] 
is exact. As $D$ is a K$3$ surface with rational double points,
$\Omega^1_D$ is torsion free and the natural map $0\to \Omega^1_D\to \widehat{\Omega}^1_D$
has cokernel supported at the singular points. 
In particular, $H^2(D,
\Omega^1_D)= H^2(D, \widehat{\Omega}^1_D)=(0)$ by Hodge symmetry. 
It is sufficient to prove that $H^3(X, \Omega^1_X)=(0)$.

The group $H^3(X,\Omega^1_X)$ is Serre dual to $H^0(X,
\Theta_X(-D))$, where $\Theta_X$ denotes $\HHom(\Omega^1_X,
\mathcal{O}_X)$. Let $\widetilde{X}\to X$ be a good resolution of
$X$. Since $X$ has terminal singularities,
$f_{\ast}(\Theta_{\widetilde{X}}\otimes K_{\widetilde{X}})=
\Theta_X \otimes K_X$. The cohomology group $H^3(\widetilde{X}, \Omega^1_{\widetilde{X}})$ 
is trivial by Hodge symmetry on the non-singular generalised Fano
$\widetilde{X}$; the result follows.
\begin{step}
The natural map $\phi \colon \Def(X,D) \to \Def(X)$ is surjective. 
\end{step}
The long exact sequence \eqref{eq:15} shows that this follows from triviality
of the group $\Ext^2(\mathcal{O}_D, \mathcal{O}_X)$. 
\end{proof}

\setcounter{step}{0}
The canonical local-to-global spectral sequence of $\Ext$ groups relates
global to local first order deformation functors and reads:

\begin{gather*}
0 \to H^1(X ,\HHom(\Omega^1_X, \mathcal{O}_X)) \to \Ext^1(\Omega^1_X,
\mathcal{O}_X) \stackrel{\alpha}\to H^0(X, \EExt^1(\Omega^1_X, \mathcal{O}_X))\to \\ \to H^2(X ,\HHom(\Omega^1_X, \mathcal{O}_X)).
\end{gather*}

The map $\alpha$ can be regarded as the homomorphism between the space of first order
global deformations and the space of first order local deformations.
If $\alpha$ is surjective, $X$ can be deformed to a non-singular
$3$-fold, which can be shown to be a generalised Fano (Remark~\ref{rem:5}). 

Define similarly the homomorphism:
\[
\alpha_{\log} \colon \Ext^1(\Omega^1_X(\log D),
\mathcal{O}_X)\to  H^0(X, \EExt^1(\Omega^1_X(\log D), \mathcal{O}_X)),
\]
which relates the spaces of global and local first order deformations
of $(X,D)$.
\setcounter{step}{0}
\begin{thm}\cite{MR1489117}
\label{thm:16}
  Let $X$ be a generalised Fano $3$-fold with terminal Gorenstein
  singularities. Then $X$ can be deformed to a non-singular generalised
  Fano $3$-fold. In particular, any Fano $3$-fold with Gorenstein
  terminal singularities is smoothable by a flat deformation.
\end{thm}
 More precisely, if $X$ has no worse than ordinary
 double points, $\alpha$ is surjective. In the general case,
 there is a ``good'' direction $\eta \in \Ext^1(\Omega^1_X(\log D),
\mathcal{O}_X)$ such that deforming $X$ along $\eta$ improves the
 singularities. After finitely many infinitesimal deformations along
 ``good'' directions,
 $X$ becomes a generalised Fano with no worse than
 ordinary double points. I give an overview of Namikawa's proof.
\begin{proof}[Sketch Proof]
Let $\Theta_X=\HHom(\Omega^1_X, \mathcal{O}_X)$ be the dual of the
cotangent sheaf and $T^i_X= \EExt^i(\Omega^1_X, \mathcal{O}_X)$.
\begin{step}
Assume that $X$ has no worse than ordinary double points; then
there is a smoothing of $X$. 
\end{step}
\begin{lem}
\label{lem:21}\cite{MR1489117,MR848512}
If $X$ has no worse than ordinary double points, the cohomology group $H^2(X,
\Theta_X)$ is trivial. In particular, $\alpha$ is surjective and $X$ is
smoothable by a flat deformation.  
\end{lem}
\begin{proof}[Proof of Lemma~\ref{lem:21}]
Let $\pi \colon \widetilde{X}\to X $ be a small resolution of $X$.
By definition, the exceptional locus of $\pi$ above an ordinary double
point $P_i$ is a rational curve $C_i$ with normal bundle
$\mathcal{N}_{C_i/ \widetilde{X}}=
\mathcal{O}_{\PS^1}(-1)\oplus\mathcal{O}_{\PS^1}(-1)$.
Friedman shows that if the only singularities are ordinary double
points, $H^0(X, R^1\pi_{\ast} \Theta_{\widetilde{X}})=(0)$ \cite{MR848512}.

The spectral sequence for
$R\pi_{\ast}\Theta_{\widetilde{X}}$ yields the exact sequence \cite{MR0345092}:
\[
H^1(\widetilde{X},\Theta_{\widetilde{X}})\to H^0(X,
R^1\pi_{\ast}\Theta_{\widetilde{X}})\to H^2(X, \Theta_X) \to
H^2(\widetilde{X}, \Theta_{\widetilde{X}}).
\] 
It is therefore sufficient to
prove that $H^2(\widetilde{X}, \Theta_{\widetilde{X}})=(0)$. 
By Serre
Duality on the non-singular $3$-fold $\widetilde{X}$,
$H^2(\widetilde{X}, \Theta_{\widetilde{X}})$ is dual to $H^1(\widetilde{X}, \Omega^1_{\widetilde{X}}\otimes
\omega_{\widetilde{X}})$.
Denote by $\widetilde{D}$ the proper transform of $D$ by $\pi$. The
map $\pi$ is small, hence $\widetilde{D}$ is a section of
$-K_{\widetilde{X}}$ and is a K$3$ surface with no worse than Du Val
singularities. 
Tensoring the residue exact sequence of $(\widetilde{X},
\widetilde{D})$ by $K_{\widetilde{X}}$ shows that the following
sequence is exact:
\begin{equation}
\label{eq:38}
0 \to \Omega^1_{\widetilde{X}}\otimes \omega_{\widetilde{X}}\to
\Omega^1_{\widetilde{X}}(\log \widetilde{D})\otimes \omega_{\widetilde{X}} 
\to  \mathcal{O}_{\widetilde{D}}\otimes \omega_{\widetilde{X}}\to 0.
\end{equation}
The line bundle $-K_{\widetilde{X}\vert \widetilde{D}}$ is nef and
big, so that
$H^{0}(\widetilde{D},\omega_{\widetilde{X}} \otimes
\mathcal{O}_{\widetilde{D}})=(0)$ by the Kawamata-Viehweg vanishing theorem. The
long exact sequence in cohomology associated to \eqref{eq:38} shows
that if
\[H^1(\widetilde{X},\Omega^1_{\widetilde{X}}(\log \widetilde{D})\otimes
\omega_{\widetilde{X}})=
H^1(\widetilde{X},\Omega^1_{\widetilde{X}}(\log
\widetilde{D})(-\widetilde{D}))=(0), \]
then $H^1(\widetilde{X}, \Omega^1_{\widetilde{X}} \otimes \omega_{\widetilde{X}})=(0)$.  
The sequence 
\begin{equation}
\label{eq:39}
0 \to \Omega^1_{\widetilde{X}}(\log \widetilde{D})(-
\widetilde{D})\to\Omega^1_{\widetilde{X}}\to
\Omega^1_{\widetilde{D}}\to 0
\end{equation}
is exact. The surface $\widetilde{D}$ has no worse than Du Val
singularities and $\Omega^1_{\widetilde{D}}$ is torsion free, so that
$\Omega^1_{\widetilde{D}}\to \widehat{\Omega}^1_{\widetilde{D}}$ is an injection. The cohomology group
$H^0(\widetilde{D},\widehat{\Omega}^1_{\widetilde{D}})$ is trivial because
$\widetilde{D}$ is a K$3$ surface and
$H^0(\widetilde{D},\widehat{\Omega}^1_{\widetilde{D}})=(0)$. As $\Pic
\widetilde{X} \to \Pic \widetilde{D}$ is an injection,
$H^1(\widetilde{X}, \Omega^1_{\widetilde{X}})\to
H^1(\widetilde{D},\widehat{ \Omega}^1_{\widetilde{D}})$ is an
injection. In addition, since it factors through
$H^1(\widetilde{D},\Omega^1_{\widetilde{D}})$, the map $H^1(\widetilde{X}, \Omega^1_{\widetilde{X}})\to
H^1(\widetilde{D},\Omega^1_{\widetilde{D}})$ is also an
injection. 

The long exact sequence in cohomology associated to
\eqref{eq:39} shows that $H^1(\widetilde{X},\Omega^1_{\widetilde{X}}(\log
\widetilde{D})(-\widetilde{D}))=(0)$; this completes the proof. 
\end{proof}

\begin{step} 
We define an invariant $\mu(P_i)$ for all $P_i\in \Sing X$, which is strictly
positive unless $P_i$ is an ordinary double point.  
\end{step}

Let $X$ be a generalised Fano $3$-fold and let $P_i \in \Sigma$ be a
singular point. As will be made clear in
Lemma~\ref{lem:24}, if $\mu>0$, there is a small deformation
$\eta$ of
$X$ which ``improves'' the singularity at $P_i$ in the following
sense: for any resolution $\widetilde{X}\to X$ of the singularity at
$P_i$, $\eta$ is not in the image of the map $\Def(\widetilde{X})\to \Def(X)$. 

I first state some results on the local cohomology groups at the
singular set of the sheaves $\Theta_X(\log D)$. Recall that
$\Sigma=\Sing(X)\cup \Sing(D)= \{P_1, \cdots , P_n\}$, that
$U=X\smallsetminus\Sigma$, and denote by $U_i$ a
Stein open neighbourhood of $P_i$ in $X$.

\begin{lem}\cite{MR1489117}
  \label{lem:33}
  \begin{enumerate}\item  
$H^2_{\{P_i\}}(U_i, \Theta_X(-\log D))= H^1(U_i\smallsetminus
  \{P_i\}, \Theta_X(-\log D))$ and both are equal to
$H^0(U_i, \EExt(\Omega^1_X(\log D),
  \mathcal{O}_X))$.
\item $H^2_{\Sigma}(X, \Theta_X(-\log D))= H^1(U, \Theta_X(-\log D))=
  \Ext^1(\Omega^1_X(\log D),\mathcal{O}_X))$.
  \end{enumerate}
\end{lem}
\begin{proof}[Sketch of proof]
Let $D_i= D\cap U_i$ and denote by $U_i^{\ast}$ (resp. $D_i^{\ast}$)
the punctured neighbourhood $U_i\smallsetminus \{P_i\}$
(resp. $D_i\smallsetminus \{P_i\}$). If $D_i$ is
empty, the result is proved by Schlessinger \cite{MR0292830}. 
The first equality is obtained by comparing the long exact sequences of
$\Ext$ groups associated to the residue exact sequences of $(U_i,
D\cap U_i)$ and $(U_i^{\ast}, D_i^{\ast})$.  

The second equality is obtained by comparing the local-to-global
spectral sequence of $\EExt$ sheaves and $\Ext$ groups with the long exact
sequence of local cohomology at the place $\Sigma$ of the sheaf $\Theta_X(-\log D)$.
\end{proof}
Let $V$ be the germ of an isolated rational singularity and let
$\widetilde{V}$ be a resolution of $V$. 
\begin{dfn}\cite{MR1489117}
\label{dfn:11}
\[\mu(V)= \dim_{\C}\Coker[\frac{1}{2i\pi}d\log\colon
H^1(\widetilde{V}, \mathcal{O}_{\widetilde{V}}^{\ast})\otimes \C \to H^1(\widetilde{V}, \Omega^1_{\widetilde{V}}) ] \]  
\end{dfn}
Namikawa shows \cite{MR1358982} that $\mu(V)$ is independent of
the chosen resolution.
\begin{lem}\cite{MR1358982}
\label{lem:20}
Let $V$ be a terminal Gorenstein singularity of dimension $3$. Then
$\mu(V)=0$ if and only if $V$ is non-singular or if $V$ is an ordinary
double point.   
\end{lem}
Let $f \colon \widetilde{V} \to V$ be a good resolution of the germ of
a terminal Gorenstein singularity and let $D\subset V$ be a Cartier divisor
with no worse than a rational double point at the singular point. Let $F= \widetilde{D} \cup E$,
where $\widetilde{D}$ is the proper transform of $D$ and $E$ is the
exceptional divisor of $f$. The divisor $F$ has simple normal
crossings and $\widetilde{D}$ is a resolution of a sufficiently small open
neighbourhood of a rational double point on $D$. 
\begin{lem}\cite{MR1489117}
\label{lem:22}
Using the above notation, the following hold:
\begin{enumerate}
\item $H^1(F, \mathcal{O}_F^{\ast})\simeq H^1(E,
  \mathcal{O}_E^{\ast})$.
\item $H^1(F, \widehat{\Omega}^1_F)\simeq H^1(E, \widehat{\Omega}^1_E)$.
\item In the natural commutative diagram
\[
\xymatrix{
H^1(\widetilde{V}, \Omega^1_{\widetilde{V}})\ar[r]^{\sigma} & H^1(F,
\widehat{\Omega}^1_F) \\
H^1(\widetilde{V}, \mathcal{O}^{\ast}_{\widetilde{V}})\otimes_{\Z}\C \ar[u]^{\phi}
\ar[r]^{\sigma'} & H^1(F, \mathcal{O}^{\ast}_F)\otimes_{\Z}\C \ar[u]_{\frac{1}{2i\pi}d\log},
}
\]
$\sigma'$ is an isomorphism and $\frac{1}{2i\pi}d\log$ is a surjection.
\end{enumerate}  
\end{lem}
\begin{proof}[Sketch proof of Lemma~\ref{lem:22}]
This Lemma follows from computations similar to the ones detailed in
the Appendix~\ref{sec:append-mixed-hodge}.
A base change diagram in the
fashion of \cite{MR972983} gives a resolution of  $\Q_V$ in terms of
$\Q_{\widetilde{V}}$ and $\Q_E$. The Mixed Hodge theory of $E$ and $F$
is then determined by a Mayer-Vietoris resolution. 
The computations of the
cohomology groups and properties of
Hodge numbers of generalised Fano $3$-folds yield the results. 
\end{proof}
\begin{step}
The invariants $\mu(P_i)$ can be used to identify a ``good'' direction
of deformations of the pair $(X,D)$.   
\end{step}
Let $f\colon \widetilde{X}\to X$ be a good
resolution of $X$ and let $F$ be the simple normal crossing divisor $E \cup
\widetilde{D}$.
\begin{cla}
There is an
injection \[\Theta_{\widetilde{X}}(-\log
F)=\HHom(\Omega^1_{\widetilde{X}}(\log F), \mathcal{O}_X)\to
\Omega^2_{\widetilde{X}}(\log F).\]
\end{cla}
Since the pair $(X,D)$ has log canonical singularities,
$K_{\widetilde{X}}+F$ is an effective divisor. The exact sequence
\[
0\to \Omega^3_{\widetilde{X}}(\log F)(-F)\to \Omega^3_{X}\to
\Omega^3_F\to 0
\]
shows that $\Omega^3_{\widetilde{X}}(\log F)\simeq
 \mathcal{O}_{\widetilde{X}}(K_{\widetilde{X}}+F)$ because
 $\Omega^3_F$ is trivial. The cup product map
 $\Omega^1_{\widetilde{X}}(\log F)\otimes\Omega^2_{\widetilde{X}}(\log
 F)\to \mathcal{O}_{\widetilde{X}}(K_{\widetilde{X}}+F)$ gives the
 desired injection.   
    
Recall that $\alpha_{\log}$ is defined as 
\[
\alpha_{\log} \colon \Ext^1(\Omega^1_X(\log D),
\mathcal{O}_X)\to  H^0(X, \EExt^1(\Omega^1_X(\log D), \mathcal{O}_X)).
\]
Lemma~\ref{lem:33} identifies the morphism $\alpha_{\log}$ with the coboundary map
of local cohomology \[H^1(U, \Theta_{X}(-\log D))\to H^2_{\Sing X}(X,
\Theta_X(-\log D)).\]

As ${\Omega^2_X(\log D)}_{\vert U} \simeq {\Theta_X(-\log D)}_{\vert
  U}$ and \[H^2_{\{P_i\}}(X, f_{\ast}\Omega^2_{\widetilde{X}}(\log
F))\simeq H^2_{\{P_i\}}(X, \Theta_X(-\log D)),\] the diagram
\[
\xymatrix{
H^1(f^{-1}(U), \Omega^2_{\widetilde{X}}(\log F)) \ar[r] & \oplus
H^2_{E_i}(\widetilde{X}, \Omega^2_{\widetilde{X}}(\log
F))\ar[r]^{\oplus \gamma_i} & H^2(\widetilde{X}, \Omega^2_{\widetilde{X}}(\log
F))\\
H^1(U,\Theta_X(-\log D))\ar[u]\ar[r]^{\alpha_{\log}} &
H^2_{\{p_i\}}(X, \Theta_X(-\log D))\ar[u]_{\oplus \tau_i}
}
\]
is commutative.
\begin{lem}\cite{MR1489117}
\label{lem:23}
Assume that $P_i\in X$ is neither non-singular nor an ordinary double
point; then $\gamma_i$ is not an injection. Moreover, $\dim \ker
\gamma_i\geq \mu(U_i)$.   
\end{lem}

\begin{proof}[Sketch proof of Lemma~\ref{lem:23}]
Set $V_i= f^{-1}(U_i)$ for $U_i$ a contractible Stein open
neighbourhood of $\{P_i\}$. 
Consider the dual map \[\gamma_i^{\ast} \colon H^1(\widetilde{X},
\Omega^1_{\widetilde{X}}(\log F)\otimes
\mathcal{O}_{\widetilde{X}}(-F))\to H^1(V_i, \Omega^1_{V_i}(\log F)\otimes
\mathcal{O}_{V_i}(-F)).\]
Namikawa shows \cite{MR1489117} that $\dim \Coker \gamma_i^{\ast}\geq
\mu(U_i)$ by studying the long exact sequences in cohomology on
$\widetilde{X}$ and $V_i$ associated to
\begin{eqnarray*}
0 \to \Omega^1_{\widetilde{X}}(\log
F)(-F)\to\Omega^1_{\widetilde{X}} \to \hat{\Omega}^1_F\to 0, \\  
0 \to \Omega^1_{V_i}(\log
F)(-F)\to\Omega^1_{V_i} \to \hat{\Omega}^1_{F\cap V_i}\to 0  
\end{eqnarray*}
and the commutative diagram:
\[
\xymatrix{
H^1(V_i, \Omega^1_{V_i})\ar[r]^{\sigma} & H^1(F_i, \hat{\Omega}^1_{F_i})\\
H^1(V_i, \mathcal{O}^{\ast}_{V_i})\otimes_{\Z}\C \ar[u]^{\phi}
\ar[r]^{\sigma'} & H^1(F, \mathcal{O}^{\ast}_{F_i})\otimes_{\Z}\C \ar[u]_{\frac{1}{2i\pi}d\log}.
}
\]  
\end{proof}

Let $\beta_i$ be the natural morphism from the space of global first order
deformations of $(V_i, F_i)$ to the space of local first order deformations
of $(U_i, D_i)$ defined as follows:
\begin{multline*}
\beta_i \colon H^1(V_i, \Theta_{V_i}(-\log F))\to
H^1(V_i \smallsetminus E_i,
\Theta_{V_i}(-\log F))\\ \simeq H^1(U_i\smallsetminus P_i, \Theta_{U_i}(-\log
D))
\simeq H^0(U_i, \EExt^1(\Omega^1_X (\log D), \mathcal{O}_X)).
\end{multline*}
\begin{lem}\cite{MR1489117}
\label{lem:24}
There is an element $\eta \in \Ext^1(\Omega^1_X, \mathcal{O}_X)$ such
that $\alpha_{\log}(\eta)_i$ is not in the image of $\beta_i$, for all
$i$ such that $P_i\in X$ is neither non-singular, nor an ordinary
double point.  
\end{lem}
\begin{proof}[Proof of Lemma~\ref{lem:24}]
By Lemma~\ref{lem:23}, if $P_i$ is neither non-singular nor an ordinary
double point, $\gamma_i$ is not injective.

There is an element
$\eta\in \Ext^1(\Omega^1_X(\log D), \mathcal{O}_X)$ such that $\tau_i
(\alpha_{\log}(\eta)_i)\neq 0$. The map $\tau_i$ factors as:
\[
H^2_{\{P_i\}}(X, \Theta_X(-\log D))\stackrel{\tau_i'}\to
H^2_{E_i}(\widetilde{X}, \Theta_{\widetilde{X}}(-\log F)) \to
H^2_{E_i}(\widetilde{X}, \Omega^2_{\widetilde{X}}(\log F)).
\]
In particular, $\tau_i'(\alpha_{\log}(\eta)_i)\neq 0$, and as the
spectral sequence of local cohomology reads
\[
H^0(X,R^1f_{\ast}\Theta_{\widetilde{X}}(-\log
F))\stackrel{\beta_i}\to H^2_{\{P_i\}}(X, \Theta_X (-\log
D))\stackrel{\tau_i'} \to H^2_{E_i}(\widetilde{X},
\Theta_{\widetilde{X}}(-\log F)),
\]
$\alpha_{\log}(\eta)_i$ is not in $\im \beta_i$.  
\end{proof}
\begin{step}
After a finite number of deformations along distinguished directions,
$X$ becomes a generalised Fano $3$-fold with no worse than ordinary
double points.  
\end{step}
Let $g_i\colon (\mathcal{U}_i, \mathcal{D}_i) \to \Def(U_i, D_i)$ be
the universal family.
Let $g\colon \mathcal{X}\to S$ be a smooth morphism and let
$\mathcal{D}= \sum \mathcal{D}_j$ be a divisor of $\mathcal{X}$ with simple normal
crossings.

\begin{dfn}
  \label{dfn:12}
The morphism $g\colon (\mathcal{X}, \mathcal{D})\to S$ is \emph{log
  smooth} if 
\begin{enumerate}
\item $\mathcal{D}_t= \sum \mathcal{D}_{j,t}$ is a divisor of $\mathcal{X}_t$ with simple
  normal crossings for each $t\in S$ and
\item for any point $p \in \mathcal{X}$, 
$g$ is locally a trivial deformation of
$(\mathcal{X}_t,\mathcal{D}_t)$ around $g(p)=t$. 
\end{enumerate}
\end{dfn}

\begin{rem}
\label{rem:31}
Notice that $g_i$ is log smooth over a non-empty Zariski open subset
$S_i^0\subset \Def(U_i, D_i)$ by Sard's theorem.
\end{rem}

Namikawa constructs iteratively a stratification of $\Def(U_i, D_i)$ into locally
closed non-singular subsets with the following properties:
\begin{enumerate}
\item 
$S^i_0 \subset \Def(U_i, D_i)$ is a Zariski open subset and $g_i$ is
log smooth over $S^i_0$.
\item
$S^i_k$ is a locally closed non-singular subset of pure codimension,
and $\codim(S_i^k, \Def(U_i, D_i))$ is strictly increasing with $k$.
\item
If $k>l$, $\overline{S^i_k}\cap S^i_l= \emptyset$.
\item $(\mathcal{U}_i, \mathcal{D}_i)$ has a simultaneous resolution
  over each $S^i_k$, that is, there exists $\mathcal{\nu}_{k}^i\colon (\mathcal{V}_i,
  \mathcal{F}_i)\to(\mathcal{U}_i,
  \mathcal{D}_i)$ a resolution over $S^i_k$, such that $g_i^k\circ
  \nu_i^k \colon (\mathcal{V}_i,\mathcal{F}_i)\to S_i^k$ is log
  smooth, where $g_i^k\colon (\mathcal{U}_i, \mathcal{D}_i)\to S_i^k$
  is the base change of $g_i$ to $S_i^k$.
\end{enumerate}
\begin{step}
 Let $X$ be a generalised Fano $3$-fold with terminal Gorenstein
 singularities. Then $X$ can be deformed to a non-singular generalised
 Fano $3$-fold.
\end{step}
Fix a stratification as above for each $P_i \in \Sing(X)$. 

Let $q_i\in
\Def(U_i, D_i)$ be the point corresponding to $(U_i, D_i)/k$ and let
$S_i^k$ be the stratum containing $q_i$. Let $\nu_i \colon V_i\to U_i$ be
the resolution induced by the log smooth simultaneous resolution of
$(\mathcal{U}_i, \mathcal{D}_i)$ over $S^i_k$. Since $\nu_i$ is an
isomorphism above $U_i \smallsetminus P_i$, the resolutions $\nu_i$ can be patched to a
global resolution $\nu\colon \widetilde{X}\to X$. Consider the divisor
$F= E \cup \widetilde{D}$, where $E$ is the exceptional divisor of
$\nu$ and $\widetilde{D}$ the proper transform of $D$. 

Pick an $\eta \in \Ext^1(\Omega^1_X(\log
D), \mathcal{O}_X)$ as in Lemma~\ref{lem:24}: since $\Def(X,D)$ is
non-singular and universal, $\eta$ determines a small deformation $g \colon
(\mathcal{X}, \mathcal{D}) \to \Delta^1_{\epsilon}$. There is a
holomorphic map $\phi_i\colon \Delta^1_{\epsilon}\to \Def(U_i, D_i)$ for
each $i$, with $\phi_i(0)=q_i$. By definition of $\eta$, if $P_i$ is
neither non-singular, nor an ordinary double point, $\im \phi_i$ is not
contained in $S_i^k$. 

Pick a suitable $t\in \Delta^1_{\epsilon}\smallsetminus \{0\}$ such that
$\phi_i(t)\in S^{k'}_i$ for some $k'<k$. Apply the same procedure
to $(\mathcal{X}_t, \mathcal{D}_t)$. This is possible because
$\Def(U_i, D_i)$ is versal at every point near $q_i$. After a finite
number of iterations, $X$ becomes a generalised Fano $3$-fold with no
worse than ordinary double points and the result follows from Lemma~\ref{lem:21}.
\end{proof}
\begin{rem}
  \label{rem:29}
Namikawa's proof shows that if $X$ is a generalised Fano $3$-fold with
no worse than ordinary double points, there is a flat (global) smoothing $f\colon
\mathcal{X}\to \Delta$.
If $X$ is a generalised Fano $3$-fold, there is a $1$-parameter flat
deformation $f \colon \mathcal{X}\to \Delta$, such that
$\mathcal{X}_t$ is a non-singular Fano $3$-fold for some $t\in
\Delta$. The construction of the small deformation shows that
$\mathcal{X}_t$ has terminal Gorenstein singularities for all $t\in \Delta$.
\end{rem}

The total space of the $1$-parameter flat deformation satisfies the following
additional properties. 
\begin{lem}\cite{jarad}
\label{lem:18} Let $\mathcal{X}\to \Delta$ be a one-parameter
smoothing of $X$. The total space $\mathcal{X}$ is normal, parafactorial and has at
most isolated Gorenstein singularities. 
\end{lem}
\begin{proof}
A point $P \in \mathcal{X}$ belongs to a fibre $\mathcal{X}_t$, that is some reduced
Cartier divisor $\mathcal{X}_t$ with no worse than terminal Gorenstein
singularities. Denote by $\mathcal{I}_t$ the ideal defining
$\mathcal{X}_t$ in $\mathcal{X}$; the conormal sequence
\begin{equation*}
\mathcal{I}_t/\mathcal{I}_t^2\to \Omega^1_{\mathcal{X}\vert \mathcal{X}_t}
\to \Omega^1_{\mathcal{X}_t}\to 0 
\end{equation*}
shows that the local embedding dimension of $\mathcal{X}$
at $P$ is $4$ or $5$. In particular, $P$ is either an isolated analytic
hypersurface singularity  or a non-singular point and $\mathcal{X}$ is
Cohen Macaulay, normal and Gorenstein. By inversion of
adjunction (\cite[Theorem 5.50]{MR1658959}), the pair $(\mathcal{X},
X)$ is purely log terminal
because $(X,0)$ is
klt. The variety $\mathcal{X}$ is therefore terminal.

As $\mathcal{X}$ is a locally complete intersection, $\mathcal{X}$ is
parafactorial by \cite[XI, Th\'eor\`eme 3.13]{MR2171939}.   
\end{proof}

\begin{rem}
  \label{rem:5}
Let $\pi \colon \mathcal{X} \to
\Delta$ be a one parameter flat deformation of a generalised Fano
$3$-fold $X= \pi^{-1}(0)$. Then $\mathcal{X}_t$ is a generalised
Fano $3$-fold for all $t$ and has the same Picard rank as $X$. 
\end{rem}
\begin{proof}[Proof of Remark~\ref{rem:5}]

This result is well known, I include a sketch of proof for convenience. 
As the total space $\mathcal{X}$ is Goresntein, $-K_{\pi}$ is Cartier
and for all $t \in \Delta$,
$\mathcal{X}_t$ is Gorenstein, so that the relative anticanonical
divisor $-K_{\pi\vert \mathcal{X}_t}= -K_{\mathcal{X}_t}$ is
Cartier.

The
divisor $-K_X= -K_{\pi \vert \mathcal{X}_0 }$ is nef, and so is
$-K_{\pi \vert \mathcal{X}_t }$ for all nearby $t$. Indeed,
Proposition~\ref{pro:2} shows that $\pi$ induces a $1$-parameter
deformation $\mathcal{Y}\to \Delta$ of $Y$. The Kuranishi spaces
$\Def(X)$ and $\Def(Y)$ are smooth complex analytic spaces; the proper
morphism $X\to Y$ is small, hence $\Def(X)\to \Def(Y)$ is a closed
immersion. For all $t\in \Delta$,
the map $\mathcal{X}_t\to \mathcal{Y}_t$ is small and contracts
$K$-trivial curves by Theorem~\ref{thm:6}. 
For all $t\in
\Delta$, $\mathcal{Y}_t$ is a terminal Gorenstein Fano $3$-fold,
because ampleness is an open condition on flat families;
$\mathcal{X}_t$ is a generalised Fano $3$-fold.  
Lemma~\ref{lem:27}
establishes that the Picard rank is constant in a smoothing of a
terminal Gorenstein generalised Fano $3$-fold. The exact same
arguments can be used in the case of the explicitly constructed
smoothing $f\colon \mathcal{X}\to \Delta$. Recall that each fibre of
$f$ has isolated hypersurface singularities. One can define vanishing
cycles $B_{i,t}$ whose cohomology is supported in degrees $0$ and $3$ and
construct a homeomorphism between $X \smallsetminus \{P_1, \cdots ,
P_n\}$ and $X_t\smallsetminus \cup B_{i,t}$. 

The Picard rank is also constant in the $1$-parameter flat deformation
of the anticanonical model $\mathcal{Y}\to \Delta$ induced by $\pi$
and for all $t\in \Delta$, $\mathcal{Y}_t$ is a terminal Gorenstein
Fano $3$-fold. The degree $-K_{\mathcal{Y}_t}^3$ is invariant in a flat
deformation, because the plurigenera are.
\end{proof}

\subsection{Deformations of extremal contractions}
\label{sec:deform-birat-oper}
Let $X$ be a generalised Fano $3$-fold and $Y$ its anticanonical
model. Denote by $f$ the anticanonical map of $X$. If $E$ is a
$\Q$-Cartier divisor of $X$ such that $\overline{E}=f(E)$ is not
$\Q$-Cartier on $Y$, denote by $Z$ the symbolic blow-up of
$\overline{E}$ on $Y$. The map $f$
factors through $Z$ in the following way:
\[
f \colon X  \stackrel{h}\to Z \stackrel{g}\to Y.
\]
In particular, when $\Pic(X/Y)=1$, $h$ is the identity. I recall some
known results about the relations between deformations of $X,Z$ and
$Y$.

\begin{pro}\cite{MR1149195}
\label{pro:2}
Let $X$ be a projective $3$-fold and let $f \colon X \to X'$ be a
proper map with connected fibres. Assume that $R^1f_{\ast}\mathcal{O}_X=0$.
Then, there exist natural morphisms $F$ and $\mathcal{F}$ that make the
  following diagram commutative:
\[
\xymatrix{
\mathcal {X} \ar[r]^{\mathcal{F}}  \ar[d] & \mathcal{X'} \ar[d]\\
\Def(X) \ar[r]^{F} & \Def(X')
}
\]
where $\mathcal{X}$ (resp. $\mathcal{X'}$) is the Kuranishi family of
$X$ (resp. of $X'$) and $\Def(X)$ (resp. $\Def(X')$) is the Kuranishi
space of $X$ (resp. of $X'$); $\mathcal{F}_{\vert X}$
coincides with $f$.
\end{pro}
In particular, there exist maps $\mathcal{G}$ and $\mathcal{H}$ that
restrict to $g$ and $h$ on the central fibre and make the following
diagram commutative:
\[
\xymatrix{
\mathcal{X}\ar[r]^{\mathcal{H}}\ar[d] &
\mathcal{Z}\ar[r]^{\mathcal{G}}\ar[d]& \mathcal{Y}\ar[d] \\
\Def(X)\ar[r]^{H} & \Def(Z)\ar[r]^{G} & \Def(Y)
.}
\]

\begin{thm}\cite{MR1149195}
\label{thm:14}
  Let $Y$ be a projective $3$-fold with canonical singularities and let $f
  \colon X \to Y$ be a projective, $\Q$-factorial terminal and crepant
  partial resolution. Then $F \colon \Def(X) \to \Def(Y)$ is a finite map.
\end{thm}
\begin{proof}[Sketch of proof]
Let $T= \Spec \C[[t]]$ be the formal disc. 
The proof relies on showing that no non-trivial deformation of $X$
over $T$ can induce a trivial deformation of $Y$ over $T$. Assume that
$F$ is not finite. There exists a deformation $\mathcal{X}/T$ such
that $F\colon \mathcal{X} \to Y\times T$. In particular, there is a
birational map $g \colon \mathcal{X} \dashrightarrow X\times T$. The map
$g$ is induced by a map $g_0$ on $X$ and
$\mathcal{X}$ is also a trivial deformation.  
\end{proof}
\begin{pro}\cite{MR1149195}
\label{pro:6}  
Let $Y$ be a terminal Gorenstein projective $3$-fold and $X \to Y$ a
small $\Q$-factorialisation.
The subspace $\im [\Def(X) \to \Def(Y)]$ is closed and independent of
the choice of small $\Q$-factorialisation.
\end{pro}
\begin{rem}
\label{rem:30}
Let $\Theta_Y$ (resp. $\Theta_Y$) be the dual of $\Omega^1_Y$
(resp. of $\Omega^1_X$) and let $C$ be the exceptional locus of $X\to
Y$. There is a natural map of functors
$\Def(X)\to \Def(Y)$ because $Y$ has terminal singularities and
$R^1\pi_{\ast}\mathcal{O}_Y=(0)$. At the level of tangent spaces, the
kernel is given by the local cohomology group
$H^1_{C}(\Theta_{X})$. This group vanishes when $C$ is a curve \cite{MR848512}.   
\end{rem} 
The following result shows that $\Q$-factoriality is an
open condition on the base of flat deformations of algebraic
$3$-folds.

\begin{thm}[Factoriality and deformations, \cite{MR1149195}]
\label{thm:13}
 Let $g \colon X \to S$ be a flat family of algebraic varieties. Assume
 that the fibres have rational singularities and that for any $s \in
 S$, $\codim(\Sing(X_s), X_s) \geq 3$. For instance, this is the case of
 a flat family of $3$-folds with isolated singularities. Then
\[
X_{\Q-\mbox{fact}}=\{x \in X  \vert g^{-1}(g(x)) \mbox{ is $\Q$-factorial}\}
\]     
is open in $X$. In particular, if a fibre is $\Q$-factorial, so are
nearby fibres. 
\end{thm}
\begin{rem}
\label{rem:48}
  Let $f\colon \mathcal{Y}\to \Delta$ be a $1$-parameter
  flat and proper deformation of a generalised Fano $3$-fold $Y$ with
  $\mathcal{Y}_t$ non-singular for some $t \in \Delta$ as
  contructed by Namikawa (see Section~\ref{sec:deform-theory-gener}).
The set $\Delta_{\Q-\mbox{fact}}= \{t\in \Delta \vert \mathcal{X}_t
  \mbox{ is $\Q$-factorial }\}$ is open.
\end{rem}
I now recall some results on deformations of extremal rays.

\begin{thm}[Flops in families, \cite{MR1149195}]
\label{thm:6}
  Let $f_{0} \colon X_{0} \to Y_{0}$ be a proper morphism between
  normal $3$-folds. Assume that $X_{0}$ has only terminal
  singularities, and that $f_0$ contracts a curve $C_0 \subset X_0$ to
  a point $Q_0 \in Y_0$. Assume moreover that every component of $C_0$
  is $K_{X_0}$-trivial. Let $X_S \to S$ be a flat deformation of
  $X_0$ over the germ of a complex space $0 \in S$.
Then,
\begin{enumerate}
\item $f_0$ extends to a contraction morphism $F_S \colon X_S \to Y_S$.
\item The flop $F_S^+ \colon X_S^+ \to Y_S $ exists and commutes with
  any base change.
\end{enumerate} 
\end{thm}
\begin{proof}[Sketch proof]
This follows from Theorem~\ref{thm:14}. For each $s \in S$, the singularities of $Y_s$ are
terminal. Using a local analytic description of $Y_s$ near a singular
point, the flop $F_S\colon
X_S^{+} \to Y_S$ can be explicitly constructed as in \cite{MR1144527}.
\end{proof}

\begin{thm}[Deformation of extremal rays, \cite{MR1149195}]
\label{thm:15}
 Let $g \colon X \to S$ be a proper flat morphism of projective
 varieties. Assume that for some $0 \in S$ the fibre $X_0$ is a
 $3$-fold with no worse than  $\Q$-factorial canonical singularities.
Let $f \colon X_0 \to X'_0$ be
 the contraction of an extremal ray $C_0 \subset X_0$. There then
 exist a proper flat morphism $h \colon X' \to S$ (by
 Proposition~\ref{pro:2}) and a factorisation $g \colon X \stackrel{f}
 \to X' \stackrel{h}\to S$. Moreover, there exists an open
 neighbourhood $U$ of $0 \in S$ such that:
 \begin{description}
 \item[(i)]if $f_0$ contracts a subset
 of $\codim \geq 2$ (resp. a divisor, resp. is a fibre space of generic
 relative dimension $k$), $f_s$ contracts a subset
 of $\codim \geq 2$ (resp. a divisor, resp. is a fibre space of generic
 relative dimension $k$) for all $s\in U$, 
\item[(ii)] $f_s$ is the
 contraction of an extremal ray.  
 \end{description}
\end{thm}

\begin{rem}
Lemma~\ref{lem:29} describes
contractions of extremal rays with Cartier exceptional divisor on
small partial $\Q$-factorialisations
of terminal Gorenstein Fano $3$-folds. This Lemma is a mild generalisation
of Cutkosky's results (Theorem~\ref{thm:8}), and shows that the same
classification of extremal rays essentially applies as in the $\Q$-factorial
terminal Gorenstein case.
   
\end{rem}
\begin{lem}
\label{lem:9}
Let $Z$ be a small partial
$\Q$-factorialisation of a terminal Gorenstein Fano $3$-fold
$Y$. Denote by $\mathcal{Z} \to S$  a proper flat $1$-parameter
deformation of $X$ over the spectrum of a complete discrete
valuation ring $S= \Spec \mathcal{O}_S$ with residue field $\C$. Let
$0$ be the closed point of $S$ and $\eta$ the generic point. Assume
that each fibre has terminal Gorenstein singularities. 
If $f_{0} \colon Z \to Z'$ is a divisorial
extremal contraction with $\Q$-Cartier exceptional divisor $E$, there
is an $S$-morphism $f \colon \mathcal{Z} \stackrel{f} \to \mathcal{Z}'$ to a projective
 $1$-parameter flat deformation of $Z'$. 
The restrictions $f_{\eta}$ and $f_0$ are extremal
 divisorial contractions. In the notation of Theorem~\ref{thm:8},
 either the contraction $f_{\eta}$ is of the same type as $f_0$, or $f_{\eta}$
 and $f_0$ are of types E$3$ and E$4$.
\end{lem}
 
\begin{proof}
I give an outline of the argument, as the assumption that $E$ is
Cartier ensures that it can be proved as in
\cite[Proposition 3.47]{MR662120}.
\setcounter{case}{0}
If $\mathcal{D}_{\eta}$ is a Cartier divisor of $\mathcal{Z}_{\eta}$, the completion
$\mathcal{D}_S$ to $S$ is an irreducible divisor which is $\Q$-Cartier
because $\mathcal{Z}$ is  parafactorial: any
proper Weil divisor on $S$ that is Cartier outside of finitely many
fibres is $\Q$-Cartier. The
specialisation map associates to $\mathcal{D}_{\eta}$ the $\Q$-Cartier
divisor $\Red(\mathcal{D}_{\eta}) =\mathcal{D}_{S}\times_S \{0\}$ of
$Z$. As $Z$ is Gorenstein, $\Red(\mathcal{D}_{\eta})$ is in fact
Cartier \cite[Lemma 6.3]{MR924674}; $\Red$ is an injective
homomorphism  $\Red \colon NS(\mathcal{Z}_{\eta})\to NS(Z)$. This
homomorphism is bijective because the Picard rank is constant
in a $1$-parameter flat deformation of a $3$-fold with isolated
hypersurface singularities. 
Similarly, if $\mathcal{C}_{\eta}$ is an effective curve, denote by
$\Red(\mathcal{C}_{\eta})=\mathcal{C}_{S}\times_S \{0\}$ the
specialisation of $\mathcal{C}_{\eta}$. The $1$-cycle
$\Red(\mathcal{C}_{\eta})$ has non-negative coefficients and $\Red$
defines an injective homomorphism $\Red \colon
\overline{NE}(\mathcal{Z}_{\eta})\to \overline{NE}(Z)$ on the cone of
effective curves. In addition, if $\mathcal{D}_{\eta}$ is a Cartier
divisor and $\mathcal{C}_{\eta}$ is a $1$-cycle,
$\mathcal{D}_{\eta}\cdot \mathcal{C}_{\eta}= \Red(\mathcal{D}_{\eta}
\cdot \mathcal{C}_{\eta})=\Red(\mathcal{D}_{\eta})
\cdot \Red(\mathcal{C}_{\eta})$ \cite{MR1644323}. 

Recall from Remark~\ref{rem:5} that for all $t\in S$, $\mathcal{Z}_t$ is a
generalised Fano $3$-fold, in particular, $NE(\mathcal{Z}_{\eta})$ is rational polyhedral
(Lemma~\ref{lem:1}). 

\begin{case} The contraction $f_0\colon Z \to Y$ is of type E$1$,
  i.e. $f_0$ contracts a surface $E$ to a curve $\Gamma \subset Y$.
\end{case}

Denote by $l$ a general fibre of the contraction $E\to
\Gamma$. Lemma~\ref{lem:29} shows that $E$ is a $\PS^1$-bundle over
$\Gamma$, $l\simeq \PS^1$ and that
$-K_Z\cdot l=-E\cdot l=1$. The curve $l$ lies on a Cartier divisor, it is a
l.c.i. variety and by adjunction on $E$, $\mathcal{N}_{\Gamma/ Z}=
\mathcal{O}_{\PS^1}(-1)\oplus \mathcal{O}_{\PS^1}$. 
Let $\mathcal{H}$ be the connected component of the Hilbert prescheme
$\Hilb_{Z/S}$ containing the class $[l]$ of $l$. The Zariski tangent
space to $\Hilb_{Z/S}$ at $[l]$ is canonically isomorphic to $H^0(Z,
\mathcal{N}_{l/Z})$ and if $H^1(Z, \mathcal{N}_{l/Z})=(0)$,
$\mathcal{H}_{0}$ is non-singular
over $S$ at $[l]$ \cite[IV, Corollaire 5.4]{MR0146040}. The component
$\mathcal{H}_0$ is non-singular at $[l]$ ($h^1(Z, \mathcal{N}_{l/Z})=0$)
and has relative dimension $1$ at $0$ ($h^0(Z,
\mathcal{N}_{l/Z})=1$). The base $S$ is complete, and hence the
component $\mathcal{H}_0$ is connected. As this is true for all but a finite number of
fibres, the component $\mathcal{H}_0$ is isomorphic to
$\Gamma$. Denote by $\mathcal{D}\subset \mathcal{Z}\times_S
\mathcal{H}$ the universal closed subscheme. The projection of
$\mathcal{D}$ to $\mathcal{Z}$ induces an embedding
$\mathcal{D}_0\simeq E$. The component $\mathcal{H}$ is non-singular
over $S$ at any point of $\mathcal{H}_0$, hence it is smooth over $S$
and $\mathcal{H}_{\eta}$ is a  non-singular projective curve.

The universal subscheme $\mathcal{D}$ is by definition flat over
$\mathcal{H}$. As $\mathcal{D}_0$ is a $\PS^1$-bundle over
$\mathcal{H}_0\simeq \Gamma$, $\mathcal{D}_{\eta}$ is also a $\PS^1$-bundle over
$\mathcal{H}_{\eta}$. Denote by $\mathcal{C}_{\eta}$ any fibre of
$\mathcal{D}_{\eta} \to \mathcal{C}_{\eta}$; $\Red(\mathcal{C}_\eta)$
is a fibre $l$ of $E\to \Gamma$; it therefore generates an extremal ray $R_0$
of
$NE(Z)\supset \Red(NE(\mathcal{Z_\eta}))$. The curve
$\mathcal{C}_{\eta}$ generates an extremal ray $R_{\eta}$ of
$\mathcal{Z}_{\eta}$. This extremal ray may be contracted because
$\mathcal{Z}_{\eta}$ is a weak Fano $3$-fold, and since
$\mathcal{D}_{\eta}\cdot \mathcal{C}_{\eta}= D\cdot l=-1$, the divisor
$\mathcal{D}_{\eta}$ is exceptional for the contraction $cont_{R_{\eta}}$. If
$\mathcal{G}_{\eta}\subset \mathcal{D}_{\eta}$ denotes a curve that
dominates $\mathcal{H}_{\eta}$ such that $\Red(\mathcal{G}_{\eta})$
dominates $\mathcal{H}_{0}$, $[\mathcal{G}_{\eta}]$ does not belong to
the class $R_{\eta}$ because $\Red(\mathcal{G}_{\eta})$ does not belong to
$R_0$. This proves that the contraction of $R_{\eta}$ is an extremal contraction of type E$1$.    
\begin{case}
 The contraction $f_0\colon Z\to Y$ is of type E$2$-E$5$,
 i.e. $f_0$ contracts a surface $E$ to a point $P\in Y_0$. 
\end{case}
The divisor $E$ is a l.c.i. subscheme because $E$
is Cartier. The cohomology groups $H^0(E, \mathcal{O}_E(E))$ and
$H^1(E, \mathcal{O}_E(E))$ are both trivial (Lemma~\ref{lem:29}),
hence the connected component $\mathcal{H}$ of the Hilbert scheme
$\Hilb_{Z/S}$ containing $[E]$ has relative dimension $0$ at $[E]$ and
$\mathcal{H}_0$ is non-singular at $[E]$. The component $\mathcal{H}$ therefore
is isomorphic to $S$. Denote by $\mathcal{D}\subset
\mathcal{Z}\times_S S$ the universal
subscheme; the projection to $\mathcal{Z}$ determines an embedding
$\mathcal{D}_0\simeq E$. As in case $1$, any effective curve lying on
$\mathcal{D}_{\eta}$ belongs to the class of an extremal ray
$R_{\eta}$, that may be contracted. The contraction $cont_{R_{\eta}}$
contracts a divisor to a point. As $\mathcal{D}_{\eta}^3=
\Red(\mathcal{D}_{\eta})^3= E^3$, the contractions $f_0$ and $cont_{R_{\eta}}$
are of the same type, except possibly if one of them is of type E$3$
and the other is of type E$4$. 
\begin{cla}
  There is an $S$-morphism $f\colon \mathcal{Z} \to \mathcal{Z'}$ that
  restricts to $f_0$ on the central fibre and to $cont_{R_{\eta}}$ on
  the generic fibre.
\end{cla}
The Picard rank is constant in the $1$-parameter flat deformation and
$\Pic(\mathcal{Z})\simeq \Pic(Z) \simeq \Pic(\mathcal{Z}_{\eta})$,
since each fiber has terminal Gorenstein singularities.
Denote by $L$ a nef invertible sheaf on $Z$ associated to $f_0$,
i.e. such that $L^{\bot} \cap NE(Z)= R_0$ and $f_0$ is the morphism
determined by $\vert mL \vert$ for $m >\!\!> 0$. The invertible sheaf
$L$ comes from an invertible sheaf $\mathcal{L}$ on $\mathcal{Z}$ and
$\mathcal{L}_{\eta}^{\bot}\cap NE(\mathcal{Z}_{\eta})= R_{\eta}$ by
construction. The Cartier divisor $\vert m \mathcal{L} \vert$
determines the desired $S$-morphism $\mathcal{Z}\to \mathcal{Z'}$,
because by the base change theorem, $H^{0}(\mathcal{Z}, m
\mathcal{L}) \otimes \C\simeq H^0(Z, mL)$. 
\end{proof}

\setcounter{case}{0}
Every divisorial step of the Minimal Model Program on a
weak$\ast$ Fano $3$-fold can be deformed to an extremal divisorial
contraction on a non-singular weak$\ast$ Fano $3$-fold with Picard
rank $2$.

Assume that the extremal contraction $\phi_i \colon X_i \to X_{i+1}$
is divisorial and denote by $E_i$ its exceptional divisor. I use the notation
introduced in the proof of Lemma~\ref{lem:12}. The generalised Fano
$3$-fold $Z_i$ is a crepant blow up of $Y_i$  along the image of the
exceptional divisor $E_i$ of $\phi_i$ by
the anticanonical map. More precisely, $Z_i$ is defined as
$\underline{\proj} \bigoplus_n f_{i\ast} \mathcal{O}(nE_i)$ and fits in the diagram:
\[
\xymatrix{
X_i \ar[d] \ar[r]^{\phi_i} & X_{i+1} \ar[d]\\
Z_i \ar[d] \ar[r]^{\psi_i} & Z_{i+1}\\
Y_i &
}
\]    
where $Z_{i+1}=Y_{i+1}$ is the anticanonical model of $X_{i+1}$. 
The following corollary is a direct consequence of
Theorem~\ref{thm:16} and of Lemma~\ref{lem:9}.
\begin{cor}
 \label{cor:9}
There is a proper $1$-parameter flat deformation $\mathcal{Z}_i$ of
$Z_i$ over $\Delta$ such that $\mathcal{Z}_{i,t}$ is a terminal
Gorenstein generalised Fano with Picard rank $2$ for all $t\in \Delta$
and $\mathcal{Z}_{i,t_0}$ is non-singular for some $t_0\in
\Delta$. There is a morphism
$\Psi_i\colon \mathcal{Z}_i \to \mathcal{Z}_{i+1}$ which restricts to
$\psi_i$ on the central fibre, and such that $\mathcal{Z}_{i+1}$ is a
$1$-parameter flat deformation of $Z_{i+1}$. In a neighbourhood of
$0\in \Delta$, $(\mathcal{Z}_{i})_{t}$  is a Picard rank $2$ terminal
Gorenstein weak Fano $3$-fold and $(\mathcal{Z}_{i+1})_{t}$ is a
Picard rank $1$ terminal Gorenstein Fano $3$-fold. The morphism $\Psi_t$ restricts to the contraction
of an extremal ray of the same type as ${\Psi_i}_0$,
except if ${\Psi_i}_{0}$ is of type E$3$ and
${\Psi_i}_{\eta}$ is of type E$4$.  
\end{cor}
\begin{rem}
The same holds for an extremal contraction of fibering type. Note
that, in a
neighbourhood of $0$, $\mathcal{Z}_{i+1,t}$ is projective because
$Z_{i+1}$ is.    
\end{rem}

\section{Takeuchi game}
\label{sec:takeuchi-game}

In this Section, I generalise some constructions studied by Takeuchi
\cite{MR1022045} to determine explicitly the extremal
divisorial contractions encountered when running the Minimal Model
Program on a weak$\ast$ Fano $3$-fold $X$. 
Iskovskikh obtained a classification of Fano $3$-folds with Picard
rank $1$ by studying the double projection of Fano
$3$-folds from lines lying on them \cite{MR463151,MR503430}. 
Takeuchi investigated
projections of Fano $3$-folds from a
point or from a general conic \cite{MR1022045}. Both
constructions are special cases of Sarkisov elementary links
\cite{MR1311348}. 
Simple numerical calculations
based on the theory of extremal rays led to a considerable refinement of
Iskovskikh's methods. I generalise Takeuchi's construction in order to
study some projections of a terminal Gorenstein Fano $3$-fold $Y$ from 
curves lying on it. 

Let $Y$ be a terminal Gorenstein Fano $3$-fold with Picard rank $1$
that does not contain a plane
and let $Z$ be a small partial $\Q$-factorialisation of $Y$. I assume
that the Picard rank of $Z$ is $2$, and that there exists on $Z$ an
extremal contraction $\phi$ with Cartier exceptional divisor. 
In
Sections~\ref{sec:elem-contr-gorenst} and \ref{sec:gener-take-constr},
I study a generalised Takeuchi
construction on $Z$. This analysis yields systems of Diophantine equations, whose
solutions describe the possible divisorial contractions encountered
when the MMP is run on $X$, a small $\Q$-factorialisation of $Y$. This approach provides a theoretical
method to construct explicit examples of non $\Q$-factorial Gorenstein
terminal Fano $3$-folds.

In Section~\ref{sec:geom-motiv-non}, I use these techniques to give a
``geometric motivation'' of $\Q$-factoriality for terminal Gorenstein
quartic $3$-folds. If $Y_4^3\subset \PS^4$ fails to be
$\Q$-factorial, by definition, $Y$ contains a surface
$\overline{E}$, which is a Weil, non $\Q$-Cartier divisor. I show that
the surface $\overline{E}$ is a plane, a quadric or is one of the
surfaces listed in the table on page \pageref{table}. 

\subsection[Elementary contractions]{Elementary contractions of terminal Gorenstein Fano $3$-folds}
\label{sec:elem-contr-gorenst}
Let $X$ be a projective non-singular $3$-fold such that $K_X$ is not
nef. Mori shows \cite{MR662120} that
$X$ admits an elementary contraction morphism $\phi$ that corresponds to an
extremal ray of $\overline{NE}(X)$; he classifies all such
contractions and shows that either $\phi$ is of fibering type, or
$\phi$ contracts a Cartier divisor $E$ to a point or a curve. In the case of
non-singular surfaces, the exceptional divisor of an extremal
contraction is a $(-1)$-curve; Mori shows that in the $3$-fold case,
$E$ either is a plane or a quadric with anti-ample normal bundle or
$E$ is a $\PS^1$-bundle over a non-singular curve. Up to minor
generalisations, the classification of divisorial extremal
contractions still holds when $X$ is assumed to have
terminal Gorenstein $\Q$-factorial singularities \cite{MR936328}. I
study extremal contractions with Cartier exceptional divisor on small
partial $\Q$-factorialisations of terminal Gorenstein Fano $3$-folds.   

\begin{setup}
Let $Z$ be a normal, terminal Gorenstein weak Fano $3$-fold and let
$Y$ be its anticanonical model. Assume that the anticanonical map $h
\colon Z \to Y$ is small, that the anticanonical ring of $Z$ is
generated in degree $1$ and that $Y$ has Picard rank $1$.

Recall from Theorem~\ref{thm:1} that this is the case unless $Z$ is
monogonal or $Y$ is birational to a special complete
 intersection $X_{2,6}\subset \PS(1^4, 2,3)$ with a node. This is
 always true if $Y$ has Picard rank $1$ and genus at least $3$.

Let $\psi \colon Z \to Z_1$ be a $K_Z$-negative extremal contraction.
Since a flipping curve $\gamma$
on a terminal variety $Z$ satisfies $-K_Z \cdot \gamma <1$
\cite{MR796252}, $\psi$ is not an isomorphism in codimension $1$.
Assume that the exceptional divisor $E$ of the
contraction $\psi$ is $\Q$-Cartier. Note that the Weil divisor
$\overline{E}= h(E)$ is not $\Q$-Cartier. If it were, as $Y$ has
Picard rank $1$, $E$ would have to be ample. Note that as $Z$ has
terminal Gorenstein singularities, a divisor is $\Q$-Cartier if and
only if it
is Cartier \cite[Lemma 6.3]{MR924674}.     

Let
$g \colon X \to Z$ be a small $\Q$-factorialisation and let $\widetilde{E}$ be
the pull back of $E$ on $X$. There is a $K_X$-negative extremal
ray on which $\widetilde{E}$ is negative. Indeed, $g$ is small, hence
$K_{X}=  g^{\ast}(K_{Z})$ and $\widetilde{E}= g^{\ast}(E)$. The divisor $\widetilde{E}$ is
covered by curves $\Gamma$ such that $K_{X}\cdot \Gamma
<0$ and $\widetilde{E}\cdot \Gamma<0$, because $E$ has this property. The $3$-fold $X$ is a
weak Fano, by Lemma~\ref{lem:1}, there is a $K_X$-negative extremal
contraction $\phi$ on $X$ with exceptional divisor $\widetilde{E}$.

Denote by $f\colon X\to
Y$ the anticanonical map of $X$. As in the proof of
Lemma~\ref{lem:12}, the following diagram is commutative because
$\widetilde{E}$ is negative on an extremal ray of $NE(X/Z_1)$: 
\begin{eqnarray}
\label{eq:40}
\xymatrix{\widetilde{E} \subset X \ar[r]^{\phi} \ar[d]_g & X_1\ar[d]^{g_1}\\
E \subset Z \ar[d]_h \ar[r]^{\psi} & Z_1\\
\overline{E} \subset Y &
}
\end{eqnarray}
In the diagram, $g \circ h= f$, $\widetilde{E} =
f^{\ast}(\overline{E})$ and $g_1$ is an isomorphism in codimension $1$. 
\end{setup}
\begin{rem}
\label{rem:17}
  In the terminology
introduced in Section~\ref{sec:categ-weak-fano}, if $Y$ does not
contain a plane $\PS^2$ with $ {-K_Y}_{\vert \PS^2}=
\mathcal{O}_{\PS^2}(1)$, then $X$ is a weak$\ast$ Fano $3$-fold.
\end{rem}
\begin{lem}
 \label{lem:29}
 \begin{enumerate}
 \item[E$1$:]If $\psi$ contracts $E$ to a curve $\Gamma$, 
 $\Gamma$ is
  locally a complete intersection and has planar singularities. The
 contraction $\psi$ is locally the blow up of the ideal sheaf
 $I_{\Gamma}$. In the local ring $\mathcal{O}_{Z_1, P}$ of any point
$P \in \Gamma$, one of the local equations of $\Gamma$ is a smooth
hypersurface near $P$.
 \item[] If $\psi$ contracts $E$ to a point $P$, then one of the
 following holds: 
\begin{itemize}   
\item[E$2$:] $(E,\mathcal{O}_{E}(E))\simeq ( \PS^2,
  \mathcal{O}_{\PS^2}(-1)))$ and $P$ is a non-singular point.
\item[E$3$:]  $(E,\mathcal{O}_{E}(E))\simeq (\PS^1 \times \PS^1,
  \mathcal{O}_{\PS^1 \times \PS^1}(-1,-1))$. 
\item[E$4$:] $(E, \mathcal{O}_{E}(-E))\simeq (Q, \mathcal{O}_{E}
  \otimes \mathcal{O}_{\PS^3}(-1))$, with  $Q$ an irreducible reduced
  singular quadric surface in $\PS^3$. 
\item[E$5$:] $(E,\mathcal{O}_{E}(E))\simeq ( \PS^2,
  \mathcal{O}_{\PS^2}(-1)))$, and $P$ is a non-Gorenstein point of
  index $2$.
\end{itemize}
\end{enumerate}
\end{lem}
\begin{proof}
\setcounter{case}{0}
The diagram \eqref{eq:40} is commutative, therefore $g_1$ maps the centre
of the contraction $\phi$ to the centre of the contraction $\psi$. If
$\phi$ contracts a divisor to a point, so does $\psi$. The map
$g_1\colon X_1\to Z_1$ is a small $\Q$-factorialisation of $Z_1$,
so that if the centre of $\phi$ is a curve $\Gamma$ and that of $\psi$ is a
point $\{P\}$, $-K_{X_1}\cdot \Gamma=0$, the curve $\Gamma$ is non-singular, $-K_{X_1}$ is nef and big and
defines a small map by the proof of Theorem~\ref{thm:1}. As in the proof
of Lemma~\ref{lem:26}, $\widetilde{E} \simeq \F_2$ or
$\widetilde{E}\simeq \PS^1\times \PS^1$. The Cartier divisor $E=g_{\ast}(\widetilde{E})$ is an
irreducible reduced quadric in both cases, and $-K_{Z\vert E}=
\mathcal{O}_E(1)$. The contraction $\psi$ is of
type E$3$ or E$4$. 

I assume that the contractions $\phi$ and $\psi$ either both
contract a divisor to a curve or both contract a divisor to a point.
 
\begin{case}
The contraction $\psi \colon Z \to Z_1$ contracts the divisor $E$ to a
curve $\Gamma$.   
\end{case}

The $3$-fold $Z$ has terminal Gorenstein singularities, in particular
its singularities are isolated.
Let $l$ be a fibre of the contraction $\psi$ that contains no
singularities of $Z$. The map $g$ is small
and by construction of $\psi$, there exists a fibre $\tilde{l}$ of $\phi$
mapping $1$-to-$1$ to $l$. Then, $-K_Z\cdot l= g^{\ast}(-K_Z)\cdot
\tilde{l}=-K_X \cdot \tilde{l}= 1 $ and $E\cdot l= g^{\ast}(E)\cdot
\tilde{l}= \widetilde{E}\cdot\tilde{l}=1$.
The contraction $\psi$ is the blow up of a
non-singular curve away from finitely many points.

The anticanonical model $Y$ of $Z$ has Picard rank $1$ and genus $g
\geq 3$, so the weak Fano $3$-fold $Z$ itself has basepoint free
anticanonical system $\vert -K_Z\vert$. Let $S$ be a general section
of $\vert -K_Z \vert$; $S$ is non-singular by Theorem~\ref{thm:1}. As $\vert -K_Z \vert=
\vert \psi^{\ast}(-K_{Z_1})-E\vert$ (since $-K_{X}=
\phi^{\ast}(-K_{X_1})-\widetilde{E}$), $S$ is mapped to $S_1$, a section of
$\vert -K_{Z_1}\vert$ that contains the curve $\Gamma$. The
contraction $\phi$ is not of type E$5$, therefore $X_1$ and $Z_1$ are
Gorenstein and $S_1$ is Cartier. 
The morphism $\psi_{\vert S}\colon S \to S_1$ is an isomorphism away
from $\Gamma$. For any fibre $l$ of $\psi$, $l$ is not
contained in $S$ because $S$ is general in $\vert -K_Z \vert$, and $l$ intersects $S$ in a
finite number of points. Over an affine neighbourhood of $P=\psi(l)
\in \Gamma$,
the morphism $\psi_{\vert S}\colon S \to S_1$ is a finite birational
morphism, which is an isomorphism away from $l$. In particular, the
singularities of $S_1$ are isolated. Further, $S_1$ is normal because
$S_1$ is a Cartier divisor in $Z$, which is Cohen Macaulay, and $S_1$
is regular in codimension $1$. The morphism $\psi_{\vert S}$ is finite
birational, hence, by Zariski's main theorem, it is an
isomorphism. The curve $\Gamma$ lies on a non-singular Cartier
divisor: it is locally a complete intersection.
Hence, the curve $\Gamma$ has planar singularities.  

\begin{case} 
The contraction $\psi \colon Z \to Z_1$ contracts the divisor $E$ to a
point $P$.
\end{case}

The morphism $\psi$ is the elementary contraction of a $K_Z$ (and $E$) negative
extremal ray. The divisor $-K_Z-E$ is Cartier, and its restriction to
$E$, $-K_{E}=(-K_Z-E)_{\vert E}$, is anti-ample. The exceptional
divisor $E$ therefore is a Gorenstein, possibly nonnormal, del Pezzo
surface.
The map $g \colon X\to Z$ is small since the anticanonical map is an
isomorphism in codimension $1$; it induces a morphism $g_{\vert
 \widetilde{E}}\colon \widetilde{E}\to E$. The morphism $g_{\vert
 \widetilde{E}}$ uniquely induces a
morphism between the normalisations of $E$ and $\widetilde{E}$:
$\tilde{g}_{\vert \widetilde{E}} \colon \widetilde{E}^{\nu} \to
E^{\nu} $. Since the map $g$ is small, $E^{\nu}$ has the same anticanonical degree
as $\widetilde{E}^{\nu}= \widetilde{E}$. Theorem~\ref{thm:8} shows
that this degree is $1,2$ or $4$. The normalisation of $E$ is either a plane or a quadric. 

By the Serre criterion, $E$ is nonnormal if and only if it is not
regular in codimension $1$. Indeed, $Z$ is Cohen Macaulay and $E$ is
Cartier: $E$ satisfies the $S_2$-condition. Any curve $C$ lying on
$\widetilde{E}$ has $-K_X \cdot C<0$ so that no $g$-exceptional curve
lies on $\widetilde{E}$ and $g_{\vert \widetilde{E}}$ is an
isomorphism outside of a finite number of points: if $E$ is not regular in codimension
$1$, neither is $\widetilde{E}$. Cutkosky's classification
(Theorem~\ref{thm:8}) shows that $\widetilde{E}$ is normal, hence $E$
is also normal and the result follows. 
\end{proof}
\begin{rem}
  \label{rem:28}
The assumption that $Z$ is a small partial $\Q$-factorialisation of a
terminal Gorenstein Fano $3$-fold is not necessary here: the proof
works when $Z$ is assumed to be a small partial $\Q$-factorialisation
of a terminal Gorenstein $3$-fold with Picard rank $1$. The surface $S$ in Case $1$ can be
replaced by a general member of $\vert nH-K_Z \vert$, where $H$ is the
linear system that determines the contraction $\psi$.
\end{rem}
\begin{lem}
 \label{lem:32}
Let $Y$ be a non $\Q$-factorial terminal Gorenstein Fano $3$-fold with
Picard rank $1$ and genus $g\geq 3$. Let $Z$ be a small partial $\Q$-factorialisation of
$Y$. Assume that there exists an extremal
divisorial contraction $\psi \colon Z \to Z_1$ that contracts a
Cartier divisor $E$ to a curve $\Gamma$.
 The
following relations hold: 
\begin{eqnarray*} 
(-K_Z)^3=(-K_Y)^3=2g-2=-K_{Z_1}^3-2((-K_{Z_1}\cdot \Gamma)+1-p_a(\Gamma)) \\
(-K_Z)^2 \cdot E=
-K_{Z_1}\cdot \Gamma+2-2p_a(\Gamma)\\
-K_{Z} \cdot E^{2}=-2+2p_a(\Gamma) \\
E^{3}=-(-K_{Z_1}\cdot \Gamma - 2 +2 p_a(\Gamma))
 \end{eqnarray*}
\end{lem}
\begin{proof}
By Lemma~\ref{lem:29}, the contraction $\psi$ is the blow up of the
 curve $\Gamma$ and $\Gamma$ is locally a complete
 intersection. The anticanonical divisors of $Z$ and $Z_1$ satisfy:
 \begin{eqnarray}
   \label{eq:18}
-K_Z=\psi^{\ast}(-K_{Z_1})-E.   
 \end{eqnarray}
The proof of Lemma~\ref{lem:29} shows that the curve
$\Gamma$ lies on a non-singular section $S_1$ of the anticanonical
linear system $\vert -K_{Z_1} \vert$, such that the proper transform
$S$ of $S_1$ on $Z$ is a non-singular section of $\vert -K_Z
\vert$. The surfaces $S$ and $S_1$ are moreover isomorphic. 
As
$\psi^{\ast} S_1 = S + E$, $E^{3}$ satisfies:
\[
E^{3}= -S\cdot E^{2}+\psi^{\ast} S_1\cdot E^{2} = -(E_{S}\cdot
E_{S})_{S}-S_1\cdot \Gamma.
\]
The curve $\Gamma$ is locally a complete intersection, and by
construction of $S$, $E_{S}\simeq \Gamma$ and $(E_{S}\cdot
E_{S})_{S}= (\Gamma\cdot \Gamma)_{S_1}$. As $S_1$ is Cartier,
$K_{S_1}= (K_{Z_1}+S_1)_{\vert S_1}$ and since $S_1$ is non-singular,
the adjunction formula for $\Gamma$ gives
$K_{\Gamma}= (K_{S_1}+\Gamma)_{\vert\Gamma}$. No correction term is
needed at any point. In particular, 
\[
 (\Gamma\cdot \Gamma)_{S_1}= -K_{S_1}\cdot \Gamma +K_{\Gamma}=
 -K_{Z_1} \cdot \Gamma -S_1 \cdot \Gamma + 2p_a(\Gamma)-2
\]
and therefore:
\[E^{3}=
K_{Z_1}\cdot \Gamma +S_1 \cdot \Gamma +2-2p_a(\Gamma)-S_1\cdot
\Gamma=K_{Z_1}\cdot \Gamma+2-2p_a(\Gamma).\]
The relations then follow from \eqref{eq:18} and from the projection
formula:
\[ 
\psi^{\ast}(-K_{Z_1})\cdot E^{2}=-(-K_{Z_1}\cdot \Gamma).
\]
\end{proof}
\begin{lem}
 \label{lem:30}
Let $Y$ be a non $\Q$-factorial terminal quartic $3$-fold.  Assume
that $Y$ does not contain a plane $\PS^2$ with
$\vert-K_Y\vert_{\PS^2}= \mathcal{O}_{\PS^2}(1)$. Let $Z$ be a small
partial $\Q$-factorialisation of $Y$ with $\rho(Z/Y)=1$. Assume that
there exists an extremal
divisorial contraction $\psi \colon Z \to Z_1$ that contracts a
Cartier divisor $E$ to a curve $\Gamma$. Denote by $i(Z_1)$ the Fano
index of $Z_1$.
The following bounds on the degree and arithmetic genus of $\Gamma$ hold: 
\begin{enumerate}
\item If $i(Z_1)=1$ and $Z_1$ has
  genus $g_1$, then $\deg(\Gamma)=g_1-4+p_a(\Gamma)$ and $p_a(\Gamma) \leq g_1-3$,
\item If $i(Z_1)=2$ and
  $(-K_{Z_1})^3= 8d$ for some $1 \leq d \leq 5$, then
  $2\deg(\Gamma)= 4d-3+p_a(\Gamma)$ and $p_a(\Gamma)= 2k-1$, where $1 \leq k \leq d+1$, 
\item If $Z_1$ is a possibly singular quadric in $\PS^4$, then $3
  \deg(\Gamma)= 24+p_a(\Gamma)$ and $p_a(\Gamma)= 3k$, for some $0 \leq  k \leq 7$,
\item If $Z_1= \PS^3$, then $4 \deg(\Gamma)=29+p_a(\Gamma) $
  and $p_a(\Gamma)= 3+4k$, for some $0 \leq k \leq 6$.
\end{enumerate}
\end{lem}
\begin{proof}
Lemmas~\ref{lem:29} and \ref{lem:32} show that $\psi$ is the inverse
of the blow up of the curve $\Gamma$ and that the following relations
hold:
\begin{eqnarray}
\label{eq:5}
(-K_Z)^3=(-K_Y)^3=4=-K_{Z_1}^3-2((-K_{Z_1}\cdot \Gamma)+1-p_a(\Gamma)) \\
(-K_Z)^2 \cdot E=-K_{Z_1}\cdot \Gamma+2-2p_a(\Gamma).
\end{eqnarray}  
The anticanonical linear system $\vert -K_Z \vert$ is basepoint free
and it is equal to $ \vert \psi^{\ast}(-K_{Z_1})-E\vert = \vert
\psi^{\ast}(-K_{Z_1}(-\Gamma)\vert$: the curve $\Gamma$ is a
scheme-theoretic intersection of members of $\vert -K_{Z_1} \vert$. 
In particular, 
\[  
-K_{Z_1}\cdot \Gamma=i(Z_1)\deg(\Gamma) \leq (-K_{Z_1})^3.
\]
The following relations and bounds on the degree of $\Gamma$
therefore hold:
\begin{enumerate}
\item If $i(Z_1)=1$ and if $Z_1$ has
  genus $g_1$, with $4 \leq g_1 \leq 10$ or $g_1=12$, then 
  \[
 \deg(\Gamma)=g_1-4+p_a(\Gamma)\quad \mbox{and}\quad
   \deg(\Gamma) \leq 2g_1-2;
  \]
\item If $i(Z_1)=2$ and if
  $(-K_{Z_1})^3= 8d$ for some $1 \leq d \leq 5$, then
  \[
  2\deg(\Gamma)= 4d-3+p_a(\Gamma)\quad \mbox{and}\quad
\deg(\Gamma) \leq 4d;  
  \]
\item If $i(Z_1))=3$, $Z_1$ is a
  possibly singular quadric in $\PS^4$, then 
  \[
3 \deg(\Gamma)= 24+p_a(\Gamma)\quad \mbox{and}\quad
\deg(\Gamma)\leq 18;    
  \]
\item If $Z_1= \PS^3$, then 
  \[
4 \deg(\Gamma)=29+p_a(\Gamma)\quad \mbox{and}\quad
\deg(\Gamma) \leq 16.    
  \]
\end{enumerate}
Now that the degree of $\Gamma$ has been bounded, the
genus of $\Gamma$ can be bounded using \eqref{eq:5}.

By assumption, $-K_Z$ is nef and big
and induces a small map, therefore $(-K_{Z})^2\cdot E>0$.   
\end{proof}

\begin{rem}
\label{rem:35}
The bound on the genus of $\Gamma$ obtained using \eqref{eq:5} is
sharper than the Castelnuovo bound when $\Gamma$ is a non-singular curve. 
\end{rem}

\subsection{A generalised Takeuchi construction}
\label{sec:gener-take-constr}
 
Takeuchi formulates numerical constraints associated to contractions
of extremal rays \cite{MR1022045}. This approach simplifies
considerably the methods of birational classification using projection
of varieties from points or lines, which lie at the heart
of Iskovskikh's classification of Fano $3$-folds \cite{MR463151,MR503430}. 
This application of the theory of extremal rays makes for a unified
treatment of projections of Fano $3$-folds from any centre.  
Takeuchi's work focuses on non-singular weak Fano $3$-folds with
Picard rank $2$ that are small $\Q$-factorialisations of terminal Fano $3$-folds
with Picard rank $1$ and defect $1$. Given $Z$, a non-singular weak
Fano $3$-fold with Picard rank $2$, there is an elementary Sarkisov
link on $Z$ involving two extremal contractions that are not
isomorphisms in codimension $1$. The numerical constraints associated
to each type of extremal contraction yield systems of Diophantine equations, whose
solutions correspond to the only possible Sarkisov elementary links. 
I generalise Takeuchi's construction, in order to classify non
$\Q$-factorial terminal quartic $3$-folds with arbitrary defect. 

\begin{setup}[A generalised Takeuchi construction]

Let $Y$ be a non $\Q$-factorial terminal Gorenstein Fano $3$-fold with
Picard rank $1$ and Fano index $1$. 
Assume that $Y$ contains neither a plane with $
{-K_Y}_{\vert \PS^2}=\mathcal{O}_{\PS^2}(1)$ nor an irreducible
reduced quadric $Q$ with ${-K_Y}_{\vert Q}= \mathcal{O}_{Q}(1)$. 

Let $f \colon X \to Y$ be a small
$\Q$-factorialisation of $Y$. The $3$-fold $X$ is
a weak$\ast$ Fano $3$-fold and has Picard rank
$\rho(X)>1$. Lemma~\ref{lem:1} shows that there is an extremal ray $R
\in NE(X)$ and that $R$ can be contracted. A Minimal
Model Program (MMP) can be run on $X$ (Section~\ref{sec:categ-weak-fano}). 

Let $\phi$ be the first
extremal contraction of the MMP that is not an
isomorphism in codimension $1$. Assume that $\phi$ is divisorial; $\phi$
necessarily contracts a divisor $E$ to a curve, since $Y$ is
of index $1$ and contains neither a plane nor a quadric. 
Taking a different small $\Q$-factorialisation $X$
of $Y$ if necessary, we may assume that the
first extremal contraction of the MMP on $X$ is $\phi\colon
X \to X_1$. 

Denote by $\overline{E} \subset Y$ the image of $E$ by the
anticanonical map. Recall from the proof of Lemma~\ref{lem:12}
that $\overline{E}$ is a Weil non $\Q$-Cartier divisor on $Y$.
Let $Z$ be the small partial $\Q$-factorialisation of $Y$
defined by:
\[Z=\PProj \bigoplus_{n\geq
  0}{f}_{\ast}\mathcal{O}_{X}(n\widetilde{E}).\] 
The $3$-fold $Z$ is a weak Fano and
has Picard rank $2$; denote by $h \colon Z \to Y$ its anticanonical map. 
If $E'$ is the image of
$E$ on $Z$, there is an extremal contraction $\psi \colon
Z \to Z_1$ that contracts $E'$ (Lemma~\ref{lem:12}). The
proof of Lemma~\ref{lem:32} shows that since $Y$ does not contain an
irreducible quadric, $\psi$ contracts the divisor $E'$ to a curve $\Gamma$.   

Theorem~\ref{thm:16} shows that there is a
$1$-parameter flat deformation $\mathcal{Z} \to \Delta$ of $Z$ such that
for each $t \neq 0$, $\mathcal{Z}_t$ is a generalised Fano
$3$-fold of Picard rank $2$. Proposition~\ref{pro:2} and
Theorem~\ref{thm:6} show that $\mathcal{Z}/\Delta$ induces a
$1$-parameter flat
deformation $\mathcal{Y}/\Delta$ of $Y$ and that $\mathcal{Z} \to \Delta$
factors through a morphism $H \colon \mathcal{Z} \to \mathcal{Y}$. The
morphism $H$ restricts to the anticanonical map on each fibre. For
each $t \neq 0$, $\mathcal{Y}_t$ is a defect $1$,
Picard rank $1$, terminal Gorenstein Fano $3$-fold. 

Corollary~\ref{cor:9} shows that there exists a $1$-parameter flat
deformation of $Z_1$ and a morphism
$\Psi$ that fit in a diagram:
\[
\xymatrix{\mathcal{Z}  \ar[r]^{\Psi} \ar[d] & \mathcal{Z}_{1} \ar[d]\\
\Delta  \ar@{=}[r]& \Delta
}
\]     
such that $\Psi$ restricts to $\psi$ on the central fibre, and such
that for all $t \in \Delta$, $\Psi_t$ is an E$1$-contraction.
 
Pick $t \in \Delta \smallsetminus \{0\}$; we may assume that
$\mathcal{Z}_t$ is a $\Q$-factorial weak Fano $3$-fold with small anticanonical map $h_t \colon
\mathcal{Z}_t \to \mathcal{Y}_t$. Indeed, for some $t\in \Delta$,
$\mathcal{Z}_t$ is non-singular and $\Q$-factoriality is an open
condition on the base. 
The Cone theorem~\ref{lem:1} shows that $\mathcal{Z}_t$ has
exactly $2$ extremal rays and that they may be contracted. One of the
extremal contractions contracts $E'_t$ (the divisor on $\mathcal{Z}_t$
mapped to $E'$ on the central fibre) to a curve
$\Gamma_t$. Denote this divisorial contraction by $\Psi_t$. On
$\mathcal{Z}_t$, a $2$-ray game yields a diagram of the form:
\[
\xymatrix{ \quad & \mathcal{Z}_t \ar[dl]_{\Psi_t} \ar[dr]^{h_t}
\ar@{<-->}[rr]^{\Phi_t} &\quad
& \widetilde{Z}_t \ar[dr]^{\alpha_t} \ar[dl]_{\widetilde{h}_t} &\quad\\
(\mathcal{Z}_{1})_t  &\quad & \mathcal{Y}_t & \quad &
(\widetilde{\mathcal{Z}_{1}})_t.}
\]
where 
\begin{enumerate}
\item $\mathcal{Z}_t$ and $\widetilde{\mathcal{Z}_t}$ are
  $\Q$-factorial weak Fano $3$-folds with Picard rank
  $2$; the anticanonical maps $h_t$ and $\widetilde{h_t}$ of
  $\mathcal{Z}_t$ and $\widetilde{\mathcal{Z}_t}$ are small,
\item $(\mathcal{Z}_1)_t$ and $\mathcal{Y}_t$ are terminal Gorenstein
  Fano $3$-folds with Picard rank $1$,
\item $\Phi_t$ is a composition of flops,
\item $\alpha_t$ is an extremal contraction that is not an isomorphism
  in codimension $1$,
\item $(\widetilde{\mathcal{Z}_1})_t$ is one of:
  \begin{enumerate}
  \item a terminal Gorenstein Fano $3$-fold with Picard rank $1$ if $\alpha_t$ is birational,
  \item $\PS^2$ if $\alpha$ is a conic bundle,
  \item $\PS^1$ if $\alpha$ is a del Pezzo fibration.
  \end{enumerate}
\end{enumerate}
\begin{cla}
 The elementary Sarkisov link on $\mathcal{Z}_t$, $t \neq 0$, induces
 an elementary Sarkisov link on the central fibre of $\mathcal{Z} \to \Delta$. 
\end{cla}
As in the proof of Lemma~\ref{lem:9}, we may assume that $\mathcal{Z}$
is a proper flat deformation over the spectrum of a complete
ring $S$ with residue field $k$, closed point ${0}$ and generic point
$\eta$. It is enough to show that if $\mathcal{Z}/S$ is a
$1$-parameter proper flat deformation of a generalised Fano $3$-fold $Z$ with Picard
rank $2$ such that the generic fibre is terminal Gorenstein and
$\Q$-factorial, an extremal contraction on the generic fibre
$\mathcal{Z}_{\eta}$ induces an $S$-morphism that restricts to the
contraction of an extremal ray on the central fibre.

Recall that the specialisation map induces an isomorphism of the
Neron Severi groups of $\mathcal{Z}_{\eta}$ and $Z$, and that
$\Pic(\mathcal{Z}_{\eta})\simeq \Pic(Z)\simeq \Pic(\mathcal{Z})$
because the singularities of each fibre are terminal and Gorenstein. 
As the cone $NE(Z)$ is rational polyhedral and
generated by extremal rays, the
specialisation of every extremal ray $R_{\eta}$ on $Z_{\eta}$, 
belongs to the class of an extremal ray $R_0$ on $Z$, and $R_0$ may be
contracted. If $R_0$ is a flopping contraction, Theorem~\ref{thm:6}
ensures that the contraction of $R_{\eta}$ also is a flopping
contraction. If $R_0$ is a $K$-negative extremal ray, so is
$R_{\eta}$ (Lemma~\ref{lem:9}). If $cont_{R_{\eta}}$ is divisorial, it
has Cartier exceptional divisor $E_{\eta}$. The specialisation
$\Red(E_{\eta})$  is ($\Q$-)Cartier and is the exceptional divisor of
$cont_{R_0}$. The contraction of $R_{\eta}$ on $\mathcal{Z}_{\eta}$
induces a projective $S$-morphism $\mathcal{Z} \to \mathcal{Z'}$ that
restricts to the contraction of $R_{\eta}$ and $R_0$. 

The flop $\Phi_t$ induces a deformation $\widetilde{\mathcal{Z}}/
\Delta$ and a morphism $\xymatrix{\Phi \colon \mathcal{Z}
  \ar@{-->}[r]& \widetilde{\mathcal{Z}}}$
that restricts to a flop on each fibre; in particular each fibre of $\widetilde{\mathcal{Z}}/
\Delta$ is a terminal Gorenstein generalised Fano $3$-fold with Picard
rank $2$. The contraction $\alpha_t$ induces a
deformation $\widetilde{\mathcal{Z}_{1}}/\Delta$
and a morphism $\alpha \colon \widetilde{\mathcal{Z}} \to
\widetilde{\mathcal{Z}_{1}}$ that restricts to an extremal
contraction. The contraction $\alpha$ has Cartier exceptional divisor
if it is divisorial.
 
I denote by $\Phi$ and $\alpha$ the restriction to the central fibre of
$\Phi$ and $\alpha$. I presume this will not lead to any confusion. On the
central fibre, there is a diagram:  
\begin{eqnarray}
\label{eq:41}
\xymatrix{ \quad & Z \ar[dl]_{\psi} \ar[dr]^{h}
\ar@{<-->}[rr]^{\Phi} &\quad
& \widetilde{Z} \ar[dr]^{\alpha} \ar[dl]_{\tilde{h}} &\quad\\
Z_1  &\quad & Y & \quad & \widetilde{Z_1}}.
\end{eqnarray}
where
\begin{enumerate}
\item $Z$ and $\widetilde{Z}$ are terminal Gorenstein weak Fano $3$-folds with Picard rank
  $2$; the anticanonical maps $h$ and $\widetilde{h}$ of $Z$ and $\widetilde{Z}$ are small,
\item $Z_1$ and $Y$ are terminal Gorenstein
  Fano $3$-folds with Picard rank $1$,
\item $\Phi$ is a composition of flops,
\item $\alpha$ is an extremal contraction that is not an isomorphism
  in codimension $1$,
\item $\widetilde{\mathcal{Z}}_1$ is one of:
  \begin{enumerate}
  \item a terminal Gorenstein Fano $3$-fold with Picard rank $1$ if $\alpha$ is birational,
  \item $\PS^2$ if $\alpha$ is a conic bundle,
  \item $\PS^1$ if $\alpha$ is a del Pezzo fibration.
  \end{enumerate}
\end{enumerate}
Fix the following notation:
\begin{description}
\item $\widetilde{E}$ is the strict transform of $E$ on $\widetilde{Z}$,
\item $e= E^{3}-\widetilde{E}^{3} $. 
\end{description}
\end{setup}
\begin{rem}
 \label{rem:4}
The purpose of this construction via deformation theory is that
if $\alpha$ is a birational contraction, its exceptional
divisor is $\Q$-Cartier. Direct constructions on $Z$ would have been
possible, but these could have involved an extremal contraction
$\alpha$ with Weil non $\Q$-Cartier exceptional divisor.
Notice also that the $3$-fold $Y$ is not
$\Q$-factorial and the map $h$ is not the identity by construction of
$Z$. In particular, the birational map $\Phi$ is not an isomorphism.  
\end{rem}
The birational map $\Phi$ is a sequence of flops, therefore:
\begin{eqnarray}
  \label{eq:17}
  (-K_Z)^2 \cdot E= (-K_{\widetilde{Z}})^2 \cdot \widetilde{E} \nonumber\\
 -K_Z \cdot E^{2}= -K_{\widetilde{Z}}\cdot \widetilde{E}^{2}\\
\widetilde{E}^{3}=E^{3}-e.\nonumber
\end{eqnarray}
\begin{lem}\cite{MR1924722}
 \label{lem:31}
In the construction \eqref{eq:41}, $e$ is a strictly positive integer.
\end{lem}
\begin{proof}
The Cartier divisor $E$ is effective and $\psi$-negative. For any
exceptional curve $\gamma$ of $\Phi$, since $\gamma$ is a flopping
curve, $E\cdot\gamma$ is strictly positive. The map $\Phi$ is an
$E$-flopping contraction and by the construction of \cite{MR986434},
$\widetilde{E}$ is Cartier:
$e$ is an integer.

Consider a common resolution of $Z$ and $\widetilde{Z}$:
\[
\xymatrix{
\quad & W\ar[dr]^{q}\ar[dl]_{p}& \quad \\
Z \ar@{<-->}[rr]^{\Phi} &\quad & \widetilde{Z}.
}
\]
The $3$-folds $Z$ and $\widetilde{Z}$ are terminal, so that
\begin{eqnarray*}
K_W= p^{\ast}(K_Z) +E_1 +F\\
K_W= q^{\ast}(K_{\widetilde{Z}}) +E_2 +G
\end{eqnarray*}
where $E_1, E_2$ are effective $p$ and $q$-exceptional divisors and $F$
(resp. $G$) is a (possibly empty) effective $p$ but not $q$ (resp. $q$ but not $p$)
exceptional divisor.

Then: \[p^{\ast}(K_Z)= q^{\ast}(K_{\widetilde{Z}}) +G + H,\] where $G$ is
effective and contains no $p$-exceptional component and $H$ is
$p$-exceptional. Since $p_{\ast}(q^{\ast}(K_{\widetilde{Z}}) +G)=
K_{Z}$, $q^{\ast}(K_{\widetilde{Z}})$ is $p$-nef, and by the standard
negativity lemma \cite[(2.5)]{MR1311348}, $H$ is effective, that is
$F=0$ and $E_2-E_1\geq 0$. 

Reversing the roles of $Z$ and $\widetilde{Z}$ shows that, as $Z$ and
$\widetilde{Z}$ are terminal, every exceptional divisor is $p$ and
$q$-exceptional.
Moreover, $p^{\ast}(K_Z)= q^{\ast}(K_{\widetilde{Z}})+E_1-E_2$ and
$q^{\ast}(K_{\widetilde{Z}})=p^{\ast}(K_Z)+E_2-E_1$, so that, by the
negativity lemma, $p^{\ast}(K_Z)= q^{\ast}(K_{\widetilde{Z}})$.

One may write
\[p_{\ast}^{-1}E= p^{\ast} E-R=
q^{\ast}(\widetilde{E})-R',\]
where $R$ and $R'$ are effective exceptional divisors for $p$ and $q$.
In particular:
\[
-p^{\ast}(E)=-q^{\ast}(\widetilde{E})+R'-R.
\]  
By the construction of the $E$-flop, $-q^{\ast}(\widetilde{E})$ is
$p$-nef. By the negativity
lemma, the divisor $R'-R$ is strictly effective because $\Phi$ is not
an isomorphism, and its push forward $p_{\ast}(R'-R)$ is effective.
Then 
\begin{eqnarray*}
\widetilde{E}^{3}= (q^{\ast}\widetilde{E}-(R'-R)
)(q^{\ast}\widetilde{E})^2=
p^{\ast}E(q^{\ast}\widetilde{E})^2\\
= p^{\ast}E(p^{\ast}E+
(R'-R))^2= E^{3}+ Ep_{\ast}(R'-R)^2 . 
\end{eqnarray*}
This concludes the proof: $-p_{\ast}(R'-R)^2$ is a non-zero
effective $1$-cycle contained in the indeterminacy locus of $\Phi$: it
has strictly positive intersection with $E$. 
\end{proof}

The Cartier divisors $-K_Z$ and $E$ are linearly
independent. The divisor $\widetilde{E}$ is a prime divisor, because
$E$ is prime and $\Phi$ is an isomorphism in codimension $1$. Let
$i(\widetilde{Z})$ be the Fano index of $\widetilde{Z}$
and $\widetilde{H}$ the uniquely determined Cartier divisor such that
$-K_{\widetilde{Z}}= i(\widetilde{Z})\widetilde{H}$. The divisors
$\widetilde{H}$ and $\widetilde{E}$ form a $\Z$-basis of $\Pic \widetilde{Z}$. 

If $\alpha$ is birational, denote by $D$ its exceptional divisor. The
divisor $D$ is Cartier, so that there exist integers $x, y$ such that:
\begin{eqnarray}
  \label{eq:21}
D=\frac{x}{i(\widetilde{Z})}(-K_{\widetilde{Z}}) - y \widetilde{E}  
\end{eqnarray}
If $\alpha$ is of fibering type, denote by $L$ the pull back of
an ample generator of $\Pic Z_1$.
The divisor $L$ is Cartier, so that there exist
integers $x, y$ such that:
\begin{eqnarray}
  \label{eq:22}
L=\frac{x}{i(\widetilde{Z})}(-K_{\widetilde{Z}}) - y \widetilde{E}.  
\end{eqnarray}
\begin{rem}
\label{rem:36}
As is noted in Remark~\ref{rem:8}, the indices of $Y$, $Z$ and
$\widetilde{Z}$ are equal. By Theorem~\ref{thm:17}, if the Fano index of
$Y$ is $4$, $Y$ is isomorphic to $\PS^3$. If the Fano index of $Y$ is $2$,
by Lemma~\ref{lem:25}
both $\alpha$ and $\psi$ are either E$2$ contractions, \'etale
conic bundles or quadric bundles. By Remark~\ref{rem:32}, if the Fano index
of $Y$ is $3$, then $\alpha$ and $\psi$ are $\PS^2$-bundles over
$\PS^1$.
\end{rem}

I now assume that $i(\widetilde{Z})=i(Z)=1$. I am mainly interested in the case of terminal quartic
$3$-folds $Y\subset \PS^4$, which have Fano index $1$; other cases can be
treated similarly.    

The morphism $\alpha$ is the contraction of an extremal ray on a
Picard rank $2$ weak Fano $3$-fold $\widetilde{Z}$. 
If the contraction $\alpha$ is divisorial, its exceptional divisor is
Cartier. Note that $\alpha$ is induced by an extremal contraction on
a weak$\ast$ Fano $3$-fold $\widetilde{X}$, which is a small
$\Q$-factorialisation of $Y$. The initial small $\Q$-factorialisation
of $Y$, $X$, and $\widetilde{X}$ are related by a sequence of flops.

The results of Lemma~\ref{lem:29} apply to $\alpha$: one can associate to
$\alpha$ numerical constraints on the intersection numbers of powers
of $D$ (resp. $L$) with $-K_{\widetilde{Z}}$. When considered
together, the constraints
associated to $\alpha$ and $\psi$ yield systems of Diophantine equations on $x,y$.

\begin{enumerate}
\item{$\alpha$ is a conic bundle.}
  \begin{cla}
    The surface $\widetilde{Z}_1$ is $\PS^2$.
  \end{cla}

As in
the proof of Lemma~\ref{lem:6}, $-K_{\widetilde{Z}_1}$, the anticanonical divisor of the
surface $\widetilde{Z}_1$, is nef. The $3$-fold
$\widetilde{Z}$ has Picard rank $2$, hence $\widetilde{Z}_1$ is
$\PS^2$. The divisor $L$
is $\alpha^{\ast}\mathcal{O}_{\PS^2}(1)$.

\begin{cla}
The integers $x$ and $y$ are positive and coprime; $y$ can only be
equal to $1$ or $2$.  
\end{cla}

The integers $x, y$ are such that 
$L\simeq x (-K_{\widetilde{Z}})-y \widetilde{E}$ and the divisor
$\widetilde{E}$ is fixed because it is the image by $\Phi$ of a fixed
divisor on $Z$. Assume that $x$ is not positive, then $\vert L\vert
\subset \vert y\widetilde{E}\vert$. Some positive multiple of $L$
defines a map to $\PS^2$, yet the linear system $\vert y \widetilde{E}
\vert$ is $0$-dimensional: this is impossible. 
Now assume that $y$ is not positive. The
linear system $\vert L \vert$ contains $\vert x
(-K_{\widetilde{Z}})\vert$, so that $L$ is big. This contradicts
$\alpha$ being of fibering type. The integers $x, y$ are
coprime because $-K_{\widetilde{Z}}$ and $\widetilde{E}$ form a
$\Z$-basis of $\Pic  \widetilde{Z}$, $L$ is prime and $L$ is not an
integer multiple of either of them.

\begin{rem}
 \label{rem:37}
This proof shows that $x, y$ are positive coprime integers when $\alpha$
is an extremal contraction of fibering type.  
\end{rem}

Denote by $l$ an effective non-singular curve that is contracted to a point by $\alpha$.
Then $-K_{\widetilde{Z}}\cdot l \leq 2$, and
$x(-K_{\widetilde{Z}}\cdot l)= L\cdot l+y\widetilde{E}\cdot
l=y\widetilde{E}\cdot l$. The divisor $\widetilde{E}$ is Cartier and
therefore $y$ divides $(-K_{\widetilde{Z}}\cdot l)x$. As $x$ and $y$
are coprime, this shows
that $y$ can only be equal to $1$ or $2$.

Let $\Delta$ be the discriminant curve of the conic bundle
$\alpha$. Recall that the curve $\Delta$ is linearly equivalent to
$-\alpha_{\ast}(-K_{\widetilde{Z}/\widetilde{Z}_1})^2$.
\[
\left \{ \begin{array}{c}
L^3=0\\
L^2 \cdot (-K_{\widetilde{Z}})=2 \\
L \cdot (-K_{\widetilde{Z}})^2 =12-\deg(\Delta) 
\end{array} \right.
\]
Recall that $g$ is the genus of $Y$, $Z$ and $\widetilde{Z}$.
These numerical constraints, together with the intersection table
\eqref{eq:17} and the numerical constraints associated to $\psi$ in
Lemma~\ref{lem:32}, yield the system of equations:
\begin{eqnarray*}
  \label{eq:23}
(2g-2)x^3
-3(-K_{Z_1}\cdot \Gamma+2-2p_a(\Gamma))x^2y
+3(2p_a(\Gamma)-2)xy^2 \nonumber \\
+(-K_{Z_1}\cdot
\Gamma-2+2p_a(\Gamma)+e)y^3=0\\
(2g-2)x^2-2(-K_{Z_1}\cdot \Gamma+2-2p_a(\Gamma))xy\nonumber\\+(2p_a(\Gamma)-2)y^2=2\\
(2g-2)x-(-K_{Z_1}\cdot \Gamma+2-2p_a(\Gamma))y= 12-\deg(\Delta)
\end{eqnarray*}
\item{$\alpha$ is a del Pezzo fibration.}
\begin{cla}
The curve $\widetilde{Z}_1$ is $\PS^1$.    
\end{cla}

If $\alpha$ is a del Pezzo fibration, by the Leray spectral sequence
and the Kawamata-Viehweg vanishing theorem,
$\widetilde{Z}_1$ is $\PS^1$. 
The divisor $L$
is $\alpha^{\ast}\mathcal{O}_{\PS^1}(1)$. Let $d$ be the degree of the generic
fibre. 
\begin{cla}
The integers $x$ and $y$ are positive and coprime; $y$ can only be
equal to $1,2$ or $3$.  
\end{cla}
As mentioned in Remark~\ref{rem:37}, $x$ and $y$ are
coprime. Denote by $l$ an effective curve of $\widetilde{Z}$ that is
mapped to a point by $\alpha$. Depending on the degree of the generic
fibre, $-K_{\widetilde{Z}}\cdot l$ can only be $1,2$ or $3$. Moreover,
$x(-K_{\widetilde{Z}}\cdot l)= L\cdot l+y\widetilde{E}\cdot
l=y\widetilde{E}\cdot l$. The divisor $\widetilde{E}$ is Cartier and
therefore $y$ divides $(-K_{\widetilde{Z}}\cdot l)x$. As $x$ and $y$
are coprime, this shows
that $y$ can only be equal to $1,2$ or $3$.
\[
\left \{
\begin{array}{c}
L^2 \cdot (-K_{\widetilde{Z}})=0\\
L^2 \cdot \widetilde{E}=0\\
L\cdot(-K_{\widetilde{Z}})^2=d
\end{array}
\right .
\]
The following sytem of equations is associated to the configuration
$(\psi, \alpha)$:
\begin{eqnarray*}
  \label{eq:24}
(2g-2)x^2-2(-K_{Z_1}\cdot
\Gamma+2-2p_a(\Gamma))xy+\nonumber\\(2p_a(\Gamma)-2)y^2=0 \nonumber
\\
(-K_{Z_1}\cdot \Gamma+2-2p_a(\Gamma))-2(2p_a(\Gamma)-2)xy \nonumber\\-(-K_{Z_1}\cdot
\Gamma-2+2p_a(\Gamma)+e)y^2=0  \nonumber\\
(2g-2)x-(-K_{Z_1}\cdot \Gamma+2-2p_a(\Gamma))y= d  \nonumber
\end{eqnarray*}
\item{$\alpha$ is a divisorial contraction.}
Recall that $Y$ contains neither a plane with $\vert -K_Y\vert_{\vert
  \PS^2}=\mathcal{O}_{\PS^2}(1)$ nor an irreducible quadric with
normal bundle $(-1)$. In particular, $Z_1$ and $\widetilde{Z}_1$ are
Gorenstein and Lemma~\ref{lem:29} shows that $\alpha$ is of type E$1$
or E$2$. The $3$-fold $Y$ has Fano index $1$ and Picard rank $1$, so it
cannot contain a plane $\PS^2$ with $-K_{Y \vert \PS^2}= \mathcal{O}_{\PS^2}(2)$. The
morphism $\alpha$ contracts a divisor $D$ to a curve $C$.

The contraction $\alpha$ is naturally induced by the contraction
$\tilde{\alpha}$ of an
extremal ray on a weak$\ast$ Fano $3$-fold $\widetilde{X}$ that is a
small $\Q$-factorialisation of $\widetilde{Z}$. The $3$-folds $X$ and
$\widetilde{X}$ are related by a sequence of flops. The exceptional
divisor $D$ of the contraction $\alpha$ is by construction Cartier on $\widetilde{Z}$
and 
\begin{equation}
  \label{eq:19}
-K_{\widetilde{Z}}= \alpha^{\ast}(-K_{\widetilde{Z}_1})-D.  
\end{equation}
The divisor $\alpha(\widetilde{E})$ is Cartier at the centre of
$\alpha$, except possibly at finitely many points. The formulae
\eqref{eq:21} and \eqref{eq:19} show that
$y$ divides $x+1$. As above, I assume that
the indices of $Y, Z$ and $\widetilde{Z}$ are equal to $1$. 
Define $k$ by $x+1=yk$.

Note that the integer $y$ is equal to the Fano index of $\widetilde{Z}_1$.

Let $C$ be the centre of the contraction $\alpha$. The curve $C$ is
locally a complete intersection and has planar singularities.  
Lemma~\ref{lem:29} shows that the Cartier divisor $D$
satisfies the following equations. 
\[
\left \{
\begin{array}{c}
(-K_{\widetilde{Z}}+D)^3=(-K_{\widetilde{Z}}+D)^2 (-K_{\widetilde{Z}})=(-K_{\widetilde{Z}_1})^3\\
(-K_{\widetilde{Z}}+D)^2 D=0\\
(-K_{\widetilde{Z}}+D)D(-K_{\widetilde{Z}})=-K_{\widetilde{Z}_1}\cdot C=i(\widetilde{Z}_1) \deg(C)\\
(-K_{\widetilde{Z}})D^2=2p_a(C)-2
\end{array}
\right .
\]
The following system of equations can be associated to the
configuration ($\psi$, $\alpha$).
\begin{eqnarray*}
  \label{eq:25}
y^2[(2g-2)k^2-2(-K_{Z_1}\cdot
\Gamma+2-2p_a(\Gamma))k+2p_a(\Gamma)-2]\nonumber
\\=-K_{\widetilde{Z_1}}^3  \nonumber \\  
(2g-2)k^2(yk-1)+(-K_{Z_1}\cdot
\Gamma+2-2p_a(\Gamma))(2k-3k^2y)\nonumber \\+ (2p_a(\Gamma)-2)(3ky-1)+(-K_{Z_1}\cdot
\Gamma-2+2p_a(\Gamma)+e)y=0  \nonumber\\
(2g-2)k(yk-1)-(-K_{Z_1}\cdot
\Gamma+2-2p_a(\Gamma))(2yk-1)\nonumber \\+ (2p_a(\Gamma)-2)y= 
\frac{i(\widetilde{Z}_1)}{y} \deg(C)  \nonumber\\
(2g-2)(yk-1)^2-2(-K_{Z_1}\cdot
\Gamma+2-2p_a(\Gamma))y(yk-1)\nonumber \\
+(2p_a(\Gamma)-2)y^2=2p_a(C)-2  \nonumber
\end{eqnarray*}

\end{enumerate}
These systems of Diophantine equations have very few solutions, once a
value is chosen for the genus of $Y$, $Z$ and
$\widetilde{Z}$. The solutions of such systems for a given value of
the genus of $Y$ exhibit all possible Sarkisov links with midpoint along $Y$.

\subsection{Geometric Motivation of non $\Q$-factoriality}
\label{sec:geom-motiv-non}
 
Let $Y_4^3 \subset \PS^4$ be a terminal Gorenstein non
$\Q$-factorial quartic $3$-fold. There exists a Weil non-Cartier
divisor on $Y$. On the one hand, well-known examples of non $\Q$-factorial quartic
$3$-folds contain planes or quadrics. On the other hand, a very
general determinantal quartic hypersurface is known to be non
$\Q$-factorial but to contain neither a plane nor a quadric. Yet, it
does contain a Bordigo surface of degree $6$. I show that $Y$ has to contain
some surface of relatively low degree. In other words, the degree of
the surface lying on $Y$ that breaks $\Q$-factoriality cannot be
arbitrarily large. 
 
Assume that $Y$ does not contain a quadric or a plane. Let $X$ be a small
$\Q$-factorialisation of $X$. The $3$-fold $X$ is a weak$\ast$ Fano
$3$-fold on which a Minimal Model Program can be run.
\[
\xymatrix{
X=X_0 \ar[r]^{\phi_0} \ar[d]& X_1 \ar[r]^{\phi_1}\ar[d] & \cdots
& X_{n-1}\ar[r]^{\phi_n} \ar[d]& X_n \ar[d] \\
Y_0 & Y_1 & \cdots & Y_{n-1} & Y_n
}
\]
At each step, $Y_i$, the anticanonical model of $X_i$, has Picard rank
$1$. Assume that the first extremal contraction $\phi \colon X \to
X_1$ is not an isomorphism in codimension $1$. The small
$\Q$-factorialisation $X \to Y$ can always be
chosen for this to be the case.
Assume that $\phi$ is divisorial. As mentioned above, $\phi$ can only
be of type E$1$ because $Y$ contains neither a
plane nor a quadric and $Y$ has Fano index $1$.  

Let $E$ be the
exceptional divisor of $\phi$ and let $\overline{E}$ be the image of
$E$ by the anticanonical map. Denote by $Z$ the Picard rank $2$
symbolic blow up of $Y$ along $\overline{E}$. Recall from
Lemma~\ref{lem:12} that there exists an extremal contraction $\psi$ that
makes the diagram
\[
\xymatrix{
X \ar[d] \ar[r]^{\phi} & X_{1} \ar[d]\\
Z \ar[d] \ar[r]^{\psi} & Z_1\\
Y
}
\]
commutative.
The contraction $\psi$ is studied in Lemmas~\ref{lem:29}, \ref{lem:32}
and \ref{lem:30}. 
Section~\ref{sec:gener-take-constr} shows that there is a
diagram of the form:
\[
\xymatrix{ \quad & Z \ar[dl]_{\psi} \ar[dr]^{h}
\ar@{<-->}[rr]^{\Phi} &\quad
& \widetilde{Z} \ar[dr]^{\alpha} \ar[dl]_{\tilde{h}} &\quad\\
Z_1  &\quad & Y & \quad & \widetilde{Z_1}},
\]
where $\alpha$ is an extremal contraction that is not an isomorphism
in codimension $1$. If $\alpha$ is divisorial, moreover, its
exceptional divisor is Cartier. Each configuration of type $(\psi,
\alpha)$ correspond to a solution of a system of Diophantine
equations determined by the types of $\psi$ and $\alpha$. The configurations listed
in the table on page \pageref{table} are the only solutions of these systems.

\begin{rem}
I have not listed solutions with $e$ negative (Lemma~\ref{lem:31}).  
\end{rem}
\begin{rem}
We can further refine the list of solutions by eliminating the
solutions such that $h^{2,1}(Z)\neq
h^{2,1}(\widetilde{Z})$. The Hodge numbers themselves are not known,
but as $Z$ and $\widetilde{Z}$ have equal defect and the same
analytic type of singularities, from \eqref{eq:8}, it is sufficient to
determine that $h^{2,1}(\mathcal{Z}_t)=
h^{2,1}(\widetilde{\mathcal{Z}}_t)$ for the elementary
Sarkisov link on the non-singular fibre. On the non-singular fibre,
$h^{2,1}(\mathcal{Z}_t)= h^{2,1}((\mathcal{Z}_1)_t)+ p_a(\Gamma_t)$
(resp.
$h^{2,1}(\widetilde{\mathcal{Z}}_t)=h^{2,1}((\widetilde{\mathcal{Z}}_1)_t)+
p_a(C_t))$).
The terminal Gorenstein Fano $3$-fold $(\mathcal{Z}_1)_t$
(resp. $(\widetilde{\mathcal{Z}}_1)_t$) is a
flat degeneration of a non-singular Fano $3$-fold $Z_1'$ (resp. $\widetilde{Z}_1'$) with
Picard rank $1$ and the same genus. The Hodge numbers $h^{2,1}((\mathcal{Z}_1)_t)$
and $h^{2,1}((\widetilde{\mathcal{Z}}_1)_t)$ might not be easily
computable, but they are bounded above by $h^{2,1}(Z_1')$ and
$h^{2,1}(Z_1')$ respectively. This is sufficient to rule
out cases $16,25$ and $32$ in the table on page \pageref{table}.
\end{rem}

\begin{nt}
 In the table on page \pageref{table}, I write $X_{2g-2}\subset \PS^{g+1}$ for Fano $3$-folds
 of Fano index $1$, $V_d$ for Fano $3$-folds of Fano index $2$ and $Q$
 for the unique terminal Gorenstein Fano $3$-fold of Fano index $1$. 

The surface $F$ is defined as the image by the anticanonical map of
$Z$ of the exceptional divisor of $\psi$. More precisely, the
exceptional divisor $E$ is mapped by the anticanonical map to a
surface $F$ of degree at most $-K_Z^2\cdot E$. Lemma~\ref{lem:32} shows that
$-K_Z^2\cdot E=-K_{Z_1}\cdot
\Gamma'+2-2p_a(\Gamma)=i(Z_1)\deg(\Gamma)+2-2p_a(\Gamma)$.
If $\alpha$
is also an E$1$ contraction with exceptional divisor $D$, $F$ is the
image of $D$ or of $E$, depending on which one has smallest degree. 
\end{nt}

\begin{thm}[Main Theorem $2$]
  \label{thm:5}
Let $Y_4^3 \subset \PS^4$ be a terminal Gorenstein quartic $3$-fold. Then one of the following holds:
\begin{enumerate} 
\item $Y$ is $\Q$-factorial.
\item $Y$ contains a plane $\PS^2$.
\item $Y$ contains an irreducible reduced quadric $Q$. 
\item  $Y$ contains an anticanonically embedded del Pezzo surface of
  degree $4$. 
\item $Y$ has a structure of Conic Bundle over $\PS^2$, $\F_0$
  or $\F_2$.
\item $Y$ contains a rational scroll $E \to C$  over a curve $C$ whose
  genus and degree appear in the table on page \pageref{table}.
\end{enumerate}  
\end{thm}

\begin{proof}
\setcounter{step}{0} 
The only thing there is to prove is that if $Y$ does not contain a
plane and if $X$, a small $\Q$-factorialisation of $Y$, admits no
divisorial contraction and does not have a structure of Conic bundle,
then $Y$ contains an anticanonically embedded del Pezzo surface of
degree $4$.

\begin{step} 
    If $Y$ has defect $1$ and if $Y$ is the midpoint of a link between two
    del Pezzo fibrations, both these fibrations have degree $4$. 
\end{step}
 
Vologodsky shows in \cite{MR1823904} that if $Y$ has defect $1$ and if
a two ray game with midpoint along $Y$ involves two del Pezzo
fibrations, then they have the same degree $d$ and $d$ is either $2$ or $4$.
Recall from Lemma~\ref{lem:35}  that a del Pezzo fibration of degree $2$ is impossible
when the anticanonical ring of $X$ is generated in degree $1$. 

\begin{step}
 If $X$ has a structure of del Pezzo fibration of degree $4$, then
 either $Y$ contains a plane, or $Y$
 contains an anticanonically embedded del Pezzo surface of degree $4$,
 and the equation of $Y$ can be written:
\[
Y=\{a_2q +b_2q'=0 \}\subset \PS^4
\]
where $a_2, b_2, q$ and $q'$ are homogeneous forms of degree $2$ on $\PS^4$.  
\end{step}

Let $F$ be a general fibre of $X \to \PS^1$. The fibre $F$ is a
non-singular del Pezzo surface of degree $4$ and $-K_F= {-K_X}_{\vert
  F}$. The linear system $\vert -K_X \vert_{\vert F}$ is a subsystem
of $\vert -K_F \vert$. The restriction of the anticanonical map to $F$
factorises as $\Phi_{\vert -K_X \vert_{\vert F}}= \nu \circ
\Phi_{\vert -K_F \vert}$ where $\nu$ is the projection from a
(possibly empty) linear subspace 
\[
\xymatrix{\PS(H^0(F, -K_F))\simeq \PS^4\ar@{-->}[r] & \PS(H^0(F,
  \vert-K_{X} \vert_{\vert F}))}
\]
associated to the inclusion of linear systems $\vert -K_X\vert_{\vert
  F} \subset \vert -K_F \vert$.
As $X$ is a weak$\ast$ Fano $3$-fold, $\vert -K_X \vert$ is basepoint
free and the image of $F$ is a surface. The morphism $\nu$ can only be
the identity or the projection either from a line not meeting $\Phi_{\vert -K_F
\vert}(F)$ or from a point not lying on $\Phi_{\vert -K_F \vert}(F)$.

Note that if $\nu$ is not the identity, as $h^0(-K_F)=h^0(-K_X)=5$, the
map $i$ appearing in the long exact sequence in cohomology 
\begin{align*}
0 \to H^{0}(X, -K_X-F)\to H^{0}(X, -K_X) \stackrel{i} \to H^{0}(F,
-K_F)\to \\ \to H^{1}(X, -K_X-F)\to 0
\end{align*}
is not surjective. In particular, $H^0(X, -K_X-F)$ is not
trivial. There is a hyperplane section of $Y$ that contains
$\Phi_{\vert -K_X \vert}(F)$. 
As this holds for the general fibre $F$, the fibration $X \to \PS^1$
is induced by a pencil of hyperplanes on $Y$.
\begin{case}
If $H^0(X, -K_X -F) \neq (0)$ for the general fibre $F$, $Y$ contains
a plane.   
\end{case}
The rational
map $\xymatrix{Y \ar@{-->}[r]& \PS^1}$ is determined by a pencil of
hyperplanes $\mathcal{H}$. Without loss of generality, we may assume
that $\mathcal{H}$ is the pencil $\mathcal{H}_{(\lambda: \mu)})= \{\lambda x_0+\mu
x_1=0\}$ for $(\lambda: \mu)\in \PS^1$. The anticanonical map $X \to Y$ is small and $X\to \PS^1$ is
a del Pezzo fibration; the map $X \to Y$ is a resolution of the base
locus of the pencil $\mathcal{H}$
on $Y$. As the map $X\to Y$ is small, the pencil $\mathcal{H}$ has a
base component on $Y$; This is only possible if $Y$ contains the plane
$\Pi= \{x_0=x_1=0\}= \Bs\mathcal{H}$ itself. 

\begin{case}
Assume that $H^0(X, -K_X-F)$ is trivial, then $Y$ contains an anticanonically embedded
non-singular del Pezzo surface $S$ of degree $4$, that is the intersection
of two quadric hypersurfaces in $\PS^4$. 
\end{case}
The equation of $S$ is $\{
q(x_0,\ldots,x_4)=q'(x_0,\ldots,x_4)=0\}$, where $q$ and $q'$ are
homogeneous forms of degree $2$. The equation of $Y$ is:
\[
Y=\{a_2q +b_2q'=0 \}\subset \PS^4
\] with $a_2$ and $b_2$ homogeneous forms of degree $2$. 
Geometrically, the structures of del Pezzo fibrations on small
$\Q$-factorialisations of $Y$ arise as the maps induced by
pencils of quadrics (eg $\mathcal{L}=\{ q,q'\}$ amd $\mathcal{M}=\{
a_2,b_2\}$) after blowing up their base locus on $Y$, which are
anticanonically embedded del Pezzo surfaces of degree $4$.  
Considering for example the unprojection of $Y$ with variables of
weight $0$:
\begin{align*}
t=\frac{q}{q'}=-\frac{b_2}{a_2}\\
t'=\frac{q}{b_2}=-\frac{q'}{a_2}  
\end{align*}
and $X$ (resp. $X'$) the blow-up of $X$ along $S$ (resp. along $S'=
\{q=b_2=0\}$), there is a diagram
\[
\xymatrix{
\quad & X \ar[dl] \ar[dr] \ar@{-->}[rr] & \quad & X'\ar[dl] \ar[dr] \\
\PS^1 & \quad & Y & \quad & \PS^1  
}
\]
 The $3$-fold $X$ (resp. $X'$) lies on $Q \times \PS^1$ (resp. $Q'
 \times \PS^1$) for $Q \subset
 \PS^4$ (resp. $Q'$) a quadric that is the proper transform of $\{a_2=0 \}$ under
 the blow-up of $\PS^4$ along $S$ (resp. $S'$). 
The $3$-fold $X$
(resp. $X'$) is the
section of a linear system $\vert 2M +2F \vert$ on $ Q \times \PS^1$
(resp. $Q'\times\PS^1$), where $M=p_1^{\ast}\mathcal{O}_Q(1)$
(resp. $M=p_1^{\ast}(\mathcal{O}_{Q'}(1))$) and $F=p_2^{\ast}\mathcal{O}_{\PS^1}(1)$. 

These have natural structures of del Pezzo
surfaces of degree $4$ that correspond to $t=\frac{t_0}{t_1}$ and
$t'=\frac{t_0}{t_1}$. The map $\xymatrix{X \ar@{-->}[r] & X'}$ is a flop in the lines
that are the preimages of the points defined by $\{q=q'=a_2=b_2=0\}$.
\end{proof}

\newpage
\begin{tabular}[t]{|c|c|c|c|c|c|c|} \hline
\label{table}
$$ &$Z_1$ & $\widetilde{Z_1}$ & $p_a(\Gamma)$& $\deg(\Gamma)$ &
$\alpha $ & max $\deg F $\\ \hline
$1$ & $X_{22}$ & $X_{22}$ & $0$ & $8$ & E$1$,
$p_a(C)=0$, $\deg(C)= 8$ & $10$\\ \hline
$2$ &$X_{22}$ & $V_5$ & $1$ & $9$ & E$1$,
$p_a(C)=1$, $\deg(C)= 9$ & $9$\\ \hline
$3$ & $X_{22}$ & $X_{22}$ & $2$ & $10$ & E$1$,
$p_a(C)=2$, $\deg(C)= 10$ & $8$ \\ \hline
$4$ & $X_{22}$ & $ \PS^{2}$ & $2$ & $10$ & Conic
Bundle,$\deg \Delta= 4$ & $8$ \\ \hline
$5$ & $X_{22}$ & $X_{12}$ & $3$ & $11$ & E$1$,
$p_a(C)=0$, $\deg(C)=3 $ & $5$ \\ \hline
$6$ & $X_{18}$ & $X_{18}$ & $0$ & $6$ & E$1$,
$p_a(C)=0$, $\deg(C)= 6$ & $8$ \\ \hline
$7$ & $X_{18}$ & $V_4$ & $1$ & $7$ & E$1$,
$p_a(C)=1$, $\deg(C)= 7$ & $7$ \\ \hline
$8$ & $X_{18}$ & $X_{18}$ & $2$ & $8$ & E$1$,
$p_a(C)=2$, $\deg(C)= 8$ & $6$ \\ \hline
$9$ & $X_{18}$ & $\PS^{2}$ & $2$ & $8$ &
Conic Bundle, $\deg \Delta = 6$ & $6$ \\ \hline
$10$ & $X_{16}$ & $Q $ & $0$ & $5$ & E$1$,
$p_a(C)=3$, $\deg(C)= 9$ & $7$ \\ \hline
$11$ & $X_{16}$ & $X_{16}$ & $1$ & $6$ & E$1$,
$p_a(C)=1$, $\deg(C)= 6$ & $6$ \\ \hline
$12$ & $X_{16}$ & $\PS^{1}$ & $1$ & $6$ & Del Pezzo fibration
of degree $6$ & $6$ \\ \hline
$13$ & $X_{16}$ & $X_{8}$ & $2$ & $7$ & E$1$,
$p_a(C)=0$, $\deg(C)= 1$ & $3$ \\ \hline
$14$ & $X_{16}$ & $V_4$ & $2$ & $7$ & E$1$,
$p_a(C)=5$, $\deg(C)= 9$ & $5$ \\ \hline
$15$ & $X_{14}$ & $X_{14}$ & $0$ & $4$ & E$1$,
$p_a(C)=0$, $\deg(C)= 4$ & $6$\\ \hline
$16$ & $X_{14}$ & $Q$ & $1$ & $5$ & E$1$,
$p_a(C)=9$, $\deg(C)= 11$ & $5$\\ \hline
$17$ & $X_{14}$ & $V_3$ & $1$ & $5$ & E$1$,
$p_a(C)=1$, $\deg(C)= 5$ & $5$\\ \hline
$18$ & $X_{12}$ & $X_{22}$ & $0$ & $3$ & E$1$,
$p_a(C)=3$, $\deg(C)= 11$ & $5$ \\ \hline
$19$ & $X_{12}$ & $\PS^{3}$ & $0$ & $3$ & E$1$,
$p_a(C)=7$, $\deg(C)= 9$ & $5$ \\ \hline
$20$ & $X_{12}$ & $X_{12}$ & $1$ & $4$ & E$1$,
$p_a(C)=1$, $\deg(C)= 4$ & $4$ \\ \hline
$21$ & $X_{12}$ & $\PS^1$ & $1$ & $4$ & Del Pezzo fibration
of degree $4$ & $4$\\ \hline
$22$ & $X_{10}$ & $X_{10}$ & $0$ & $2$ & E$1$,
$p_a(C)=0$, $\deg(C)= 2$ & $4$\\ \hline
$23$ & $X_{10}$ & $V_5$ & $0$ & $2$ & E$1$,
$p_a(C)=7$, $\deg(C)= 12$ & $4$\\ \hline
$24$ & $X_{10}$ &  $V_2$ & $1$ & $3$ & E$1$,
$p_a(C)=1$, $\deg(C)= 3$ & $3$ \\ \hline
$25$ & $X_{10}$ &  $\PS^3$ & $1$ & $3$ & E$1$,
$p_a(C)=15$, $\deg(C)= 11$ & $3$ \\ \hline
$26$ & $X_{8}$ & $\PS^{2}$ & $0$ & $1$ & Conic Bundle, $\deg
\Delta = 7$ & $3$ \\ \hline
$27$ & $X_{8}$ & $X_{16}$ & $0$ & $1$ & E$1$,
$p_a(C)=2$, $\deg(C)=7$ & $3$ \\ \hline
$28$ & $V_2$ & $\PS^1$ & $1$ & $3$ & Del Pezzo fibration
of degree $6$ & $6$\\ \hline
$29$ & $V_2$ & $X_{16}$ & $1$ & $3$ & E$1$, $p_a(C)=1$,
$\deg(C)=6$ & $6$\\ \hline
$30$ & $V_3$ & $\PS^{3}$ & $3$ & $12$ & E$1$, $p_a(C)=3$,
$\deg(C)=8$ & $20$\\ \hline
$31$ & $V_5$ & $\PS^1 $ & $9$ & $13$ & Del Pezzo fibration
of degree $6$ & $6$ \\ \hline
$32$ & $Q$ & $X_{22}$ & $12$ & $12$ & E$1$, $p_a(C)=0$,
$\deg(C)=8$ & $10$\\ \hline
\end{tabular}

\bibliography{bibphd}
\end{document}